\documentclass[12pt,reqno]{amsart}
\usepackage{amsmath}
\usepackage{amssymb}
\usepackage{verbatim}
\usepackage{graphicx}
\usepackage[font=footnotesize]{caption}
\usepackage{subcaption}
\usepackage{epstopdf}
\usepackage{multirow}
\usepackage{url}
\usepackage{makecell}
\usepackage{tabularx}
\usepackage{longtable}
\usepackage{graphicx}
\usepackage{tikz}
\usepackage{pgfplots}
\usepackage{xcolor}
\usepackage{pgfplots}
\usepackage{pgfplotstable}
\allowdisplaybreaks

\usepackage[toc]{appendix}
\usepackage[capitalise]{cleveref}
\usepackage{color}
\allowdisplaybreaks

\usepackage[top=0.5in,bottom=0.5in,left=0.7in,right=0.7in]{geometry}

\newtheorem{lemma}{Lemma}

\newtheorem{theorem}[lemma]{Theorem}
\newtheorem{remark}[lemma]{Remark}

\newtheorem{example}[lemma]{Example}

\newcommand{\C}{\mathbb{C}}    
\newcommand{\N}{\mathbb{N}}    
\newcommand{\NN}{\mathbb{N}_0} 
\newcommand{\R}{\mathbb{R}}    

\newcommand{\nv}{\vec{n}}
\newcommand{\bo}{\mathcal{O}}
\newcommand{\Op}{\Omega_+}
\newcommand{\Om}{\Omega_-}

\newcommand{\ka}{\textsf{k}}
\newcommand{\BB}{\Gamma}
\newcommand{\ia}{\textsf{i}}
\newcommand{\wa}{\textsf{w}}
\newcommand{\ra}{\textsf{r}}
\newcommand{\B}{\mathcal{B}}
\newcommand{\id}{\mathbf{I}_d}

\DeclareMathOperator*{\argmin}{arg\,min}

\newcommand{\be}{ \begin{equation} }
	\newcommand{\ee}{ \end{equation} }
\newcommand{\odd}{\operatorname{odd}}
\newcommand{\ind}{\Lambda}

\numberwithin{equation}{section}
\numberwithin{lemma}{section}

\begin{document}
	
	\title[Sixth Order Compact Finite Difference Method for 2D Helmholtz Equations]{Sixth Order Compact Finite Difference Method for 2D Helmholtz Equations with Singular Sources and Reduced Pollution Effect}
	
	\author{Qiwei Feng}
	\author{Bin Han}
	\address{Department of Mathematical and Statistical Sciences, University of Alberta, Edmonton, Alberta, Canada T6G 2G1.}
	\email{qfeng@ualberta.ca, bhan@ualberta.ca}
	
	\author{Michelle Michelle}
	\address{Department of Mathematics, Purdue University, West Lafayette, IN, USA 47907.}
	\email{mmichell@purdue.edu}
	
	\thanks{Research supported in part by
		Natural Sciences and Engineering Research Council (NSERC) of Canada under Grant RGPIN-2019-04276, NSERC Postdoctoral Fellowship, and Alberta Innovates and Alberta Advanced Education}
	
	
	\makeatletter \@addtoreset{equation}{section} \makeatother
	
	\begin{abstract}
	   	Due to its highly oscillating solution, the Helmholtz equation is numerically challenging to solve. To obtain a reasonable solution, a mesh size that is much smaller than the reciprocal of the wavenumber is typically required (known as the pollution effect). High order schemes are desirable, because they are better in mitigating the pollution effect. In this paper, we present a high order compact finite difference method for 2D Helmholtz equations with singular sources, which can also handle any possible combinations of boundary conditions (Dirichlet, Neumann, and impedance) on a rectangular domain. Our method achieves a sixth order consistency for a constant wavenumber, and a fifth order consistency for a piecewise constant wavenumber. To reduce the pollution effect, we propose a new pollution minimization strategy that is based on the average truncation error of plane waves. Our numerical experiments demonstrate the superiority of our proposed finite difference scheme with reduced pollution effect to several state-of-the-art finite difference schemes, particularly in the critical pre-asymptotic region where $\ka h$ is near $1$ with $\ka$ being the wavenumber and $h$ the mesh size.
	\end{abstract}
	
	
	\maketitle
	
	\pagenumbering{arabic}
	
	\section{Introduction and motivations} 
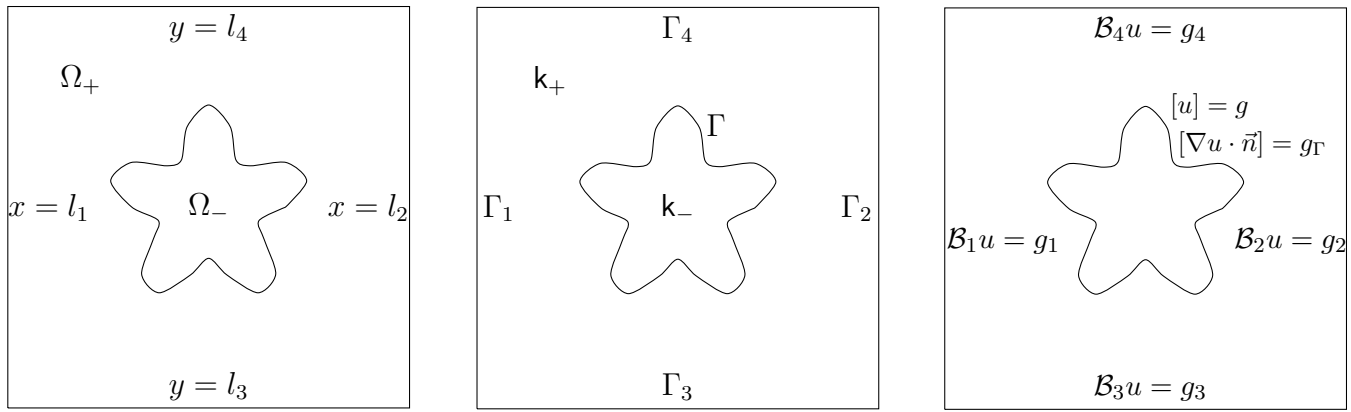
\begin{figure}[htbp]
	\centering	
	\begin{subfigure}[b]{0.3\textwidth}
			\hspace{-0.9cm}
	\begin{tikzpicture}[scale = 0.85]
		\draw[domain =0:360,smooth]
		plot({(1.2+0.4*sin(5*\x))*cos(\x)}, {(1.2+0.4*sin(5*\x))*sin(\x)});
		\draw
		(-pi, -pi) -- (-pi, pi) -- (pi, pi) -- (pi, -pi) --(-pi,-pi);
		
	    \node (A) at (0,0) {$\Omega_{-}$};
        \node (A) at (-2,2) {$\Omega_{+}$};
		
		\node (A) at (-2.5,0) {$x=l_1$};
		
		\node (A) at (2.5,0) {$x=l_2$};
		
		\node (A) at (0,-2.8) {$y=l_3$};
		
		\node (A) at (0,2.8) {$y=l_4$};
	\end{tikzpicture}
\end{subfigure}
	\begin{subfigure}[b]{0.3\textwidth}
		\begin{tikzpicture}[scale = 0.85]
			\draw[domain =0:360,smooth]
			 plot({(1.2+0.4*sin(5*\x))*cos(\x)}, {(1.2+0.4*sin(5*\x))*sin(\x)});
			\draw
			(-pi, -pi) -- (-pi, pi) -- (pi, pi) -- (pi, -pi) --(-pi,-pi);

			\node (A) at (0.6,1.3) {$\Gamma$};
			
			\node (A) at (-2.8,0) {$\BB_1$};

			\node (A) at (2.8,0) {$\BB_2$};

			\node (A) at (0,-2.8) {$\BB_3$};
		
			\node (A) at (0,2.8) {$\BB_4$};
			
		    \node (A) at (0,0) {$\ka_{-}$};
			\node (A) at (-2,2) {$\ka_{+}$};
		\end{tikzpicture}
		\end{subfigure}
	\begin{subfigure}[b]{0.3\textwidth}
			\hspace{0.4cm}
	\begin{tikzpicture}[scale = 0.85]
		\draw[domain =0:360,smooth]
		plot({(1.2+0.4*sin(5*\x))*cos(\x)}, {(1.2+0.4*sin(5*\x))*sin(\x)});
		\draw
		(-pi, -pi) -- (-pi, pi) -- (pi, pi) -- (pi, -pi) --(-pi,-pi);
		\node (A) at (1,1.6) {{\footnotesize{$\left[u\right]=g$}}};
		\node (A) at (1.6,1) { {\footnotesize{ $\left[\nabla  u \cdot \nv \right]=g_{\Gamma}$}}};
		\node (A) at (-2.3,-0.5) {{\small{ $\B_1u = g_1$}}};
		\node (A) at (2.2,-0.5) {{\small{ $\B_2u = g_2$}}};
		\node (A) at (0,-2.8) {{\small{ $\B_3u = g_3$}}};
		\node (A) at (0,2.8) {{\small{ $\B_4u = g_4$}}};
	\end{tikzpicture}
	\end{subfigure}
	\caption
	{The illustration for the model problem \eqref{Qeques2}, where  $\B_1,\ldots,\B_4  \in \{\id,\frac{\partial }{\partial \nv},\frac{\partial }{\partial \nv}- \ia \ka \id\}$, and $\ka_{\pm}$ represents the  wavenumber $\ka$ in $\Omega_{\pm}$.}
	\label{fig:model:problem}
\end{figure}

	In this paper, we study the 2D Helmholtz equation, which is a time-harmonic wave propagation model, with a singular source term along a smooth interface curve and mixed boundary conditions. The Helmholtz equation appears in many applications such as electromagnetism \cite{BS05,N01}, geophysics \cite{C16, DL19, EOV06, FG17}, ocean acoustics \cite{JKPS11}, and photonic crystals \cite{FLCG21}. Let $\Omega:=(l_1, l_2)\times(l_3, l_4)$ and $\psi$ be a smooth two-dimensional function. Consider a smooth curve $\Gamma:=\{(x,y)\in \Omega: \psi(x,y)=0\}$, which partitions $\Omega$ into two subregions:
	$\Op:=\{(x,y)\in \Omega\; :\; \psi(x,y)>0\}$ and $\Om:=\{(x,y)\in \Omega\; : \; \psi(x,y)<0\}$. The model problem (see \cref{fig:model:problem} for an illustration) is defined as follows:
	\begin{equation} \label{Qeques2}
		\begin{cases}
			\Delta u+{\ka}^2u=f &\text{in $\Omega \setminus \Gamma$},\\	
			\left[u\right]=g, \quad \left[\nabla  u \cdot \nv \right]=g_{\Gamma} &\text{on $\Gamma$},\\
			\B_1 u =g_1 \text{ on } \BB_1:=\{l_{1}\} \times (l_3,l_4), & \B_2 u =g_2 \text{ on } \BB_2:=\{l_{2}\} \times (l_3,l_4),\\
			\B_3 u =g_3 \text{ on } \BB_3:=(l_1,l_2) \times \{l_{3}\}, &
			\B_4 u =g_4 \text{ on } \BB_4:=(l_1,l_2) \times \{l_{4}\},
		\end{cases}
	\end{equation}
	where $\partial \Omega=\cup_{i=1}^4\Gamma_i$, $\ka$ is the wavenumber, $f$ is the source term, and for any point $(x_0,y_0)\in \Gamma$,
	\begin{align}
		[u](x_0,y_0) & :=\lim_{(x,y)\in \Op, (x,y) \to (x_0,y_0) }u(x,y)- \lim_{(x,y)\in \Om, (x,y) \to (x_0,y_0) }u(x,y),\label{jumpCD1}\\
		[ \nabla  u \cdot \nv](x_0,y_0) & :=  \lim_{(x,y)\in \Op, (x,y) \to (x_0,y_0) } \nabla  u(x,y) \cdot \nv- \lim_{(x,y)\in \Om, (x,y) \to (x_0,y_0) } \nabla  u(x,y) \cdot \nv, \label{jumpCD2}
	\end{align}
where $\nv$ is the unit normal vector of $\Gamma$ pointing towards $\Op$.
	In \eqref{Qeques2}, the boundary operators $\B_1,\ldots,\B_4  \in \{\id,\frac{\partial }{\partial \nv},\frac{\partial }{\partial \nv}- \ia \ka \id\}$, where
$\id$ corresponds to the Dirichlet boundary condition (sound soft boundary condition for the identical zero boundary datum), $\frac{\partial }{\partial \nv}$ corresponds to the Neumann boundary condition (sound hard boundary condition for the identical zero boundary datum), and $\frac{\partial }{\partial \nv}-\ia \ka \id$ (with $\ia$ being the imaginary unit) corresponds to the impedance boundary condition.
Moreover, the Helmholtz equation of \eqref{Qeques2} with $g=0$ is equivalent to finding the weak solution $u\in H^1(\Omega)$ of
$\Delta u+{\ka}^2u=f+g_{\Gamma}\delta_{\Gamma}$ in $\Omega$, where $\delta_{\Gamma}$ is the Dirac distribution along the interface curve $\Gamma$.
	
	The Helmholtz equation is challenging to solve numerically due to several reasons. The first is due to its highly oscillatory solution, which necessitates the use of a very small mesh size $h$ in many discretization methods. Taking a mesh size $h$ proportional to the reciprocal of the wavenumber $\ka$ is not enough to guarantee that a reasonable solution is obtained or a convergent behavior is observed. The mesh size $h$ employed in a standard discretization method often has to be much smaller than the reciprocal of the wavenumber $\ka$. In the literature, this phenomenon is referred to as the pollution effect, which has close ties to the numerical dispersion (or a phase lag). The situation is further exacerbated by the fact that the discretization of the Helmholtz equation typically yields an ill-conditioned coefficient matrix. Taken together, one typically faces an enormous ill-conditioned linear system when dealing with the Helmholtz equation, where standard iterative schemes fail to work \cite{EG11}.
	
	To gain a better insight on how the mesh size requirement is related to the wavenumber, we recall some relevant findings on the finite element method (FEM) and finite difference method (FDM). Two common ways to quantify the pollution
	effect are through the error analysis and the dispersion analysis. In the FEM literature, the former route is typically used and the analysis is applicable even for unstructured meshes. The authors in \cite{MS11} considered the interior impedance problem and discovered that the quasi-optimality in the $hp$-FEM setting can be achieved by choosing a polynomial degree $p$ and a mesh size $h$ such that $p \ge C \log(\ka)$ (for some positive $C$ independent of $\ka$, $h$, $p$) and $\ka h/p$ is small enough.
	The authors in \cite{DW15} found that for sufficiently small $\ka^{2p+1}h^{2p}$, the leading pollution term in an upper bound of the standard Sobolev $H^{1}$-norm is $\ka^{2p+1}h^{2p}$. This coincides with the numerical dispersion studied in \cite{A04,IB06}. Meanwhile, the pollution effect in the FDM setting is studied via the dispersion analysis. That is, we analyze the difference between the true and numerical wavenumbers. For order $2$ FDMs, \cite{CCFW13,CWY11} found that $\ka^3 h^2 \le C$ (for some positive $C$ independent of $\ka,h$) is required to obtain a reasonable solution. Meanwhile, for order $4$ FDM, \cite{DL19} found that $\ka^5 h^4 \le C$ (for some positive $C$ independent of $\ka,h$) is required to obtain a reasonable solution. From the previous discussion, it is clear that the grid size requirement for high order schemes is less stringent than low order ones. This is why we presently focus on the construction of a high order scheme.
	
	A lot of research effort has been invested in developing ways to cope with the enormous ill-conditioned linear system arising from a discretization of the Helmholtz equation. Various preconditioners and domain decomposition methods have been developed over the years (see the review paper \cite{GZ19} and references therein). Many variants of FEM/Galerkin/variational methods have been explored. For example, \cite{FW09,FW11} relaxed the inter-element continuity condition and imposed penalty terms on jumps across the element edges. These penalty terms can be tuned to reduce the pollution effect. A class of Trefftz methods, where the trial and test functions consist of local solutions to the underlying (homogeneous) Helmholtz equation, were considered in \cite{HMP16} and references therein. The inter-element continuity in this kind of methods typically cannot be strongly imposed. Unfortunately, the pollution effect still persists in the $h$-refinement setting of these Trefftz methods. A closely related method, called the generalized FEM or the partition of unity FEM, has been explored. It involves multiplying solutions to the homogeneous Helmholtz equation (e.g. plane waves) with elements of a chosen partition of unity, which then serve as the trial and test functions. In recent years, multiscale FEM has also become an appealing alternative to deal with the pollution effect \cite{P17}. From the perspective of FDM, one common approach is to use a scheme (preferably of high order) with minimum dispersion. We start with a stencil having a given accuracy order with some free parameters. Afterwards, we plug in the plane wave solution into the scheme and minimize the ratio between the true and numerical wavenumbers by forming an overdetermined linear system with respect to a set of discretized angles and a range of $\frac{2\pi}{\ka h}$ (i.e., the number of points per wavelength). Such a procedure has been used
	in \cite{CCFW13, CWY11, DL19, SAB14, WX18}. The resulting stencils have accuracy orders 2 in \cite{CCFW13, CWY11},
	4 in \cite{DL19}, and 6 in \cite{WX18}. The number of points used in the proposed stencil varies from 9 in \cite{CCFW13, WX18}, 13 in \cite{DL21b}, and both 17 and 25 in \cite{DL19}. Other studies on FDMs that do not explicitly consider the numerical dispersion are \cite{BTT11} (a 4th order compact FDM on polar coodinates), \cite{BTT11b} (a 4th order compact FDM), \cite{TGGT13} (a 6th order compact FDM), and \cite{ZWG19} (a 6th order FDM with non-compact stencils for corners and boundaries). The authors in \cite{FLQ11, PHL21, Z10} proposed finite difference schemes with at most fourth consistency order
	for model problems similar to ours. A characterization of the pollution effect in terms of eigenvalues was done in \cite{DV21}. The authors in \cite{CGX21} showed that the order of the numerical dispersion matches the order of the finite difference scheme for all plane wave solutions. Furthermore, by using an asymptotic analysis and modifying the wavenumber, 
	they derived an FDM stencil for vanishing source terms whose accuracy for all plane wave solutions is of order 6. It is widely accepted that the pollution effect in standard discretizations arising from FEMs and FDMs cannot be eliminated for 2D and higher dimensions \cite{BS00}. However, in 1D, we can obtain pollution free FDMs \cite{HMW21,WW14}, which are used to solve special 2D Helmholtz equations \cite{HMW21,WW17}.
	
	There are only a few papers that deal with high order finite difference discretizations of mixed boundary conditions. The authors \cite{BTT11b,LP21,TGGT13,WX18} considered a square domain with Dirichlet and Robin boundary conditions. Furthermore, a method to handle flux type boundary conditions was discussed in \cite{LP21}. The method of difference potentials was studied in \cite{BTT13,MTT12} to handle a domain with a smooth nonconforming boundary and mixed boundary conditions.
	
	From the theoretical standpoint, as long as an impedance boundary condition appears on one of the boundary sides, the solution to  \eqref{Qeques2} exists and is unique as studied in \cite{GS20}. When an impedance boundary condition is absent, we shall avoid wavenumbers that lead to nonuniqueness. The rigorous stability analysis of the problem in \eqref{Qeques2} with $g=g_{\Gamma}=0$ was also done in \cite{HM21,H07}. For the situation where $g,g_{\Gamma} \neq0$, the well-posedness, regularity, and stability were rigorously studied in \cite{MS19}.
	
	\subsection{Main contributions of this paper} We derive fifth and sixth order compact finite difference schemes with reduced pollution effect to solve \eqref{Qeques2}-\eqref{jumpCD2} given the following assumptions:
	\begin{itemize}
		\item[(A1)] The solution $u$ and the source term $f$ have uniformly continuous partial derivatives of (total) orders up to seven and six respectively in each of the subregions $\Op$ and $\Om$. However, both $u$ and $f$ may be discontinuous across the interface $\Gamma$.
		\item[(A2)] $\ka$ is piecewise constant.
		\item[(A3)] The interface curve $\Gamma$ is smooth in the sense that for each $(x^*,y^*)\in \Gamma$, there exists a local parametric equation: $\gamma: t\in \R \rightarrow \Gamma$  such that $\gamma(t^*)=(x^*,y^*)$ and $\|\gamma'(t^*)\|_{2}\ne 0$ for some $t^*\in \R$.
		\item[(A4)] The one-dimensional functions $g$ and $g_{\Gamma}$ have uniformly continuous derivatives of (total) orders up to seven and six respectively on the interface $\Gamma$.
		\item[(A5)] Each of the functions $g_1,\ldots,g_4$ has uniformly continuous derivatives of (total) order up to seven on the boundary $\partial \Omega$.
	\end{itemize}
	Our proposed compact finite difference scheme attains the maximum overall consistency order (in the context of methods relying on Taylor expansions and our sort of techniques) everywhere on the domain with the shortest stencil support for the problem \eqref{Qeques2}-\eqref{jumpCD2}.
	Similar to \cite{FHM21a,FHM21b}, our approach is based on a critical observation regarding the inter-dependence of high order derivatives of the underlying solution. When constructing a discretization stencil, we start with a general expression that allows us to recover all possible fifth and sixth consistency order finite difference schemes. The former is for piecewise constant wavenumbers, while the latter is for constant wavenumbers. When the wavenumber is constant, this general expression is critical for the next step, where we determine the remaining free parameters in the stencil by using our new pollution minimization strategy that is based on the average truncation error of plane waves. Our method differs from existing dispersion minimization methods in the literature in several ways. First, our method does not require us to compute the numerical wavenumber. Second, we use our pollution minimization procedure in the construction of all interior, boundary, and corner stencils. This is in stark contrast to the common approach in the literature, where the dispersion is minimized only in the interior stencil. The effectiveness of our pollution minimization strategy is evident from our numerical experiments. Our proposed compact finite difference scheme with reduced pollution effect outperforms several state-of-the-art finite difference schemes in the literature, particularly in the pre-asymptotic critical region where $\ka h$ is near 1. When a large wavenumber $\ka$ is present, this means that our proposed finite difference scheme is more accurate than others at a computationally feasible grid size.
	
	For the Helmholtz interface problem in \eqref{Qeques2} with a constant $\ka$, we derive a seventh consistency order compact finite difference scheme to handle nonzero jump functions at the interface. On the other hand, for a piecewise constant $\ka$, we derive a fifth consistency order compact finite difference scheme to handle nonzero jump functions at the interface. I.e., such a scheme is used at an irregular point near the interface $\Gamma$; the stencil centered at this point overlaps with both $\Op$ and $\Om$ subregions.
	
	We provide a comprehensive treatment of mixed inhomogeneous boundary conditions. In particular, our approach is capable of handling all possible combinations of Dirichlet, Neumann, and impedance boundary conditions for the 2D Helmholtz equation defined on a rectangular domain. For each corner, we explicitly provide a $4$-point stencil with at least sixth order of consistency and reduced pollution effect. For each side, we explicitly give a $6$-point stencil with at least sixth order of consistency and reduced pollution effect. To the best of our knowledge, our present work is the first paper to comprehensively study the construction of corner and boundary finite difference stencils for all possible combinations of boundary conditions on a rectangular domain. Unlike the common technique used in the literature, no ghost or artificial points are introduced in our construction.
	
	Since our proposed finite difference scheme is compact, the linear system arising from the discretization is sparse. The stencils themselves have a nice structure in that their coefficients are symmetric and take the form of polynomials of $\ka h$. Also, the coefficients in our interior stencil are simpler compared to \cite{CGX21}, as they are polynomials of degree 6, while those in \cite{CGX21} are of degree 16. Hence, the process of assembling the coefficient matrix is highly efficient. Furthermore, for a fixed constant wavenumber $\ka$ and for any given interface and boundary data, the coefficient matrix of our linear system does not change; only the vector on the right-hand side of the linear system changes. 
	
	\subsection{Organization of this paper}
 	In \cref{sec:sixord}, we explain how our proposed compact finite difference schemes are developed (fifth order for piecewise constant wavenumbers, and sixth order with reduced pollution effect for constant wavenumbers). We start our discussion by constructing the interior finite difference stencil, followed by boundary and corner stencils, and finally interface finite difference stencils. In \cref{sec:numerical}, we present several numerical experiments to demonstrate the performance of our fifth and sixth order compact finite difference schemes. In \cref{sec:proofs}, we present the proofs of theorems in \cref{sec:sixord}.
	
	\section{Stencils for sixth order compact finite difference schemes with reduced pollution effect using uniform Cartesian grids}
	\label{sec:sixord}

	We follow the same setup as in \cite{FHM21a,FHM21b}. As stated in the introduction, let $\Omega:=(l_1,l_2)\times (l_3,l_4)$. Without loss of generality, we assume $l_4-l_3=N_0 (l_2-l_1)$ for some $N_0 \in \N$. For any positive integer $N_1\in \N$, we define $N_2:=N_0 N_1$ and so the grid size is  $h:=(l_2-l_1)/N_1=(l_4-l_3)/N_2$.
	Let
	\be \label{xiyj}
	x_i:=l_1+i h, \quad i=0,\ldots,N_1, \quad \text{and} \quad y_j:=l_3+j h, \quad j=0,\ldots,N_2.
	\ee
	Our focus of this section is to develop our compact finite difference schemes on uniform Cartesian grids (fifth order for piecewise constant wavenumbers, and sixth order with reduced pollution effect for constant wavenumbers). Recall that a compact 9-point stencil centered at $(x_i,y_j)$ contains nine points $(x_i+kh, y_j+lh)$ for $k,l\in \{-1,0,1\}$. Define
	\begin{align*}
		& d_{i,j}^+:=\{(k,\ell) \; : \;
		k,\ell\in \{-1,0,1\}, \psi(x_i+kh, y_j+\ell h)> 0\}, \quad \mbox{and}\\
		& d_{i,j}^-:=\{(k,\ell) \; : \;
		k,\ell\in \{-1,0,1\}, \psi(x_i+kh, y_j+\ell h)\le 0\}.
	\end{align*}
	Thus, the interface curve $\Gamma:=\{(x,y)\in \Omega \; :\; \psi(x,y)=0\}$ splits the nine points in our compact stencil into two disjoint sets  $\{(x_{i+k}, y_{j+\ell})\; : \; (k,\ell)\in d_{i,j}^+\} \subseteq \Op$ and
	$\{(x_{i+k}, y_{j+\ell})\; : \; (k,\ell)\in d_{i,j}^-\} \subseteq \Om \cup \Gamma$. We refer to a grid/center point $(x_i,y_j)$ as
	\emph{a regular point} if  $d_{i,j}^+=\emptyset$ or $d_{i,j}^-=\emptyset$.
	The center point $(x_i,y_j)$ of a stencil is \emph{regular} if all its nine points are completely in $\Op $ (hence $d_{i,j}^-=\emptyset$) or in $\Om\cup \Gamma$ (i.e., $d_{i,j}^+=\emptyset$).
	Otherwise, the center point $(x_i,y_j)$ of a stencil is referred to as \emph{an irregular point} if both $d_{i,j}^+$ and $d_{i,j}^-$ are nonempty.
	
	Now, let us pick and fix a base point $(x_i^*,y_j^*)$ inside the open square $(x_i-h,x_i+h)\times (y_j-h,y_j+h)$, which can be written as
	\be \label{base:pt}
	x_i^*:=x_i-v_0h  \quad \mbox{and}\quad y_j^*:=y_j-w_0h  \quad \mbox{with}\quad
	-1<v_0, w_0<1.
	\ee
	We shall use the following notations:
	\be \label{ufmn}
	u^{(m,n)}:=\frac{\partial^{m+n} u}{ \partial^m x \partial^n y}(x_i^*,y_j^*)
	\quad\mbox{and}\quad
	f^{(m,n)}:=\frac{\partial^{m+n} f}{ \partial^m x \partial^n y}(x_i^*,y_j^*),
	\ee
	which are used to represent their $(m,n)$th partial derivatives at the base point $(x_i^*,y_j^*)$. For boundary functions $g_1,g_2,g_3,g_4$, we define that
	\[
	g_1^{(n)}:=\frac{d^n g_1}{dy^n}(y_j), \qquad g_2^{(n)}:=\frac{d^n g_2}{dy^n}(y_j), \qquad g_3^{(n)}:=\frac{d^n g_3}{dx^n}(x_i),\qquad g_4^{(n)}:=\frac{d^n g_4}{dx^n}(x_i).
	\]

We define $(u_h)_{i,j}$ to be the value of the numerical approximation  $u_h$ of the exact solution $u$ of the  Helmholtz interface problem \eqref{Qeques2}, at the grid point $(x_i, y_j)$.

	Define $\NN:=\N\cup\{0\}$, the set of all nonnegative integers.
Given $L\in \NN$, we define
\be \label{Sk}
\ind_{M+1}:=\{(m,n)\in \NN^2 \; : \; m+n\le M+1\}, \qquad M+1\in \NN,
\ee
	\be \label{indV12}
	\ind_{M+1}^{ 1}:=\{(m,n)\in \ind_{M+1}\; :   m=0,1\},\qquad 
\ind_{M+1}^{ 2}:=\ind_{M+1}\setminus \ind_{M+1}^{ 1}.
\ee

	For the sake of brevity, we also define
\be\label{Hel:GM1mnxy}
G_{M+1,m,n}(x,y):= \sum_{p=0}^{\lfloor \frac{M+1-m-n}{2} \rfloor}	 \sum_{\ell=p}^{p+\lfloor \frac{n}{2}\rfloor}
\frac{(-1)^{\ell}x^{m+2\ell} y^{n+2p-2\ell}}{(m+2\ell)!(n+2p-2\ell)!}  {\ell \choose p} {\ka}^{2p},
\ee
\be\label{Hel:HM1mnxy}
H_{M+1,m,n}(x,y):=\sum_{p=0}^{\lfloor \frac{M-1-m-n}{2} \rfloor} \sum_{\ell=1+p}^{1+\lfloor p+\frac{n}{2}\rfloor}  \	(-1)^{\ell-1} {\ell-1 \choose p}
{\ka}^{2p} 	\frac{x^{m+2\ell} y^{2p+n+2-2\ell}}{(m+2\ell)!(2p+n+2-2\ell)!}.
\ee

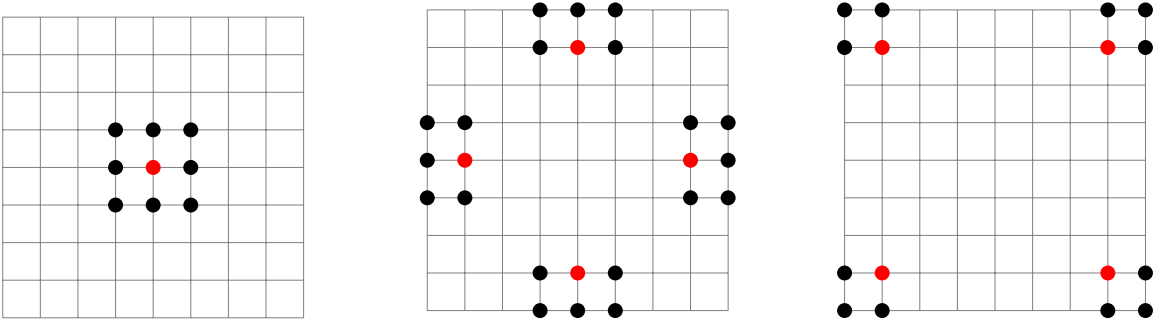
\begin{figure}[htbp]
	\hspace{0.8cm}
	\begin{subfigure}[b]{0.3\textwidth}
		\begin{tikzpicture}[scale = 2]
			\draw[help lines,step = 0.25]
			(-1,-1) grid (1,1);
			\node at (-0.25,0.25)[circle,fill,inner sep=2pt,color=black]{};
			\node at (-0.25,0)[circle,fill,inner sep=2pt,color=black]{};
			\node at (-0.25,-0.25)[circle,fill,inner sep=2pt,color=black]{};
			\node at (0,0.25)[circle,fill,inner sep=2pt,color=black]{};
			\node at (0,0)[circle,fill,inner sep=2pt,color=red]{};
			\node at (0,-0.25)[circle,fill,inner sep=2pt,color=black]{};
			\node at (0.25,0.25)[circle,fill,inner sep=2pt,color=black]{};
			\node at (0.25,0)[circle,fill,inner sep=2pt,color=black]{};
			\node at (0.25,-0.25)[circle,fill,inner sep=2pt,color=black]{};
		\end{tikzpicture}
	\end{subfigure}
	\begin{subfigure}[b]{0.3\textwidth}
		\begin{tikzpicture}[scale = 2]
			\draw[help lines,step = 0.25]
			(-1,-1) grid (1,1);
			\node at (-1,0.25)[circle,fill,inner sep=2pt,color=black]{};
			\node at (-1,0)[circle,fill,inner sep=2pt,color=black]{};
			\node at (-1,-0.25)[circle,fill,inner sep=2pt,color=black]{};
			\node at (-0.75,0.25)[circle,fill,inner sep=2pt,color=black]{};
			\node at (-0.75,0)[circle,fill,inner sep=2pt,color=red]{};
			\node at (-0.75,-0.25)[circle,fill,inner sep=2pt,color=black]{};
			\node at (0.75,0.25)[circle,fill,inner sep=2pt,color=black]{};
			\node at (0.75,0)[circle,fill,inner sep=2pt,color=red]{};
			\node at (0.75,-0.25)[circle,fill,inner sep=2pt,color=black]{};
			\node at (1,0.25)[circle,fill,inner sep=2pt,color=black]{};
			\node at (1,0)[circle,fill,inner sep=2pt,color=black]{};
			\node at (1,-0.25)[circle,fill,inner sep=2pt,color=black]{};
			\node at (0.25,-1)[circle,fill,inner sep=2pt,color=black]{};
			\node at (0,-1)[circle,fill,inner sep=2pt,color=black]{};
			\node at (-0.25,-1)[circle,fill,inner sep=2pt,color=black]{};
			\node at (0.25,-0.75)[circle,fill,inner sep=2pt,color=black]{};
			\node at (0,-0.75)[circle,fill,inner sep=2pt,color=red]{};
			\node at (-0.25,-0.75)[circle,fill,inner sep=2pt,color=black]{};
			\node at (0.25,0.75)[circle,fill,inner sep=2pt,color=black]{};
			\node at (0,0.75)[circle,fill,inner sep=2pt,color=red]{};
			\node at (-0.25,0.75)[circle,fill,inner sep=2pt,color=black]{};
			\node at (0.25,1)[circle,fill,inner sep=2pt,color=black]{};
			\node at (0,1)[circle,fill,inner sep=2pt,color=black]{};
			\node at (-0.25,1)[circle,fill,inner sep=2pt,color=black]{};
		\end{tikzpicture}
	\end{subfigure}
	\begin{subfigure}[b]{0.3\textwidth}
		\begin{tikzpicture}[scale = 2]
			\draw[help lines,step = 0.25]
			(-1,-1) grid (1,1);
			\node at (-1,-1)[circle,fill,inner sep=2pt,color=black]{};
			\node at (-1,-0.75)[circle,fill,inner sep=2pt,color=black]{};
			\node at (-0.75,-1)[circle,fill,inner sep=2pt,color=black]{};
			\node at (-0.75,-0.75)[circle,fill,inner sep=2pt,color=red]{};
			\node at (1,-1)[circle,fill,inner sep=2pt,color=black]{};
			\node at (1,-0.75)[circle,fill,inner sep=2pt,color=black]{};
			\node at (0.75,-1)[circle,fill,inner sep=2pt,color=black]{};
			\node at (0.75,-0.75)[circle,fill,inner sep=2pt,color=red]{};
			\node at (-1,1)[circle,fill,inner sep=2pt,color=black]{};
			\node at (-1,0.75)[circle,fill,inner sep=2pt,color=black]{};
			\node at (-0.75,1)[circle,fill,inner sep=2pt,color=black]{};
			\node at (-0.75,0.75)[circle,fill,inner sep=2pt,color=red]{};
			\node at (1,1)[circle,fill,inner sep=2pt,color=black]{};
			\node at (1,0.75)[circle,fill,inner sep=2pt,color=black]{};
			\node at (0.75,1)[circle,fill,inner sep=2pt,color=black]{};
			\node at (0.75,0.75)[circle,fill,inner sep=2pt,color=red]{};
		\end{tikzpicture}
	\end{subfigure}
	\caption
	{A compact 9-point scheme for the interior point (left), compact 6-point schemes for boundary side points (middle) and compact 4-point schemes for corner points (right). Red points are center points. }
	\label{HEL:fig:all:cases}
\end{figure}

In the next two subsections, we shall explicitly present our stencils having at least sixth order consistency with reduced pollution effect for interior, boundary and corner points (see three panels in \cref{HEL:fig:all:cases} for illustrations).

	\subsection{Regular points (interior)}
	\label{subsec:regular}

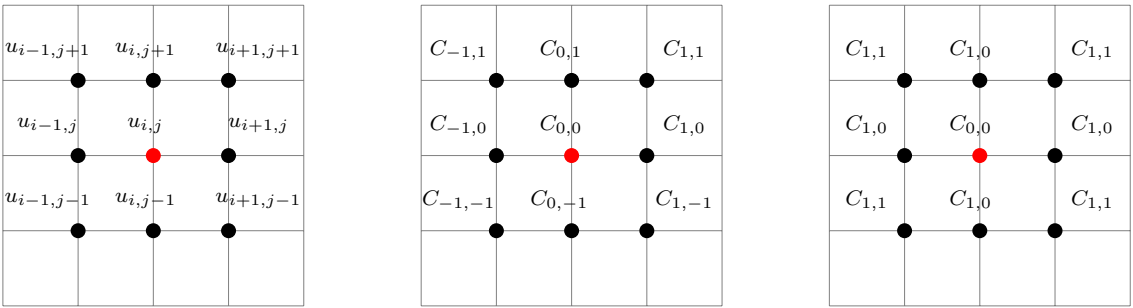
\begin{figure}[htbp]
		\hspace{0.8cm}
		\begin{subfigure}[b]{0.3\textwidth}
			\begin{tikzpicture}[scale = 2]
				\draw[help lines,step = 0.5]
				(-1,-1) grid (1,1);
				\node at (-0.5,0.5)[circle,fill,inner sep=2pt,color=black]{};
				\node at (-0.5,0)[circle,fill,inner sep=2pt,color=black]{};
				\node at (-0.5,-0.5)[circle,fill,inner sep=2pt,color=black]{};
				\node at (0,0.5)[circle,fill,inner sep=2pt,color=black]{};
				\node at (0,0)[circle,fill,inner sep=2pt,color=red]{};
				\node at (0,-0.5)[circle,fill,inner sep=2pt,color=black]{};
				\node at (0.5,0.5)[circle,fill,inner sep=2pt,color=black]{};
				\node at (0.5,0)[circle,fill,inner sep=2pt,color=black]{};
				\node at (0.5,-0.5)[circle,fill,inner sep=2pt,color=black]{};
				\node (A) at (-0.7,0.7) {{\tiny{$u_{i-1,j+1}$}}};
				\node (A) at (-0.7,0.2) {{\tiny{$u_{i-1,j}$}}};
				\node (A) at (-0.7,-0.3) {{\tiny{$u_{i-1,j-1}$}}};
				\node (A) at (-0.05,0.7) {{\tiny{$u_{i,j+1}$}}};
				\node (A) at (-0.05,0.2) {{\tiny{$u_{i,j}$}}};
				\node (A) at (-0.05,-0.3) {{\tiny{$u_{i,j-1}$}}};
				\node (A) at (0.7,0.7) {{\tiny{$u_{i+1,j+1}$}}};
				\node (A) at (0.7,0.2) {{\tiny{$u_{i+1,j}$}}};
				\node (A) at (0.7,-0.3) {{\tiny{$u_{i+1,j-1}$}}};
			\end{tikzpicture}
		\end{subfigure}	
		\begin{subfigure}[b]{0.3\textwidth}
			\begin{tikzpicture}[scale = 2]
				\draw[help lines,step = 0.5]
				(-1,-1) grid (1,1);
				\node at (-0.5,0.5)[circle,fill,inner sep=2pt,color=black]{};
				\node at (-0.5,0)[circle,fill,inner sep=2pt,color=black]{};
				\node at (-0.5,-0.5)[circle,fill,inner sep=2pt,color=black]{};
				\node at (0,0.5)[circle,fill,inner sep=2pt,color=black]{};
				\node at (0,0)[circle,fill,inner sep=2pt,color=red]{};
				\node at (0,-0.5)[circle,fill,inner sep=2pt,color=black]{};
				\node at (0.5,0.5)[circle,fill,inner sep=2pt,color=black]{};
				\node at (0.5,0)[circle,fill,inner sep=2pt,color=black]{};
				\node at (0.5,-0.5)[circle,fill,inner sep=2pt,color=black]{};
				\node (A) at (-0.75,0.7) {{\tiny{$C_{-1,1}$}}};
				\node (A) at (-0.75,0.2) {{\tiny{$C_{-1,0}$}}};
				\node (A) at (-0.75,-0.3) {{\tiny{$C_{-1,-1}$}}};
				\node (A) at (-0.07,0.7) {{\tiny{$C_{0,1}$}}};
				\node (A) at (-0.07,0.2) {{\tiny{$C_{0,0}$}}};
				\node (A) at (-0.08,-0.3) {{\tiny{$C_{0,-1}$}}};
				\node (A) at (0.75,0.7) {{\tiny{$C_{1,1}$}}};
				\node (A) at (0.75,0.2) {{\tiny{$C_{1,0}$}}};
				\node (A) at (0.75,-0.3) {{\tiny{$C_{1,-1}$}}};
			\end{tikzpicture}
		\end{subfigure}
		\begin{subfigure}[b]{0.3\textwidth}
			\begin{tikzpicture}[scale = 2]
				\draw[help lines,step = 0.5]
				(-1,-1) grid (1,1);
				\node at (-0.5,0.5)[circle,fill,inner sep=2pt,color=black]{};
				\node at (-0.5,0)[circle,fill,inner sep=2pt,color=black]{};
				\node at (-0.5,-0.5)[circle,fill,inner sep=2pt,color=black]{};
				\node at (0,0.5)[circle,fill,inner sep=2pt,color=black]{};
				\node at (0,0)[circle,fill,inner sep=2pt,color=red]{};
				\node at (0,-0.5)[circle,fill,inner sep=2pt,color=black]{};
				\node at (0.5,0.5)[circle,fill,inner sep=2pt,color=black]{};
				\node at (0.5,0)[circle,fill,inner sep=2pt,color=black]{};
				\node at (0.5,-0.5)[circle,fill,inner sep=2pt,color=black]{};
				\node (A) at (-0.75,0.7) {{\tiny{$C_{1,1}$}}};
				\node (A) at (-0.75,0.2) {{\tiny{$C_{1,0}$}}};
				\node (A) at (-0.75,-0.3) {{\tiny{$C_{1,1}$}}};
				\node (A) at (-0.06,0.7) {{\tiny{$C_{1,0}$}}};
				\node (A) at (-0.06,0.2) {{\tiny{$C_{0,0}$}}};
				\node (A) at (-0.06,-0.3) {{\tiny{$C_{1,0}$}}};
				\node (A) at (0.75,0.7) {{\tiny{$C_{1,1}$}}};
				\node (A) at (0.75,0.2) {{\tiny{$C_{1,0}$}}};
				\node (A) at (0.75,-0.3) {{\tiny{$C_{1,1}$}}};
			\end{tikzpicture}
		\end{subfigure}
		\caption
		{ The 9 values of $u$ centered at $(x_i,y_j) \notin  \partial\Omega$ (left), the general (middle) and symmetric (right) 9-point schemes centered at $(x_i,y_j) \notin  \partial\Omega$ for \cref{thm:regular:all:Ckl,thm:regular:interior}. Red points are center points. }
		\label{symmetric:Ckl:9:6}
	\end{figure}

	In this subsection, we state one of our main results on a sixth consistency order  compact finite difference scheme (with reduced pollution effect) centered at a regular point $(x_i,y_j)$ and $(x_i,y_j) \notin \partial \Omega$. We let $(x_i,y_j)$ be the base point $(x_i^*,y_j^*)$ by setting $v_0=w_0=0$ in \eqref{base:pt}. 
	%
	In the following \cref{thm:regular:all:Ckl}, we find a general expression for all possible discretization 9-point symmetric stencils centered at $(x_{i},y_j) \notin \partial \Omega$ achieving the sixth order of consistency. The proof of the following theorem is deferred to \cref{sec:proofs}. See the right panel of \cref{symmetric:Ckl:9:6}.
	 \begin{theorem}\label{thm:regular:all:Ckl}
	Let a grid point $(x_i,y_j)$ be a regular point, i.e., either $d_{i,j}^+=\emptyset$ or $d_{i,j}^-=\emptyset$ and $(x_i,y_j) \notin \partial \Omega $.
	Then the following compact 9-point symmetric stencil centered at $(x_{i},y_j)$ (see  the right panel of \cref{symmetric:Ckl:9:6})
	\be\label{stencil:all:free:Ckl}
	h^{-2}	\mathcal{L}_h u_h :=
	h^{-2}
	\sum_{k=-1}^1 \sum_{\ell=-1}^1 C_{k,\ell}(u_{h})_{i+k,j+\ell}
	=\sum_{(m,n)\in \ind_{\tilde{M}-1}} f^{(m,n)}J_{m,n}
	\ee
with
$J_{m,n}:=h^{-2} \sum_{k=-1}^1 \sum_{\ell=-1}^1 C_{k,\ell} H_{	 \tilde{M}+1,m,n}(kh, \ell h)$
has the sixth consistency order for $\Delta u+{\ka}^2u=f$ at the point $(x_i,y_j)$
if and only if the $7$th-degree polynomials of $\ka h$ of the stencil coefficients are given by
{\small{
		\be \label{CVmn:All}
		\begin{aligned}
			& C_{-1,-1}=C_{-1,1}=C_{1,-1}=C_{1,1}, \qquad C_{-1,0}=C_{0,-1}=C_{0,1}=C_{1,0},\\
			& C_{1,1} =1+(-240c_{2}+15c_{4}-120c_{6}  +480c_{10}+120c_{11} +480c_{9})\ka h +(1/15+4c_{1}+2c_{5}-8c_{7}\\
			&\qquad \quad   -2c_{8} -8c_{3})(\ka h)^2 +(-12c_{2}+c_{4}-6c_{6}+24c_{10}+6c_{11} +24c_{9})(\ka h)^3 +c_{1}(\ka h)^4 +c_{2}(\ka h)^5 \\
			&\qquad  \quad +c_{3}(\ka h)^6 + c_{9}(\ka h)^7+\bo((\ka h)^{8}),
			\\
			& C_{1,0} = 4+(-960c_{2}+60c_{4}  -480c_{6}+1920c_{10}+480c_{11}+1920c_{9})\ka h   +(1/15+16c_{1}+8c_{5}\\
			&\qquad  \quad-32c_{7}-8c_{8}  -32c_{3})(\ka h)^2   +c_{4}(\ka h)^3  +c_{5}(\ka h)^4 +c_{6}(\ka h)^5 +c_{7}(\ka h)^6 +c_{10}(\ka h)^7+\bo((\ka h)^{8}),
			\\
			& C_{0,0} = -20  +(4800c_{2}-300c_{4} +2400c_{6}-9600c_{10}-2400c_{11}-9600c_{9})\ka h +(82/15-80c_{1}\\
			&\qquad \quad -40c_{5}+160c_{7}+40c_{8}   +160c_{3})(\ka h)^2 +(-1392c_{2}+82c_{4}-696c_{6}+2784c_{10}+696c_{11} \\
			&\qquad \quad +2784c_{9})(\ka h)^3 +(-3/10 +20c_{1}+8c_{5}-48c_{7}-12c_{8}  -48c_{3})(\ka h)^4 +(92c_{2}-9c_{4}/2 +44c_{6}\\
			&\qquad \quad -192c_{10}-48c_{11}-192c_{9})(\ka h)^5 +c_{8}(\ka h)^6+c_{11}(\ka h)^7+\bo((\ka h)^{8}),
		\end{aligned}
\ee}}
as $h \rightarrow 0$ and $\tilde{M}\ge 6$, where $c_1,\ldots,c_{11} \in \R$ are free parameters, and the polynomials $H_{M+1,m,n}$ are defined in \eqref{Hel:HM1mnxy}.

%
%
		\end{theorem}

Generally, the pollution effect comes from two sources: the PDE $\Delta u+{\ka}^2u=0$ itself and the source term $f$ (i.e., highly varying or oscillating $f$).
Because the source term $f$ is known, we can reduce the pollution effect from the source term by increasing $\tilde{M}$ in \eqref{stencil:all:free:Ckl}.
To reduce the pollution effect from the PDE $\Delta u+{\ka}^2u=0$,
we minimize the average truncation error of plane waves over the free parameters $c_1,c_2,\ldots,c_{11}$ in \eqref{CVmn:All} to obtain a scheme with reduced pollution effect as follows:

	\be \label{Value:C1:C11}
	\begin{aligned}
	& c_{1} =\frac{ 303}{2^{18}}, \qquad c_{2} =-\frac{3}{2^{20}}, \qquad c_{3} =\frac{ 13}{2^{20}},\qquad c_{4} =-\frac{7}{2^{16}},\qquad c_{5} =-\frac{ 3027}{2^{20}},\\
	& c_{6} =\frac{ 5}{2^{20}}, \qquad c_{7} =-\frac{ 73}{2^{20}},\qquad c_{8} =\frac{ 4173}{2^{19}}, \qquad c_{9} = c_{10} = c_{11}=0.
	\end{aligned}
	\ee
Plugging these particular choices of $c_1,\ldots,c_{11}$ in \eqref{Value:C1:C11} into \cref{thm:regular:all:Ckl},
we have a compact 9-point scheme having
the sixth consistency order and reduced pollution effect in the following \cref{thm:regular:interior}. The proof of the following theorem is deferred to \cref{sec:proofs}.
 %
 %
	 \begin{theorem}\label{thm:regular:interior}
		Let a grid point $(x_i,y_j)$ be a regular point, i.e., either $d_{i,j}^+=\emptyset$ or $d_{i,j}^-=\emptyset$ and $(x_i,y_j) \notin \partial \Omega $.  Then the following compact 9-point symmetric stencil centered at $(x_{i},y_j)$ (see  the right panel of \cref{symmetric:Ckl:9:6})
		 \be\label{stencil:regular:interior:V}
		h^{-2}	\mathcal{L}_h u_h :=
			h^{-2}
			\sum_{k=-1}^1 \sum_{\ell=-1}^1 C_{k,\ell}(u_{h})_{i+k,j+\ell}
		=h^{-2} \sum_{(m,n)\in \ind_{6}} f^{(m,n)}\bigg( \sum\limits_{k=-1}^1 \sum\limits_{\ell=-1}^1 C_{k,\ell}
		H_{8,m,n}(kh, \ell h) \bigg),
		\ee
		achieves the sixth consistency order  for $\Delta u+{\ka}^2u=f$ at the point $(x_i,y_j)$ with reduced pollution effect, where
\begin{align}
	\nonumber
	& C_{-1,-1}=C_{-1,1}=C_{1,-1}=C_{1,1}, \quad C_{-1,0}=C_{0,-1}=C_{0,1}=C_{1,0},\\
	\nonumber
	& C_{1,1} := 1-\frac{195}{2^{17}}\ka h+\bigg(\frac{1}{15}-\frac{8709}{2^{19}}\bigg)(\ka h)^2-\frac{53}{2^{19}}(\ka h)^3+\frac{303}{2^{18}}(\ka h)^4-\frac{3}{2^{20}}(\ka h)^5+\frac{13}{2^{20}}(\ka h)^6,\\
	\label{stencil:Cv}
	& C_{1,0} :=  4-\frac{195}{2^{15}}\ka h+\bigg(\frac{1}{15}-\frac{8709}{2^{17}}\bigg)(\ka h)^2-\frac{7}{2^{16}}(\ka h)^3-\frac{3027}{2^{20}}(\ka h)^4+\frac{5}{2^{20}}(\ka h)^5-\frac{73}{2^{20}}(\ka h)^6,\\
	\nonumber
	& C_{0,0} := -20+\frac{975}{2^{15}}\ka h+\bigg(\frac{43545}{2^{17}}+\frac{82}{15}\bigg)(\ka h)^2-\frac{1061}{2^{17}}(\ka h)^3-\bigg(\frac{3}{10}+\frac{3039}{2^{15}}\bigg)(\ka h)^4 \\
	\nonumber
	&\qquad \quad +\frac{7}{2^{14}}(\ka h)^5 +\frac{4173}{2^{19}}(\ka h)^6,
\end{align}
$H_{8,m,n}(x,y)$ is defined in \eqref{Hel:HM1mnxy}, and $(x_i^*,y_j^*)=(x_i,y_j)$.
	\end{theorem}
	
	\subsection{Boundary and corner points}
	Recall that
	$ \BB_1:=\{l_{1}\} \times (l_3,l_4)$,  $\BB_2:=\{l_{2}\} \times (l_3,l_4)$,
	$ \BB_3:=(l_1,l_2) \times \{l_{3}\}$, and $\BB_4:=(l_1,l_2) \times \{l_{4}\}$. Similar to the ideas of \cref{thm:regular:all:Ckl,thm:regular:interior}, we first construct all possible compact stencils with the sixth or seventh consistency order and then minimize the average truncation error of plane waves over the free parameters of stencils to reduce the pollution effect. We discuss how to find a compact scheme centered at $(x_i,y_j) \in \partial \Omega=\cup_{i=1}^4 \Gamma_i$  in this subsection.

	\subsubsection{Boundary points}
	We first discuss in detail how the left boundary (i.e., $(x_i,y_j) \in \BB_1 := \{l_{1}\} \times (l_3,l_4)$) stencil is constructed. The stencils for the other three boundaries can afterwards be obtained by symmetry. If $\B_1 u = u =g_1$ on $\BB_1$, then the left boundary stencil can be directly obtained from \eqref{stencil:regular:interior:V}-\eqref{stencil:Cv} in \cref{thm:regular:interior} by replacing $(u_h)_{0,j-1}$, $(u_h)_{0,j}$, and $(u_h)_{0,j+1}$ with $g_1(y_{j-1})$, $g_1(y_{j})$, and $g_1(y_{j+1})$ respectively, where $y_j \in (l_3,l_4)$, and moving terms involving these known boundary values to the right-hand side of \eqref{stencil:regular:interior:V}. The other three boundary sides are dealt in a similar straightforward fashion if a Dirichlet boundary condition is present. On the other hand, the stencils for the other two boundary conditions are not trivial at all. 
	%
	%
	For the sake of presentation, we define
	\be\label{GH1234}
	\begin{split}
&	G^1_{n} :=  -h^{-1}\sum\limits_{k=0}^1 \sum\limits_{\ell=-1}^1 C_{k,\ell} G_{8,1,n}(kh, \ell h),\qquad
	H^1_{m,n} := h^{-1}\sum\limits_{k=0}^1 \sum\limits_{\ell=-1}^1 C_{k,\ell} H_{8,m,n}(kh, \ell h),\\
	&G^2_{n}: =h^{-1}\sum_{k=-1}^0 \sum_{\ell=-1}^{1} {C}_{-k,\ell} G_{8,1,n}(kh, \ell h), \qquad H^2_{m,n}: =h^{-1}\sum_{k=-1}^0 \sum_{\ell=-1}^{1} {C}_{-k,\ell} H_{8,m,n}(kh, \ell h), \\	
	& G^3_{n} := -h^{-1} \sum_{k=-1}^1 \sum_{\ell=0}^{1} {C}_{\ell,k} G_{8,1,n}( \ell h,kh),\qquad H^3_{m,n} := h^{-1}\sum_{k=-1}^1 \sum_{\ell=0}^{1} {C}_{\ell,k} H_{8,n,m}(\ell h,kh),\\
	& G^4_{n}:= h^{-1} \sum_{k=-1}^1 \sum_{\ell=-1}^0 {C}_{-\ell,k} G_{8,1,n}( \ell h, kh),\qquad H^4_{m,n} := h^{-1}\sum_{k=-1}^1 \sum_{\ell=-1}^0 {C}_{-\ell,k} H_{8,n,m}( \ell h, kh),
	\end{split}
	\ee
where	$G_{8,m,n}(x,y)$ and $H_{8,m,n}(x,y)$  are defined in \eqref{Hel:GM1mnxy} and \eqref{Hel:HM1mnxy}.

The following theorem provides the explicit $6$-point stencil of  consistency order at least six with reduced pollution effect for the left boundary operator $\mathcal{B}_1 \in \{\frac{\partial }{\partial \nv}- \ia \ka \id,\frac{\partial }{\partial \nv}\}$. The proof of the following result is deferred to \cref{sec:proofs}.
	
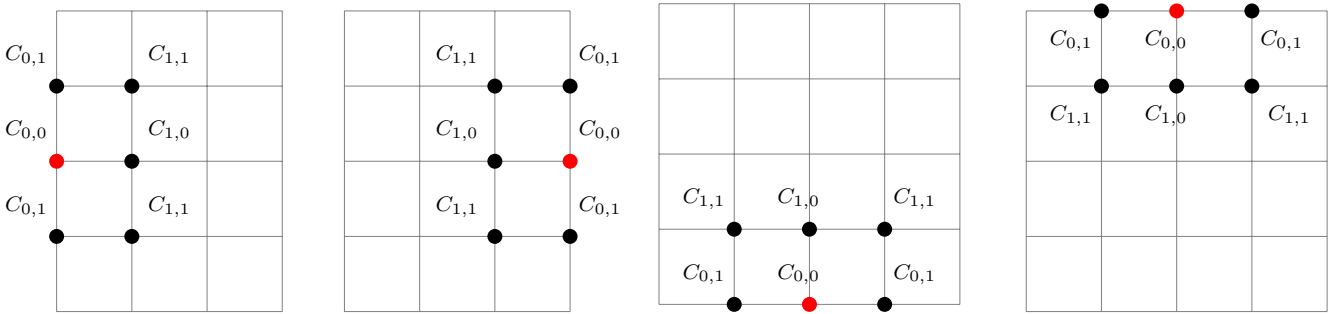
\begin{figure}[htbp]
		\begin{subfigure}[b]{0.2\textwidth}
			\hspace{-1.2cm}
			\begin{tikzpicture}[scale = 2]
				\draw[help lines,step = 0.5]
				(-1,-1) grid (0.5,1);
				\node at (-1,0.5)[circle,fill,inner sep=2pt,color=black]{};
				\node at (-1,0)[circle,fill,inner sep=2pt,color=red]{};
				\node at (-1,-0.5)[circle,fill,inner sep=2pt,color=black]{};
				\node at (-0.5,0.5)[circle,fill,inner sep=2pt,color=black]{};
				\node at (-0.5,0)[circle,fill,inner sep=2pt,color=black]{};
				\node at (-0.5,-0.5)[circle,fill,inner sep=2pt,color=black]{};
				\node (A) at (-1.2,0.7) {{\tiny{$C_{0,1}$}}};
				\node (A) at (-1.2,0.2) {{\tiny{$C_{0,0}$}}};
				\node (A) at (-1.2,-0.3) {{\tiny{$C_{0,1}$}}};
				\node (A) at (-0.25,0.7) {{\tiny{$C_{1,1}$}}};
				\node (A) at (-0.25,0.2) {{\tiny{$C_{1,0}$}}};
				\node (A) at (-0.25,-0.3) {{\tiny{$C_{1,1}$}}};
			\end{tikzpicture}
		\end{subfigure}
		\hspace{-0.3cm}
		\begin{subfigure}[b]{0.2\textwidth}
			\begin{tikzpicture}[scale = 2]
				\draw[help lines,step = 0.5]
				(-0.5,-1) grid (1,1);
				\node at (0.5,0.5)[circle,fill,inner sep=2pt,color=black]{};
				\node at (0.5,0)[circle,fill,inner sep=2pt,color=black]{};
				\node at (0.5,-0.5)[circle,fill,inner sep=2pt,color=black]{};
				\node at (1,0.5)[circle,fill,inner sep=2pt,color=black]{};
				\node at (1,0)[circle,fill,inner sep=2pt,color=red]{};
				\node at (1,-0.5)[circle,fill,inner sep=2pt,color=black]{};
				\node (A) at (0.25,0.7) {{\tiny{$C_{1,1}$}}};
				\node (A) at (0.25,0.2) {{\tiny{$C_{1,0}$}}};
				\node (A) at (0.25,-0.3) {{\tiny{$C_{1,1}$}}};
				\node (A) at (1.2,0.7) {{\tiny{$C_{0,1}$}}};
				\node (A) at (1.2,0.2) {{\tiny{$C_{0,0}$}}};
				\node (A) at (1.2,-0.3) {{\tiny{$C_{0,1}$}}};
			\end{tikzpicture}		
		\end{subfigure}
		\hspace{0.3cm}
		\begin{subfigure}[b]{0.2\textwidth}
			\begin{tikzpicture}[scale = 2]
				\draw[help lines,step = 0.5]
				(-1,-1) grid (1,1);
				\node at (-0.5,-0.5)[circle,fill,inner sep=2pt,color=black]{};
				\node at (-0.5,-1)[circle,fill,inner sep=2pt,color=black]{};
				\node at (0,-0.5)[circle,fill,inner sep=2pt,color=black]{};
				\node at (0,-1)[circle,fill,inner sep=2pt,color=red]{};
				\node at (0.5,-0.5)[circle,fill,inner sep=2pt,color=black]{};
				\node at (0.5,-1)[circle,fill,inner sep=2pt,color=black]{};
				\node (A) at (-0.7,-0.3) {{\tiny{$C_{1,1}$}}};
				\node (A) at (-0.7,-0.8) {{\tiny{$C_{0,1}$}}};	
				\node (A) at (-0.07,-0.3) {{\tiny{$C_{1,0}$}}};
				\node (A) at (-0.07,-0.8) {{\tiny{$C_{0,0}$}}};
				\node (A) at (0.7,-0.3) {{\tiny{$C_{1,1}$}}};
				\node (A) at (0.7,-0.8) {{\tiny{$C_{0,1}$}}};	
			\end{tikzpicture}
		\end{subfigure}
		\hspace{1cm}
		\begin{subfigure}[b]{0.2\textwidth}
			\begin{tikzpicture}[scale = 2]
				\draw[help lines,step = 0.5]
				(-1,-1) grid (1,1);
				\node at (-0.5,1)[circle,fill,inner sep=2pt,color=black]{};
				\node at (-0.5,0.5)[circle,fill,inner sep=2pt,color=black]{};
				\node at (0,1)[circle,fill,inner sep=2pt,color=red]{};
				\node at (0,0.5)[circle,fill,inner sep=2pt,color=black]{};
				\node at (0.5,1)[circle,fill,inner sep=2pt,color=black]{};
				\node at (0.5,0.5)[circle,fill,inner sep=2pt,color=black]{};
				\node (A) at (-0.7,0.8) {{\tiny{$C_{0,1}$}}};
				\node (A) at (-0.7,0.3) {{\tiny{$C_{1,1}$}}};
				\node (A) at (-0.07,0.8) {{\tiny{$C_{0,0}$}}};
				\node (A) at (-0.07,0.3) {{\tiny{$C_{1,0}$}}};
				\node (A) at (0.7,0.8) {{\tiny{$C_{0,1}$}}};
				\node (A) at (0.75,0.3) {{\tiny{$C_{1,1}$}}};
			\end{tikzpicture}
		\end{subfigure}
		\caption
		{ The symmetric compact 6-point scheme centered at $(x_i,y_j) =(x_0,y_j)\in  \Gamma_1$ (first), the symmetric 6-point scheme centered at $(x_i,y_j) =(x_{N_1},y_j) \in  \Gamma_2$ (second), the symmetric 6-point scheme centered at $(x_i,y_j) =(x_i,y_0) \in  \Gamma_3$ (third), and the symmetric 6-point scheme centered at $(x_i,y_j) =(x_i,y_{N_2}) \in  \Gamma_4$ (fourth) for \cref{thm:regular:Robin:1}. Red points are center points. }
		\label{symmetric:Ckl:6}
	\end{figure}

	 \begin{theorem}\label{thm:regular:Robin:1}
	 Consider the following compact 6-point stencil  centered at $(x_0,y_j)\in  \BB_1:=\{l_{1}\} \times (l_3,l_4)$ (see the first panel of \cref{symmetric:Ckl:6}) for $\mathcal{B}_1 u = g_1$ on $\BB_1$ with
		$\mathcal{B}_1 \in \{\frac{\partial }{\partial \nv}- \ia \ka \id,\frac{\partial }{\partial \nv}\}$:
		\be	 \label{stencil:regular:interior:V:Robin:1}
		h^{-1}	\mathcal{L}_h u_h  :=
		\begin{aligned}	
				h^{-1}
		\sum_{k=0}^1 \sum_{\ell=-1}^1 C_{k,\ell}(u_{h})_{k,j+\ell}
		\end{aligned} =	\sum_{(m,n)\in \ind_{6}} f^{(m,n)}H^1_{m,n}+\sum_{n=0}^{7}g_{1}^{(n)}G^1_{n},
		\ee
		where  $G^1_{n}$ and $H^1_{m,n}$ are defined in \eqref{GH1234}, $(x_i^*,y_j^*)=(x_0,y_j)$.
		\begin{itemize}
			\item[(1)] For $\B_1=\frac{\partial }{\partial \nv}- \ia \ka \id$, the coefficients in \eqref{stencil:regular:interior:V:Robin:1} are given by
	{\footnotesize{
			\be \label{stencil:CB1:R}
			\begin{aligned}
				& C_{0,-1}=C_{0,1},\qquad C_{1,-1}=C_{1,1},\\
				&C_{1,1} = 1-120\bigg(\frac{237+433\ia}{2^{17}}-\frac{2\ia}{225}\bigg)\ka h  +\bigg(\frac{99}{2^{9}}-\frac{979\ia}{2^{13}}-\frac{48}{135}\bigg)(\ka h)^2  +\frac{1017-410\ia}{2^{15}} (\ka h)^3   -\frac{112 +49\ia}{2^{15}}(\ka h)^4,
				\\
				& C_{0,1} = 2-\bigg(15 \cdot \frac{237+433\ia}{2^{13}}-\frac{29\ia}{15}\bigg)\ka h  +\bigg(3\cdot \frac{1679}{ 2^{14}}-\frac{3205\ia}{ 2^{14}}-\frac{11}{18}\bigg)(\ka h)^2  +\bigg(  \frac{2841+7271\ia}{ 2^{16}} -\frac{17\ia}{90}\bigg)(\ka h)^3\\
				& \qquad \quad
				+\frac{807+798\ia}{2^{15}}(\ka h)^4,
				\\
				& C_{1,0} =4
				-\bigg( 15\cdot \frac{237+433\ia}{2^{12}}-\frac{58\ia}{15}\bigg)\ka h+\bigg( 3\cdot \frac{1679}{ 2^{13}} -\frac{3205\ia}{ 2^{13}}-\frac{49}{45}\bigg)(\ka h)^2 +\bigg(  3\cdot \frac{631}{2^{15}} + \frac{5539\ia}{ 2^{15}} -\frac{\ia}{5} \bigg)(\ka h)^3 \\
				& \qquad \quad +\frac{87+397\ia}{2^{15}}(\ka h)^4,\\
				& C_{0,0} = -10
				+\bigg( 75 \cdot \frac{237+433\ia}{2^{13}}-\frac{58\ia}{15}\bigg)\ka h +\bigg( 3\cdot \frac{2081}{ 2^{13}} -\frac{2297\ia}{ 2^{13}}+\frac{1}{45} \bigg)(\ka h)^2 -\bigg(  5 \cdot \frac{907\ia}{2^{13}}-\frac{7\ia}{9}+\frac{3723}{ 2^{14}} \bigg)(\ka h)^3\\
				& \qquad \quad-\bigg(\frac{347+148\ia}{2^{12}}-\frac{1}{20}\bigg)(\ka h)^4.
			\end{aligned}
			\ee	}}
			Then the finite difference scheme in \eqref{stencil:regular:interior:V:Robin:1} achieves sixth consistency order  for $\B_1u=\frac{\partial u}{\partial \nv}- \ia \ka u=g_1$ at the point $(x_0,y_j) \in \BB_1$ with reduced pollution effect.
			\item[(2)] For $\B_1=\frac{\partial }{\partial \nv}$, the coefficients in \eqref{stencil:regular:interior:V:Robin:1} are given by
			\be\label{stencil:CB1:N}
\begin{aligned}
	& C_{0,-1}=C_{0,1},\qquad C_{1,-1}=C_{1,1},\qquad C_{1,1} = 1+\bigg(\frac{163}{2^{14}}+\frac{1}{15}\bigg)(\ka h)^2+\frac{99}{2^{15}}(\ka h)^4, \\
	& C_{0,1} = 2+\bigg(\frac{163}{2^{13}}+\frac{1}{30}\bigg)(\ka h)^2-\frac{35}{2^{16}}(\ka h)^4,
	\quad C_{1,0} =4+\bigg(\frac{ 163}{2^{12}}+\frac{1}{15}\bigg)(\ka h)^2-\frac{35}{2^{15}}(\ka h)^4, \\
	&C_{0,0} = -10-5\bigg( \frac{163}{2^{13}}-\frac{41}{75}\bigg)(\ka h)^2+\bigg(\frac{425}{2^{14}}-\frac{3}{20}\bigg)(\ka h)^4.
\end{aligned}
\ee
			Then the finite difference scheme in \eqref{stencil:regular:interior:V:Robin:1} achieves seventh consistency order  for $\B_1u=\frac{\partial u}{\partial \nv}=g_1$ at the point $(x_0,y_j) \in \BB_1$ with reduced pollution effect.
		\end{itemize}		
	\end{theorem}

	By symmetry, we can directly state the stencils for the other three boundary sides. Same order of consistency results as in \cref{thm:regular:Robin:1} hold. Recall the definitions of $\{C_{k,\ell}\}_{k\in\{0,1\},\ell\in\{-1,0,1\}}$ in \eqref{stencil:CB1:R}-\eqref{stencil:CB1:N}, and $G^i_{n}$, $H^{i}_{m.n}$ for $i=2,3,4$ in \eqref{GH1234}. The compact 6-point stencil for $\mathcal{B}_2 u = g_2$ on  $\BB_2$ (see the second panel of \cref{symmetric:Ckl:6}) with
	$\mathcal{B}_2 \in \{\frac{\partial }{\partial \nv}- \ia \ka \id,\frac{\partial }{\partial \nv}\}$ centered at $(x_{N_1},y_j) \in  \BB_2:=\{l_{2}\} \times (l_3,l_4)$ is
	\[
	h^{-1}	\mathcal{L}_h u_h :=h^{-1} \sum_{k=-1}^0 \sum_{\ell=-1}^{1} {C}_{-k,\ell}(u_{h})_{N_1+k,j+\ell}
	=	\sum_{(m,n)\in \ind_{6}} f^{(m,n)}H^2_{m,n}+\sum_{n=0}^{7}g_{2}^{(n)}G^2_{n},
	\]
	where $(x_i^*,y_j^*)=(x_{N_1},y_j)$. 
	The compact 6-point stencil for $\mathcal{B}_3 u = g_3$ on $ \BB_3:=(l_1,l_2) \times \{l_{3}\}$ (see the  third panel of \cref{symmetric:Ckl:6})  with $\mathcal{B}_3 \in \{\frac{\partial }{\partial \nv}- \ia \ka \id,\frac{\partial }{\partial \nv}\}$ centered at $(x_i,y_0) \in \BB_3$ is
	\[
h^{-1}	\mathcal{L}_h u_h :=
	h^{-1}\sum_{k=-1}^1 \sum_{\ell=0}^{1} {C}_{\ell,k}(u_{h})_{i+k,\ell}
	=	\sum_{(m,n)\in \ind_{6}} f^{(m,n)}H^3_{m,n}+\sum_{n=0}^{7}g_{3}^{(n)}G^3_{n},
	\]
	where $(x_i^*,y_j^*)=(x_i,y_0)$. 
	The compact 6-point stencil for $\mathcal{B}_4 u = g_4$ on $ \BB_4:=(l_1,l_2) \times \{l_{4}\}$ (see the fourth panel of \cref{symmetric:Ckl:6})  with
	$\mathcal{B}_4 \in \{\frac{\partial }{\partial \nv}- \ia \ka \id,\frac{\partial }{\partial \nv}\}$ centered at  $(x_i,y_{N_2}) \in \BB_4$ is
	\[
h^{-1}	\mathcal{L}_h u_h :=
h^{-1}	\sum_{k=-1}^1 \sum_{\ell=-1}^0 {C}_{-\ell,k}(u_{h})_{i+k,N_2+\ell}
	=	\sum_{(m,n)\in \ind_{6}} f^{(m,n)}H^4_{m,n}+\sum_{n=0}^{7}g_{4}^{(n)}G^4_{n},
	\]
	where $(x_i^*,y_j^*)=(x_i,y_{N_2})$. 
	
	\subsubsection{Corner points} \label{sssec:corner}
\begin{figure}[htbp]
	\centering	
	\resizebox{0.5\textwidth}{!}{
		\begin{tikzpicture}
			\draw
			(-pi, -pi) -- (-pi, pi) -- (pi, pi) -- (pi, -pi) --(-pi,-pi);
			\node (A) at (-5,0) {$\B_1u := \tfrac{\partial u}{\partial \nv}-\ia \ka u=g_1$};
			\node (A) at (4.4,0) {$\B_2u:=u=g_2$};
			\node (A) at (0,-3.5) {$\B_3u:=\tfrac{\partial u}{\partial \nv}=g_3$};
			\node (A) at (0,3.5) {$\B_4u := \tfrac{\partial u}{\partial \nv}-\ia \ka u=g_4$};
			\node (A) at (-2.5,0) {$x=l_1$};
			\node (A) at (2.5,0) {$x=l_2$};
			\node (A) at (0,-2.8) {$y=l_3$};
			\node (A) at (0,2.8) {$y=l_4$};
		\end{tikzpicture}
	}
	\caption
	{An illustration for the boundary configuration considered in \cref{thm:corner:1,thm:corner:2} in \cref{sssec:corner}.}
	\label{fig:config}
	\end{figure}
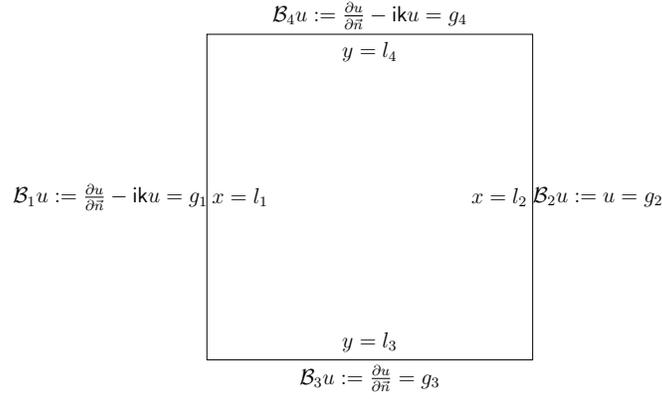

	For clarity of presentation, consider the following boundary configuration (see \cref{fig:config}).
	\begin{equation*}
	\begin{aligned}
	&\B_1u:=\tfrac{\partial u}{\partial \nv}- \ia \ka u=g_1 \;\; \text{on} \;\;  \BB_1:=\{l_{1}\} \times (l_3,l_4),
	&& \quad \B_2u:=u=g_2 \;\; \text{on} \;\;   \BB_2:=\{l_{2}\} \times (l_3,l_4), \\
	&\B_3u:=\tfrac{\partial u}{\partial \nv}=g_3 \;\; \text{on} \;\;  \BB_3:=(l_1,l_2) \times \{l_{3}\},
	&& \quad \B_4u:=\tfrac{\partial u}{\partial \nv}- \ia \ka u=g_4 \;\; \text{on} \;\; \BB_4:=(l_1,l_2) \times \{l_{4}\}.
	\end{aligned}
	\end{equation*}
	See \cref{fig:config}.

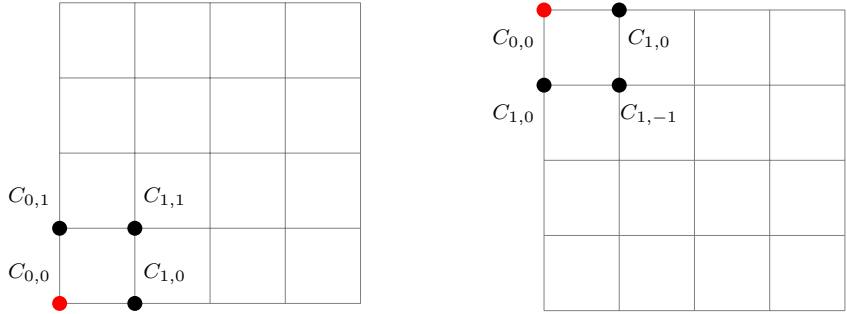
\begin{figure}[htbp]
		\hspace{-2cm}
		\begin{subfigure}[b]{0.3\textwidth}
			\begin{tikzpicture}[scale = 2]
				\draw[help lines,step = 0.5]
				(-1,-1) grid (1,1);
				\node at (-1,-0.5)[circle,fill,inner sep=2pt,color=black]{};
				\node at (-1,-1)[circle,fill,inner sep=2pt,color=red]{};
				\node at (-0.5,-0.5)[circle,fill,inner sep=2pt,color=black]{};
				\node at (-0.5,-1)[circle,fill,inner sep=2pt,color=black]{};
				\node (A) at (-1.2,-0.3) {{\tiny{$C_{0,1}$}}};
				\node (A) at (-1.2,-0.8) {{\tiny{$C_{0,0}$}}};	
				\node (A) at (-0.3,-0.3) {{\tiny{$C_{1,1}$}}};
				\node (A) at (-0.3,-0.8) {{\tiny{$C_{1,0}$}}};	
			\end{tikzpicture}
		\end{subfigure}
		\begin{subfigure}[b]{0.3\textwidth}
			\hspace{0.8cm}
			\begin{tikzpicture}[scale = 2]
				\draw[help lines,step = 0.5]
				(-1,-1) grid (1,1);
				\node at (-1,1)[circle,fill,inner sep=2pt,color=red]{};
				\node at (-1,0.5)[circle,fill,inner sep=2pt,color=black]{};
				\node at (-0.5,1)[circle,fill,inner sep=2pt,color=black]{};
				\node at (-0.5,0.5)[circle,fill,inner sep=2pt,color=black]{};
				\node (A) at (-1.2,0.8) {{\tiny{$C_{0,0}$}}};
				\node (A) at (-1.2,0.3) {{\tiny{$C_{1,0}$}}};
				\node (A) at (-0.3,0.8) {{\tiny{$C_{1,0}$}}};	
				\node (A) at (-0.3,0.3) {{\tiny{$C_{1,-1}$}}};
			\end{tikzpicture}
		\end{subfigure}
		\caption
		{The general 4-point scheme centered at $(x_i,y_j)=(x_0,y_0)=(l_1,l_3)$ for \cref{thm:corner:1} (left), and the symmetric 4-point scheme centered at $(x_i,y_j)=(x_0,y_{N_2})=(l_1,l_4)$ for \cref{thm:corner:2} (right). Red points are center points.  }
		\label{symmetric:Ckl:4}
	\end{figure}

	The corners coming from other boundary configurations can be handled in a similar way. When a corner involves at least one Dirichlet boundary condition, we can use \cref{thm:regular:Robin:1} and subsequent remarks to handle it.  In what follows, we discuss in detail how the bottom and top left stencils are constructed. The following two theorems provide the compact 4-point stencils of consistency order at least six with reduced pollution effect for the left corners. Their proofs are deferred to \cref{sec:proofs}.
	%
	%
	\begin{theorem}  \label{thm:corner:1}
	 Consider the following compact 4-point stencil centered at the corner point $(x_0,y_0)$ (see \cref{fig:config} and the  left panel of \cref{symmetric:Ckl:4}):
		\be \label{stencil:corner:1}
		\begin{aligned}
	h^{-1}	\mathcal{L}_h u_h  :=
		\begin{aligned}	
				h^{-1}
\sum_{k=0}^1 \sum_{\ell=0}^1 C_{k,\ell}(u_{h})_{k,\ell}
		\end{aligned}
		= \sum_{(m,n)\in \ind_{6}} f^{(m,n)}J_{m,n} + \sum_{n=0}^{7}g_{1}^{(n)}J_{g_{1},n} + \sum_{n=0}^{7}g_{3}^{(n)}J_{g_{3},n},
		\end{aligned}	
		\ee
		where
		{\footnotesize{
			\be \label{CR1:Corner:1}
	\begin{aligned}
		&C_{1,1} =    1
		-\bigg(15 \cdot \frac{112+219\ia }{2^{13}}-\frac{16\ia }{15}\bigg) \ka h
		+\bigg( \frac{961-419\ia}{ 2^{12}} -\frac{16}{45} \bigg) (\ka h)^2
		 +\frac{721+282\ia}{2^{15}}(\ka h)^3
		-\frac{181+51\ia }{2^{15}}(\ka h)^4,
		\\
		&C_{0,1} = 2
		-\bigg(15\cdot \frac{112+219\ia }{2^{12}}-\frac{29\ia }{15}\bigg)\ka h
		+\bigg( \frac{3187}{ 2^{13}} -5\cdot \frac{67\ia }{2^{11}}-\frac{11}{18}\bigg)(\ka h)^2
        +\bigg(\frac{1059+4899\ia }{ 2^{15}} -\frac{17\ia }{90}\bigg)(\ka h)^3 \\
        & \qquad \quad
		+ \frac{507+606\ia }{2^{15}}(\ka h)^4,
		\\
		&C_{1,0} = 2-\bigg( 15\cdot \frac{112+219\ia}{2^{12}}-\frac{29\ia}{15}\bigg)\ka h
		+\bigg(      \frac{3187}{ 2^{13}} -5\cdot \frac{67\ia}{2^{11}}-\frac{49}{90}\bigg)(\ka h)^2
		+\bigg(3\cdot \frac{1341\ia}{2^{15}}-\frac{\ia}{10} +\frac{611}{ 2^{15}}\bigg)(\ka h)^3 \\
		& \qquad \quad
		+\frac{-208+105\ia}{2^{15}}(\ka h)^4
		\\
		&C_{0,0} = -5+\bigg( 75\cdot  \frac{112+219\ia}{2^{13}}-\frac{29\ia}{15} \bigg)\ka h+ \bigg( \frac{1559-1522\ia}{2^{13}}
		+\frac{1}{90}\bigg)(\ka h)^2 -\bigg(\frac{3759}{2^{15}}+\frac{4239\ia}{2^{14}}-\frac{7\ia}{18}\bigg)(\ka h)^3 \\
		& \qquad \quad
		-\bigg(  \frac{775+324\ia}{2^{15}}-\frac{1}{40}\bigg) (\ka h)^4,
	\end{aligned}
	\ee	}}
	$(x_i^*,y_j^*)=(x_0,y_0)$, and $J_{m,n}$, $J_{g_1,n}$, $J_{g_3,n}$ are defined in \eqref{corner:1:righthand} with $\tilde{M}=8$. Then, the finite difference scheme in \eqref{CR1:Corner:1} achieves sixth consistency order for $\B_1u:=\frac{\partial u}{\partial \nv}- \ia \ka u=g_1$ and $\B_3u:=\frac{\partial u}{\partial \nv}=g_3$ at the point $(x_0,y_0)$ with reduced pollution effect.
	\end{theorem}


	\begin{theorem}  \label{thm:corner:2}
		Consider the following compact 4-point stencil centered at the corner point $(x_0,y_{N_2})$ (see \cref{fig:config} and the  right panel of \cref{symmetric:Ckl:4}):
		\be \label{stencil:corner:2}
			\begin{aligned}
		h^{-1}	\mathcal{L}_h u_h :=
			\begin{aligned}	
							h^{-1}
			\sum_{k=0}^1 \sum_{\ell=-1}^0 C_{k,\ell}(u_{h})_{k,N_2+\ell}
			\end{aligned}
			= \sum_{(m,n)\in \ind_{6}} f^{(m,n)}J_{m,n} + \sum_{n=0}^{7}g_{1}^{(n)}J_{g_{1},n} + \sum_{n=0}^{7}g_{4}^{(n)}J_{g_{4},n},
			\end{aligned}
		\ee
		where
		{\small{		
		\be \label{CR2:Corner:2}
	\begin{aligned}
		& C_{0,-1}=C_{1,0},\\
		&C_{1,-1} =
		1- \bigg( 3 \cdot \frac{293}{ 2^{13}} -\frac{5}{47} \frac{  16381\ia}{ 2^{11}}-\frac{2}{47} \frac{3467\ia}{315} \bigg) \ka h
		-  \bigg( \frac{5339 \ia}{ 5\cdot  2^{15}}+\frac{111547}{141 \cdot 2^{11}}  +\frac{100}{1269} \bigg)(\ka h)^2 \\
		& \qquad \quad
		-\frac{3+898\ia}{2^{14}}(\ka h)^3
		- \frac{1}{20}  \frac{1220+1281\ia}{2^{15}}  (\ka h)^4,
		\\
		&C_{1,0} =  2
		-  \bigg(  3\cdot \frac{293}{2^{12}} -\frac{5}{47} \frac{16381\ia}{ 2^{10}}- \frac{2}{47} \frac{3973\ia}{315}\bigg) \ka h
		-  \bigg(   \frac{1823\ia}{5 \cdot 2^{14}}+ \frac{15601}{141\cdot 2^{8}}  + \frac{10979}{88830} \bigg)(\ka h)^2 \\
		& \qquad \quad
		+  \bigg(  \frac{25}{ 2^{13}} -\frac{1}{47}\frac{3089\ia}{3\cdot 2^{8}}- \frac{1}{47}\frac{2581\ia}{1890}  \bigg)(\ka h)^3
		+  \bigg(  \frac{903\ia}{5\cdot  2^{15}} +  \frac{36461}{47\cdot 2^{15}} - \frac{79}{29610}  \bigg)(\ka h)^4,
		\\
		&C_{0,0} =  -5
		+  \bigg(  \frac{4}{47} \frac{16501\ia}{315} -\frac{25}{47} \frac{16381\ia}{ 2^{11}}+ 15\cdot \frac{293}{ 2^{13}} \bigg) \ka h
		-   \bigg(   \frac{92849 \ia}{ 5\cdot 2^{15}} +  \frac{1113127}{141\cdot 2^{11}} - \frac{23}{10} \frac{3151}{8883} \bigg)(\ka h)^2 \\
		& \qquad \quad
		-   \bigg(  \frac{1}{47} \frac{16691\ia}{945} -\frac{5}{47} \frac{165463\ia}{6\cdot 2^{12}} +  5\cdot \frac{539}{ 2^{14}} \bigg)(\ka h)^3 +    \bigg(    \frac{28811\ia}{5\cdot 2^{17}}  +    \frac{1342939}{141\cdot 2^{15}} -    \frac{2321}{40\cdot 8883}  \bigg)(\ka h)^4,
	\end{aligned}			
	\ee}}
$(x_i^*,y_j^*)=(x_0,y_{N_2})$, and $J_{m,n}$, $J_{g_1,n}$, $J_{g_4,n}$ can similarly be obtained as in \eqref{corner:1:righthand} with $\tilde{M}=8$. Then, the finite difference scheme in \eqref{stencil:corner:2} achieves seventh consistency order  accuracy for $\B_1u:=\frac{\partial u}{\partial \nv}- \ia \ka u=g_1$ and $\B_4u:=\frac{\partial u}{\partial \nv}- \ia \ka u=g_4$ at the point $(x_0,y_{N_2})$ with reduced pollution effect.
	\end{theorem}
	
	In the final subsection, we shall explicitly present our stencils  for irregular points, which have at least fifth order consistency.
	
	\subsection{Irregular points}
	
	Let $(x_i,y_j)$ be an irregular point (i.e., both $d_{i,j}^+$ and $d_{i,j}^-$ are nonempty, see the right panel of \cref{fig:irregular} for an illustration) and let us take a base point $(x^*_i,y^*_j)\in \Gamma \cap (x_i-h,x_i+h)\times (y_j-h,y_j+h)$. 
	By \eqref{base:pt}, we have
	\begin{equation} \label{base:pt:gamma}
		x_i^*=x_i-v_0h  \quad \mbox{and}\quad y_j^*=y_j-w_0h  \quad \mbox{with}\quad
		-1<v_0, w_0<1 \quad \mbox{and}\quad (x_i^*,y_j^*)\in \Gamma.
	\end{equation}
	Let $u_{\pm}$, $f_{\pm}$, $\ka_{\pm}$ represent the solution $u$, source term $f$, and wavenumber $\ka$ in $\Omega_{\pm}$.
	Similar to \eqref{ufmn}, the following notations are used
	\begin{align*}
		& u_{\pm}^{(m,n)}:=\frac{\partial^{m+n} u_{\pm}}{ \partial^m x \partial^n y}(x^*_i,y^*_j),\qquad f_{\pm}^{(m,n)}:=\frac{\partial^{m+n} f_{\pm}}{ \partial^m x \partial^n y}(x^*_i,y^*_j).
	\end{align*}
	
	Since the interface curve $\Gamma$ is smooth, the solution $u$ and the source term $f$ are assumed to be piecewise smooth, we can extend $u_+$ and $f_+$ on $\Op$ into smooth functions in a neighborhood of $(x_i^*,y_j^*)$. The same applies to $u_-$ and $f_-$ on $\Om$. As in \cite{FHM21a,FHM21b}, we assume that we have a parametric equation for $\Gamma$ on the base point $(x_i^*,y_j^*)$. I.e.,
	\be \label{parametric}
	x=r(t),\quad y=s(t), \quad \mbox{for} \quad (x,y)\in \Gamma,
	\ee
	where $r(t)$ and $s(t)$ are smooth functions. For $(x_i^*,y_j^*)\in \Gamma$, there exists a $t_k^*$ such that
	\be\label{t:star}
	x^*_i=r(t_k^*),\quad y^*_j=s(t_k^*), \quad \mbox{and} \quad (r'(t_k^*))^2+(s'(t_k^*))^2\ne0 \quad \mbox{for some } k\in \N.
	\ee
	\begin{theorem}\label{thm:interface}
		Let $u$ be the solution to the Helmholtz interface problem in \eqref{Qeques2} and the base point $(x_i^*,y_j^*)\in \Gamma$ be parameterized by \eqref{parametric} and \eqref{t:star}.
		Then
		\be	\label{transmission:Hel}
		\begin{split}
			 u_-^{(m',n')}&=\sum_{(m,n)\in \ind_{M}^1} T^{u_{+}}_{m',n',m,n}u_+^{(m,n)}+\sum_{(m,n)\in \ind_{M-2}} \left(T^{+}_{m',n',m,n} f_+^{(m,n)}
			+ T^{-}_{m',n',m,n} f_{-}^{(m,n)}\right)  \\
			&\quad +\sum_{p=0}^{M} T^{g}_{m',n',p} g^{(p)} +\sum_{p=0}^{M-1} T^{g_{\Gamma}}_{m',n',p} g^{(p)}_{\Gamma},\qquad \forall\; (m',n')\in \ind_{M}^{1},
		\end{split}
		\ee
		where
		\[
g^{(p)}:=\frac{1}{p!} \frac{d^p}{dt^p} (g(t))
\bigg|_{t=t_k^*},\qquad p=0,1,\dots,M,
\]
\[
g^{(p)}_{\Gamma}:=\frac{1}{p!} \frac{d^p}{dt^p} \left(g_\Gamma(t)\sqrt{(r'(t))^2+(s'(t))^2}\right)
\bigg|_{t=t_k^*},\qquad p=0,1,\dots,M-1,
\]
	all the transmission coefficients $T^{u_{+}}, T^{\pm}, T^{g}, T^{g_{\Gamma}}$ are uniquely determined by
	$r^{(p)}(t_k^*)$, $s^{(p)}(t_k^*)$ for $p=0,\ldots,M$ and ${\ka}_{\pm}$. In particular, if $\ka_{+}=\ka_{-}\ge 0$, then
	\be\label{transmission:Hel:equal}
	T^{u_{+}}_{m',n',m,n}=
	\begin{cases}
		1, &\text{if $m=m'$ and $n=n'$},\\
		0, &\text{else}.
	\end{cases}
	\ee
For the general case, where $\ka_{+},\ka_{-}\ge 0$, we have
\be\label{transmission:Hel:not:equal}
T^{u_{+}}_{m',n',m,n}=
\begin{cases}
	1, &\text{if $(m,n)=(m',n') \in \{(0,0), (0,1), (1,0)\}$},\\
	0, &\text{if $(m,n)=(0,0)$ and $(m',n') \in \{(0,1), (1,0)\}$},\\
	0, &\text{if $m+n>m'+n'$}.
\end{cases}
\ee
	\end{theorem}

	Next, we state the compact 9-point stencil for interior irregular points in two separate cases: the special case $\ka_{+}=\ka_{-}$ with seventh consistency order or the general case $\ka_{+},\ka_{-}\ge 0$ with fifth consistency order.
	
	\begin{theorem}\label{fluxtm2}
	Let $(x_i-v_0h,y_j-w_0h)=(x_i^*,y_j^*)\in \Gamma$ with $-1<v_0,w_0<1$. Suppose $\ka_{+}=\ka_{-} \ge 0$.
		The following compact 9-point stencil  centered at the interior irregular point $(x_i,y_j)$ (see the left panel of \cref{fig:irregular})
		\be\label{irregular:same:k}
		\begin{split}
			h^{-1}\mathcal{L}_h u_h :&=
			h^{-1}
		\sum_{k=-1}^1 \sum_{\ell=-1}^1 C_{k,\ell}(u_{h})_{i+k,j+\ell}\\
			&=\sum_{(m,n)\in \ind_{5} } f_+^{(m,n)}J^{+}_{m,n} + \sum_{(m,n)\in \ind_{5}} f_-^{(m,n)}J^{-}_{m,n} +\sum_{p=0}^7 g^{(p)}J^{g}_p + \sum_{p=0}^6 g^{(p)}_{\Gamma}J^{g_{\Gamma}}_p,
		\end{split}
	\ee
		where $\{C_{k,\ell}\}_{k,\ell\in\{-1,0,1\}}$ are defined in \eqref{stencil:Cv},
		\be \label{irregular:same:k:right}
			\begin{split}
				& J^{\pm}_{m,n}:=
				 J_{m,n}^{\pm,0}+J^{\pm,T}_{m,n}, \quad J^{\pm,0}_{m,n}:=h^{-1}\sum_{(k,\ell)\in d_{i,j}^\pm} C_{k,\ell} H^{\pm}_{7,m,n}((v_0 +k)h,(w_0+\ell) h),\\
				&   J^{\pm,T}_{m,n}:=h^{-1}
				\sum_{(m',n')\in \ind_{7}^{1}} I^{-}_{m',n'} T^{\pm}_{m',n',m,n}, \quad I^{-}_{m,n}:=\sum_{(k,\ell)\in d_{i,j}^-}
				C_{k,\ell} G^{-}_{7,m,n}((v_0+k)h,(w_0+\ell) h), \\
				& J^{g}_p:=h^{-1}
				\sum_{(m',n')\in \ind_{7}^{1}} I^{-}_{m',n'} T^{g}_{m',n',p}, \qquad
				J^{g_{\Gamma}}_p:=h^{-1}
				\sum_{(m',n')\in \ind_{7}^{1}} I^{-}_{m',n'} T^{g_{\Gamma}}_{m',n',p},
			\end{split}
		\ee
		$H^{\pm}_{7,m,n}$, $G^{-}_{7,m,n}$ are defined in \eqref{Hel:GM1mnxy}-\eqref{Hel:HM1mnxy} with $\ka$ being replaced by $\ka_{\pm}$, and $T^{\pm}_{m',n',m,n}$, $T^{g}_{m',n',p}$, $T^{g_{\Gamma}}_{m',n',p}$ are transmission coefficients in \eqref{transmission:Hel}, achieves seventh consistency order  for $\left[u\right]=g$ and $\left[\nabla  u \cdot \nv \right]=g_{\Gamma}$ on $\Gamma$.
	\end{theorem}
		\begin{theorem}\label{fluxtm3}
		Let $(x_i-v_0h,y_j-w_0h)=(x_i^*,y_j^*)\in \Gamma$ with $-1<v_0,w_0<1$.
		Suppose $\ka_{+},\ka_{-}\ge 0$. The following compact 9-point stencil centered at the interior irregular point $(x_i,y_j)$ (see the middle panel of \cref{fig:irregular})
		\be\label{irregular:diff:k}
\begin{split}
			h^{-1}\mathcal{L}_h u_h :&=
			h^{-1}
			\sum_{k=-1}^1 \sum_{\ell=-1}^1 C_{k,\ell}(u_{h})_{i+k,j+\ell}\\
			&=\sum_{(m,n)\in \ind_{3} } f_+^{(m,n)}J^{+}_{m,n} + \sum_{(m,n)\in \ind_{3}} f_-^{(m,n)}J^{-}_{m,n} +\sum_{p=0}^5 g^{(p)}J^{g}_p + \sum_{p=0}^4 g^{(p)}_{\Gamma}J^{g_{\Gamma}}_p,
		\end{split}
\ee
		where $\{C_{k,\ell}\}_{k,\ell\in\{-1,0,1\}}$ are obtained by solving \eqref{System:Ckl} with $M=5$, 
		%
		\be \label{irregular:diff:k:right}
	\begin{split}
			& J^{\pm}_{m,n}:=
			J_{m,n}^{\pm,0}+J^{\pm,T}_{m,n}, \quad J^{\pm,0}_{m,n}:=h^{-1}\sum_{(k,\ell)\in d_{i,j}^\pm} C_{k,\ell} H^{\pm}_{5,m,n}((v_0 +k)h,(w_0+\ell) h),\\
			&   J^{\pm,T}_{m,n}:=h^{-1}
			\sum_{(m',n')\in \ind_{5}^{1}} I^{-}_{m',n'} T^{\pm}_{m',n',m,n}, \quad I^{-}_{m,n}:=\sum_{(k,\ell)\in d_{i,j}^-}
			C_{k,\ell} G^{-}_{5,m,n}((v_0+k)h,(w_0+\ell) h), \\
			& J^{g}_p:=h^{-1}
			\sum_{(m',n')\in \ind_{5}^{1}} I^{-}_{m',n'} T^{g}_{m',n',p}, \qquad
			J^{g_{\Gamma}}_p:=h^{-1}
			\sum_{(m',n')\in \ind_{5}^{1}} I^{-}_{m',n'} T^{g_{\Gamma}}_{m',n',p},
			\end{split}
	\ee
		$H^{\pm}_{5,m,n}$, $G^{-}_{5,m,n}$ are defined in \eqref{Hel:GM1mnxy}-\eqref{Hel:HM1mnxy} with $\ka$ being replaced by $\ka_{\pm}$, and $T^{\pm}_{m',n',m,n}$, $T^{g}_{m',n',p}$, $T^{g_{\Gamma}}_{m',n',p}$ are transmission coefficients in \eqref{transmission:Hel}, achieves fifth consistency order  for $\left[u\right]=g$ and $\left[\nabla  u \cdot \nv \right]=g_{\Gamma}$ on $\Gamma$.
	\end{theorem}
Depending on how the interface curve partitions the 9 points in it, there exist many configurations for the scheme in \eqref{irregular:diff:k}. For the system of linear equations $Ax=b$ with infinitely many solutions, the MATLAB Package $\texttt{mldivide}(A,b)$ can automatically choose free parameters to be 0.
To solve all cases in \eqref{System:Ckl}, we choose
\be
\begin{split}
	& M=5,\quad C_{k,\ell}:=\sum_{p=0}^{M} c_{k,\ell,p}(\max(\ka_{+},\ka_{-}) h)^p, \quad c_{k,\ell,p} \in \R,\\
	& c_{0,0,0}=-20, \quad c_{-1,-1,0}=c_{-1,1,0}=c_{1,-1,0}=c_{1,1,0}=1, \quad c_{-1,0,0}=c_{1,0,0}=c_{0,-1,0}=c_{0,1,0}=4,
\end{split}
\ee
and use  the MATLAB Package $\texttt{mldivide}(A,b)$ to solve \eqref{System:Ckl}.
\begin{figure}[htbp]
	\centering	
	\hspace{-0.9cm}
	\begin{subfigure}[b]{0.3\textwidth}
		\begin{tikzpicture}[scale = 0.85]
			\draw[domain =0:360,smooth]
			plot({(1.2+0.4*sin(5*\x))*cos(\x)}, {(1.2+0.4*sin(5*\x))*sin(\x)});
			\draw
			(-pi, -pi) -- (-pi, pi) -- (pi, pi) -- (pi, -pi) --(-pi,-pi);

			\node (A) at (0.6,1.3) {$\Gamma$};
			\node (A) at (0,-2) {$\ka_{+}=\ka_{-}$};
			\node (A) at (0,0) {$\ka_{-}$};
			\node (A) at (-2,2) {$\ka_{+}$};
		\end{tikzpicture}
	\end{subfigure}
	\begin{subfigure}[b]{0.3\textwidth}
		\hspace{0.4cm}
		\begin{tikzpicture}[scale = 0.85]
			\draw[domain =0:360,smooth]
			plot({(1.2+0.4*sin(5*\x))*cos(\x)}, {(1.2+0.4*sin(5*\x))*sin(\x)});
			\draw
			(-pi, -pi) -- (-pi, pi) -- (pi, pi) -- (pi, -pi) --(-pi,-pi);

			\node (A) at (0.6,1.3) {$\Gamma$};
			\node (A) at (0,-2) {$\ka_{+}\ne \ka_{-}$};
			\node (A) at (0,0) {$\ka_{-}$};
			\node (A) at (-2,2) {$\ka_{+}$};
		\end{tikzpicture}
	\end{subfigure}
	\begin{subfigure}[b]{0.3\textwidth}
		\hspace{0.8cm}
		\begin{tikzpicture}[scale = 1.33]
			\draw[help lines,step = 1]
			(0,0) grid (4,4);
			\node at (1,1)[circle,fill,inner sep=2pt,color=black]{};
			\node at (1,2)[circle,fill,inner sep=2pt,color=black]{};
			\node at (1,3)[circle,fill,inner sep=2pt,color=black]{};	
			\node at (2,1)[circle,fill,inner sep=2pt,color=black]{};
			\node at (2,2)[circle,fill,inner sep=2pt,color=red]{};
			\node at (2,3)[circle,fill,inner sep=2pt,color=black]{};	
			\node at (3,1)[circle,fill,inner sep=2pt,color=black]{};
			\node at (3,2)[circle,fill,inner sep=2pt,color=black]{};
			\node at (3,3)[circle,fill,inner sep=2pt,color=black]{};
			\node (A) at (0.8,1.5) {$\Gamma$};	
			\node (A) at (2.5,1.7) {$(x_i,y_j)$};
			\draw[line width=1.5pt, blue]  plot [smooth,tension=0.8]
			coordinates {(0,1.1) (1,1.4) (2,2.4) (2.7,4)};
			\node (A) at (1.5,2.5) {$\ka_{+}$};	
			\node (A) at (2.5,0.5) {$\ka_{-}$};	
		\end{tikzpicture}
	\end{subfigure}
	\caption
	{The illustrations for \cref{thm:interface,fluxtm2,fluxtm3}. Note that the red point is the center point $(x_i,y_j)$ and the blue curve is the interface curve $\Gamma$.}
	\label{fig:irregular}
\end{figure}
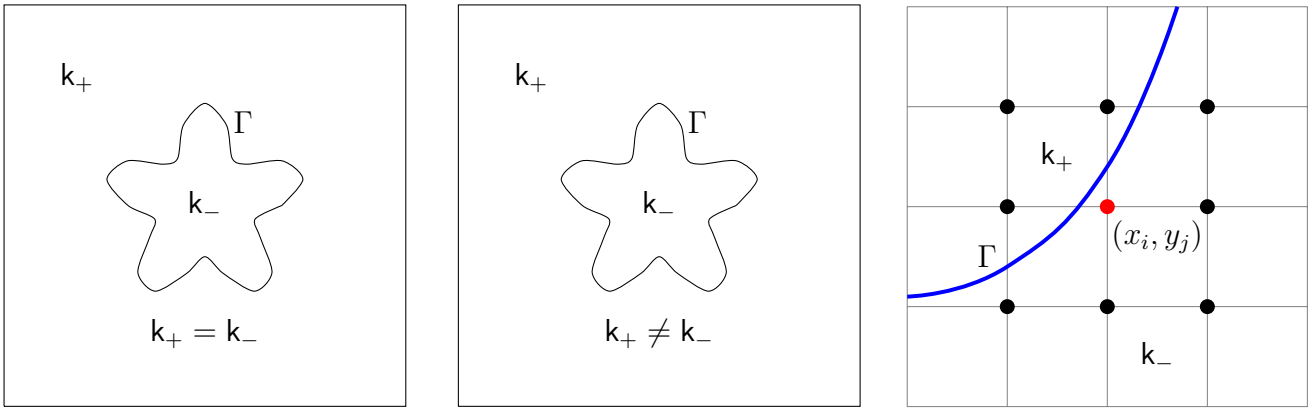
\begin{remark} 
	If we replace $\ka$ by $\ia \ka$, we will propose a discrete maximum principle preserved scheme for \eqref{Qeques2} in the future.
\end{remark}
	\section{Numerical experiments}
	\label{sec:numerical}
	
	In this section, we let $\Omega:=(l_1,l_2)^2$. For a given $J\in \NN$, we define
	$h:=(l_2-l_1)/N_1$ with $N_1:=2^J$. Recall the definition of $(x_i, y_j)$ in \eqref{xiyj}. Let $u(x,y)$ be the exact solution of \eqref{Qeques2} and $(u_{h})_{i,j}$ be the numerical solution at $(x_i, y_j)$ using the mesh size $h$.
	We shall evaluate our proposed finite difference scheme in the $l_2$ norm by the relative error
	$\frac{\|u_{h}-u\|_{2}}{\|u\|_{2}}$ if the exact solution $u$ is available,
	and
	by the error ${\|u_{h}-u_{h/2}\|_{2}}$ if the exact solution is not known, where
	\begin{align*}
		&\|u_{h}-u\|_{2}^2:= h^2\sum_{i=0}^{N_1}\sum_{j=0}^{N_1} \left((u_h)_{i,j}-u(x_i,y_j)\right)^2, \qquad \|u_{h}-u_{h/2}\|_{2}^2:= h^2\sum_{i=0}^{N_1}\sum_{j=0}^{N_1} \left((u_{h})_{i,j}-(u_{h/2})_{2i,2j}\right)^2.
	\end{align*}
In addition we also provide results for the infinity norm of the errors given by:
\[
\|u_h-u\|_\infty
:=\max_{0\le i\le N_1, 0\le j\le N_1} \left|(u_h)_{i,j}-u(x_i,y_j)\right|,
\quad
\|u_{h}-u_{h/2}\|_\infty:=\max_{0\le i\le N_1,0\le j\le N_1} \left|(u_{h})_{i,j}-(u_{h/2})_{2i,2j}\right|.
\]
 All eight examples presented below verify the theoretical findings of \cref{sec:sixord}. 
	\subsection{Numerical examples with no interfaces}
\begin{figure}[htbp]
		\centering
		\begin{subfigure}[b]{0.3\textwidth}
			 \includegraphics[width=6cm,height=5cm]{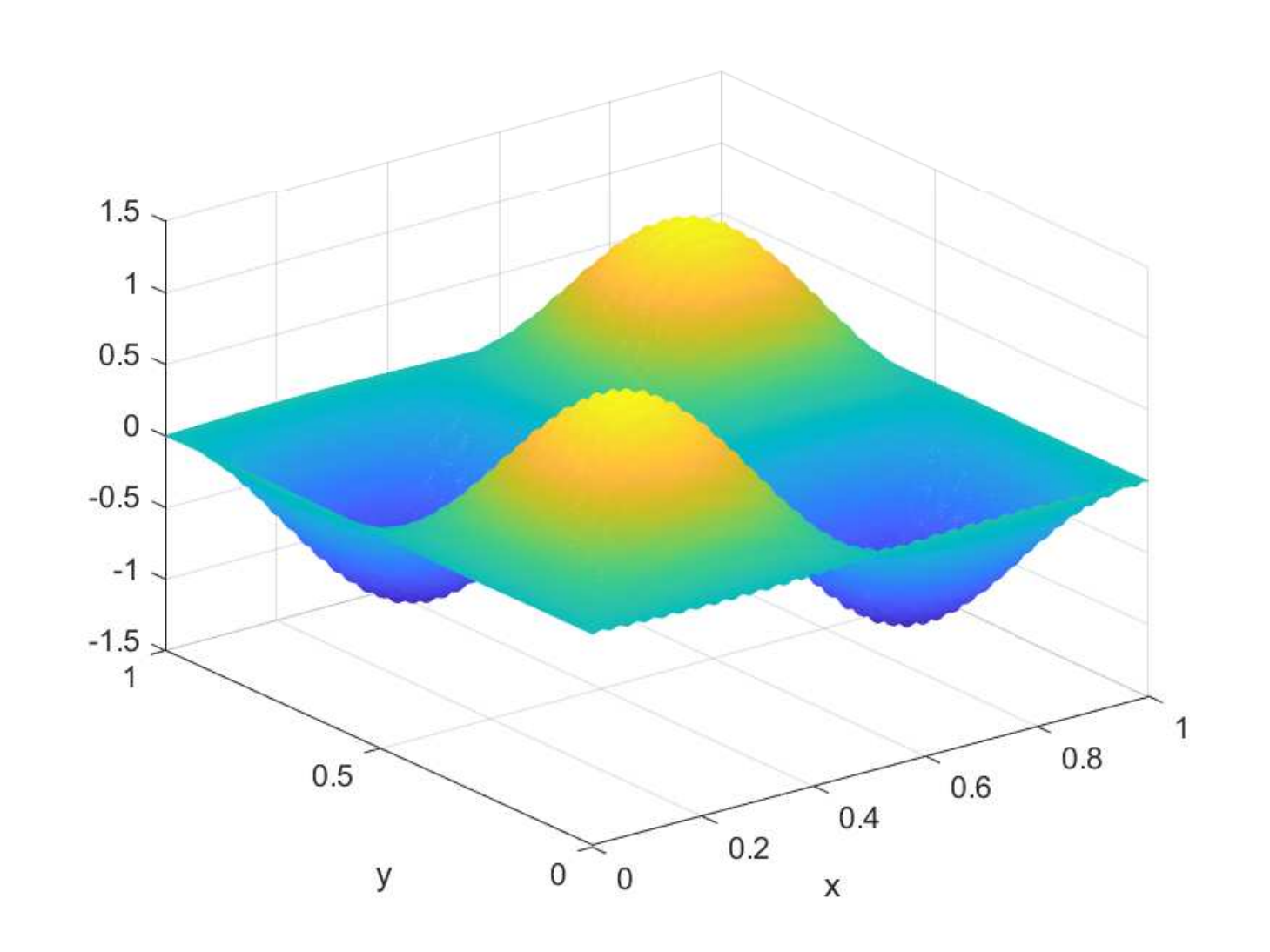}
		\end{subfigure}
		\begin{subfigure}[b]{0.3\textwidth}
			 \includegraphics[width=6cm,height=5cm]{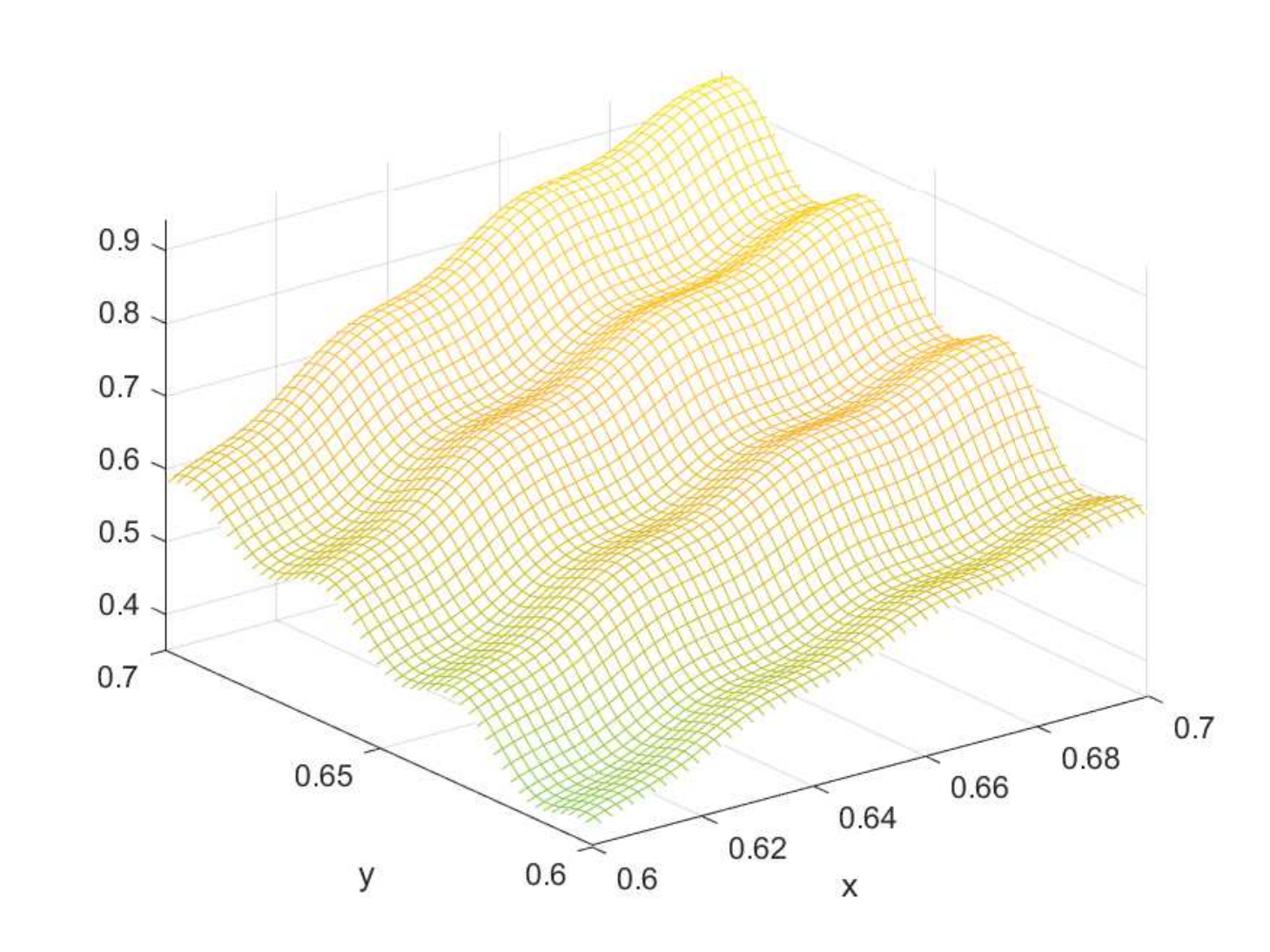}
		\end{subfigure}
		\begin{subfigure}[b]{0.3\textwidth}
			 \includegraphics[width=6cm,height=5cm]{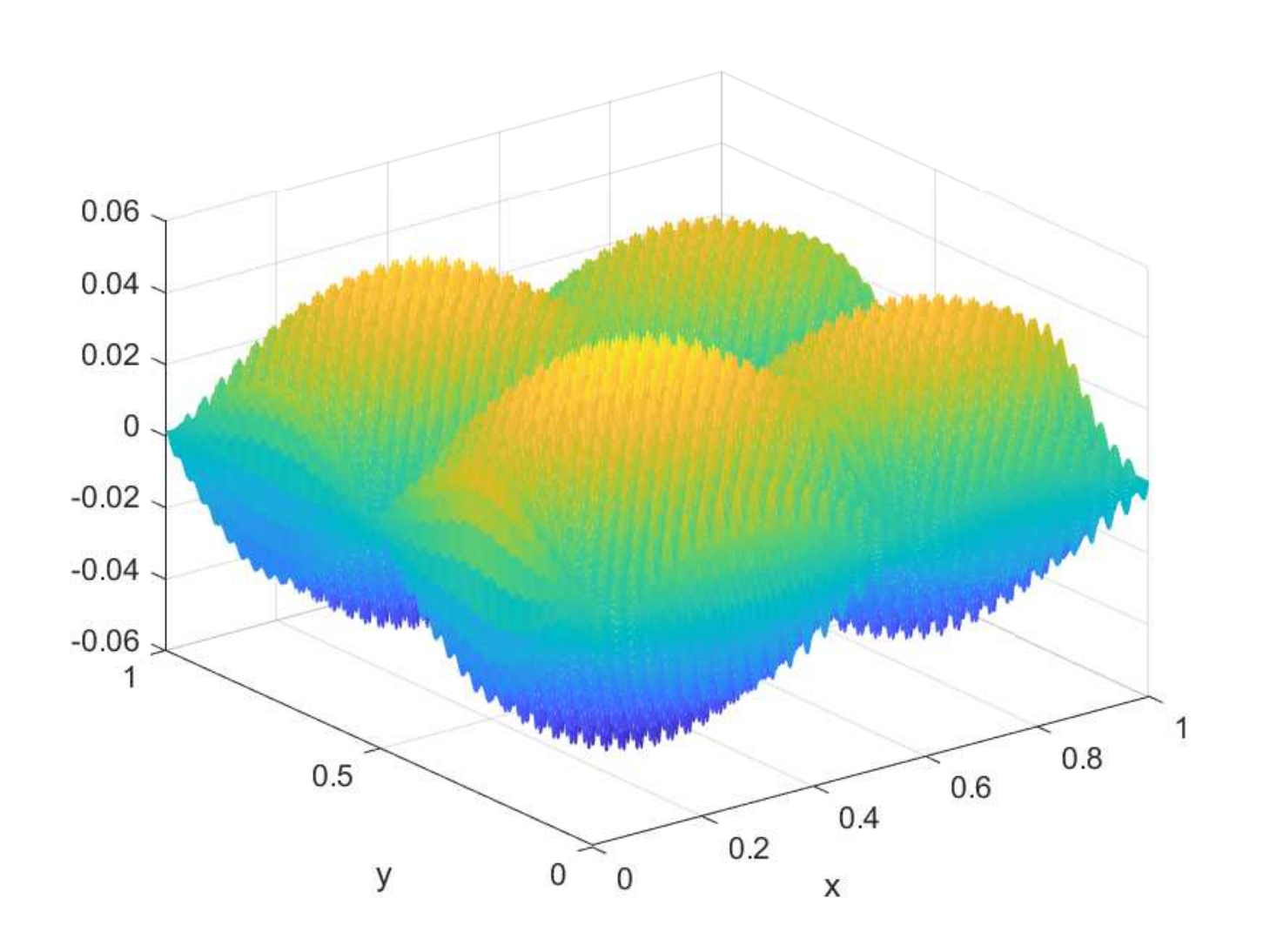}
		\end{subfigure}
		\begin{subfigure}[b]{0.3\textwidth}
			 \includegraphics[width=6cm,height=5cm]{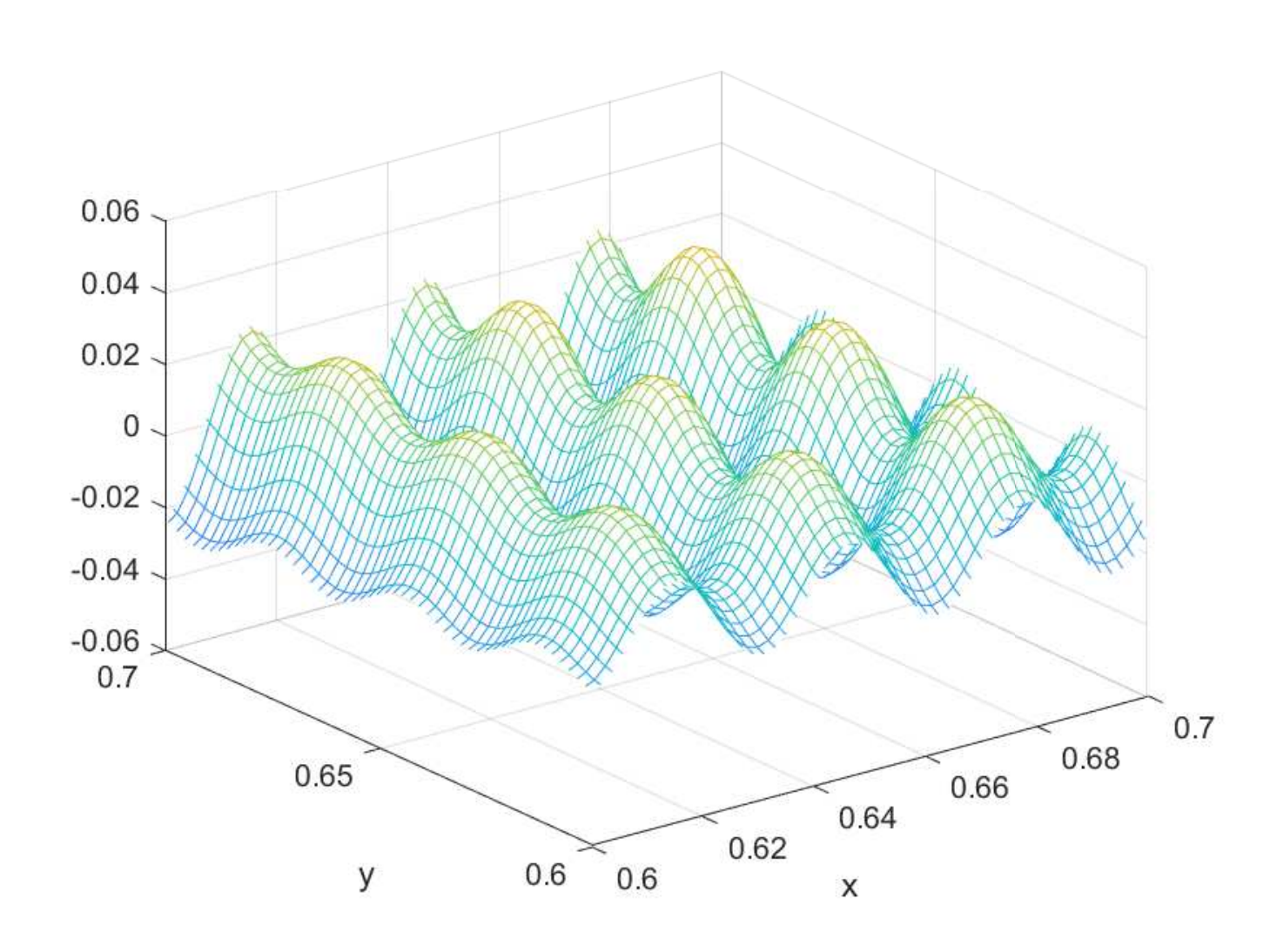}
		\end{subfigure}
		\begin{subfigure}[b]{0.3\textwidth}
			 \includegraphics[width=6cm,height=5cm]{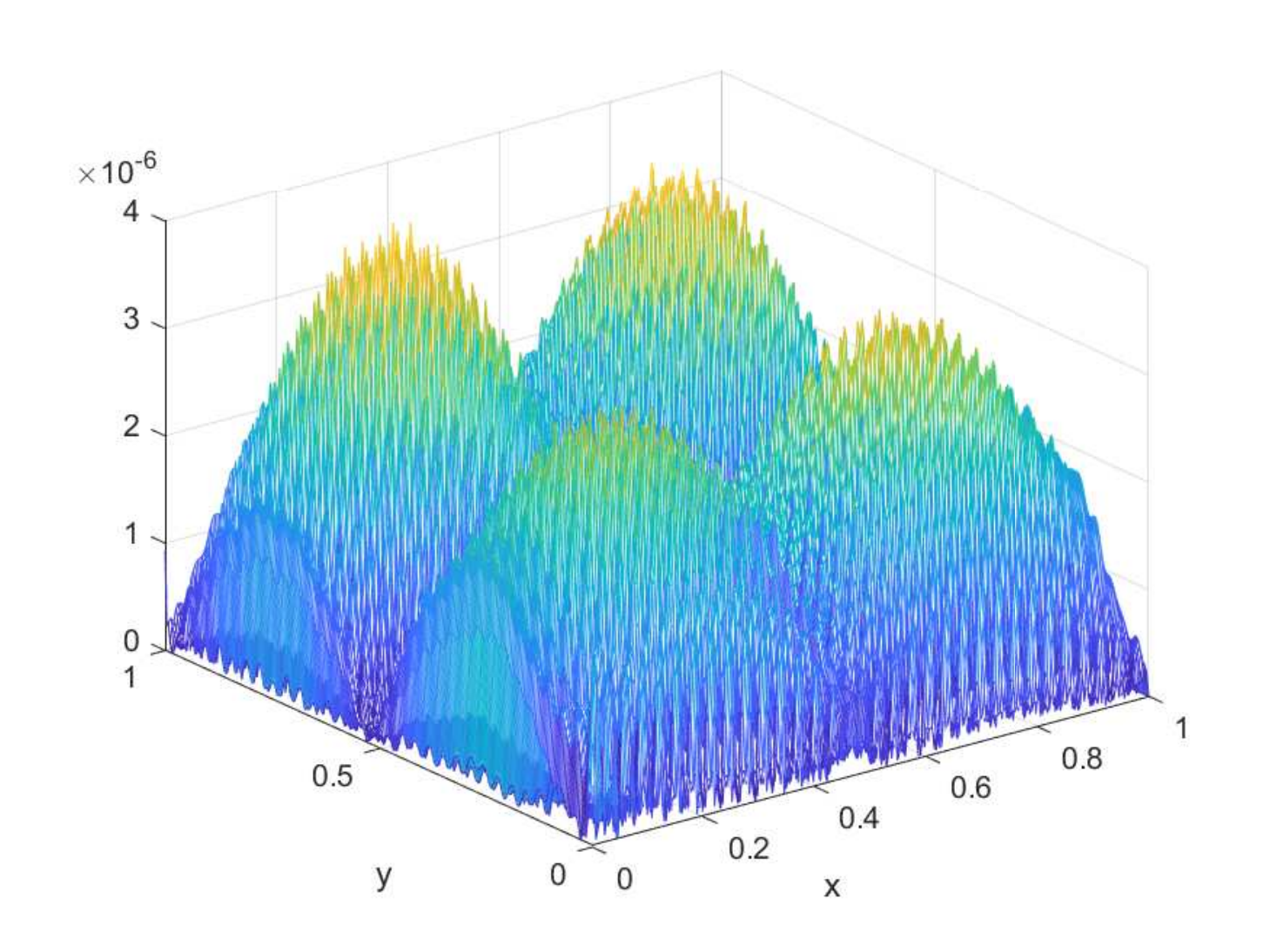}
		\end{subfigure}
		\begin{subfigure}[b]{0.3\textwidth}
			 \includegraphics[width=6cm,height=5cm]{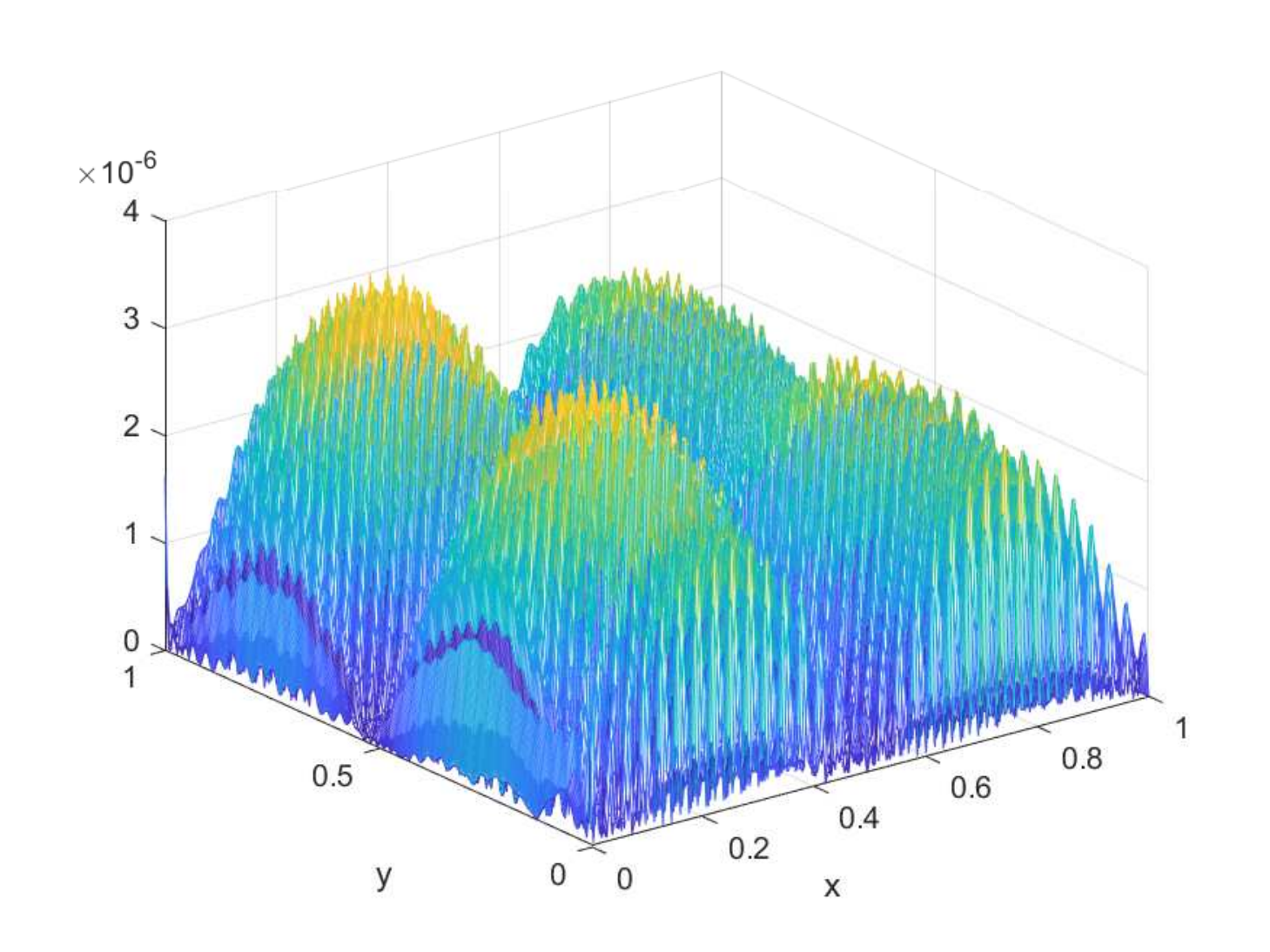}
		\end{subfigure}
		\caption
		{ \cref{ex4} with $\ka=200$ and $h=1/2^8$. First row: $\Re(u_{h/2})$ in $[0,1]^2$ (left), $\Re(u_{h/2})$ in $[0.6,0.7]^2$ (middle), $\Im(u_{h/2})$ in $[0,1]^2$  (right). Second row: $\Im(u_{h/2})$ in $[0.6,0.7]^2$ (left), $|\Re(u_h-u_{h/2})|$ in $[0,1]^2$ (middle), $|\Im(u_h-u_{h/2})|$ in $[0,1]^2$ (right).}
		\label{fig:U:graph:k200}
	\end{figure}

\begin{figure}[htbp]
	\centering
	\begin{subfigure}[b]{0.3\textwidth}
		 \includegraphics[width=6cm,height=5cm]{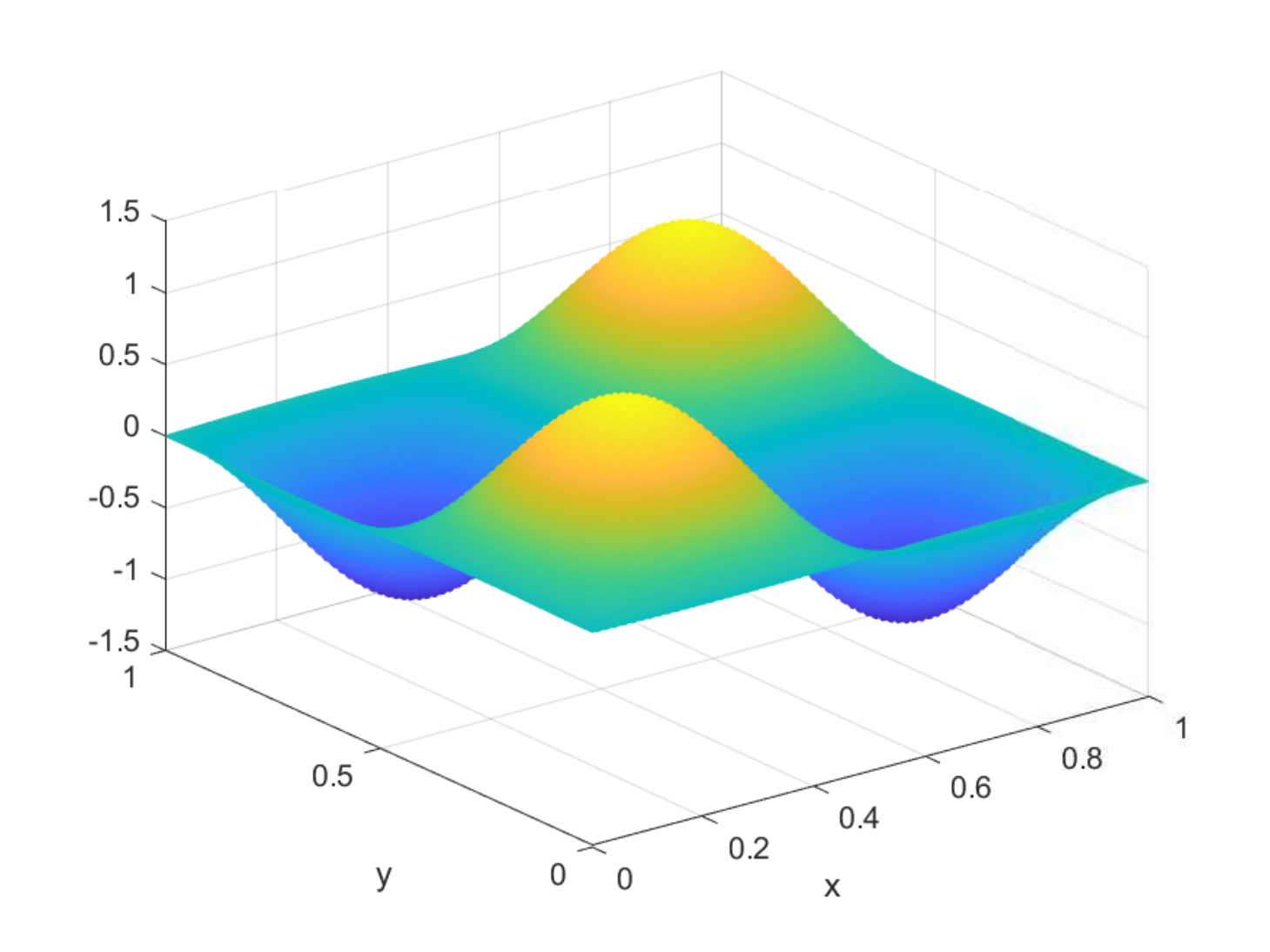}
	\end{subfigure}
	\begin{subfigure}[b]{0.3\textwidth}
		 \includegraphics[width=6cm,height=5cm]{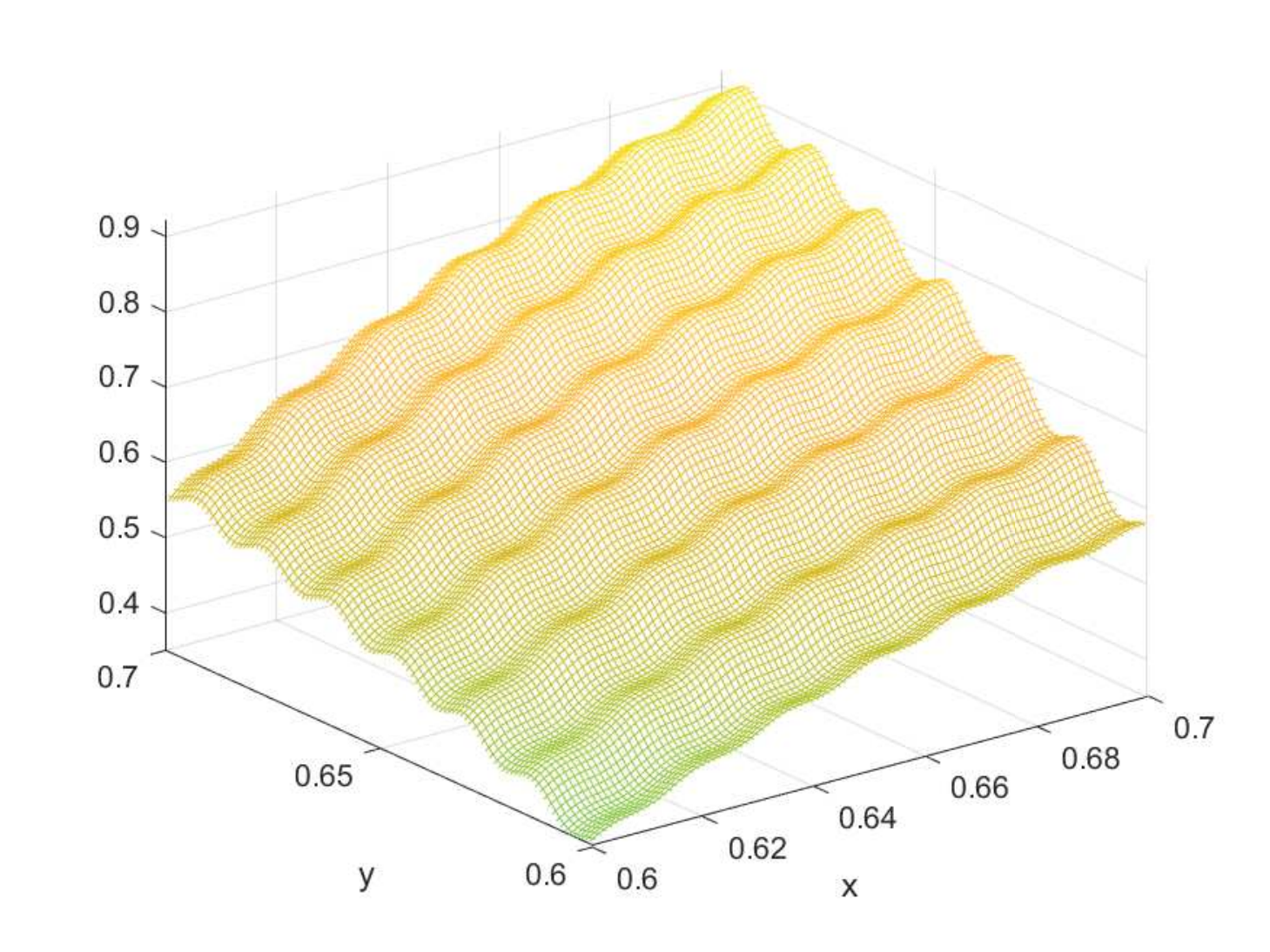}
	\end{subfigure}
	\begin{subfigure}[b]{0.3\textwidth}
		 \includegraphics[width=6cm,height=5cm]{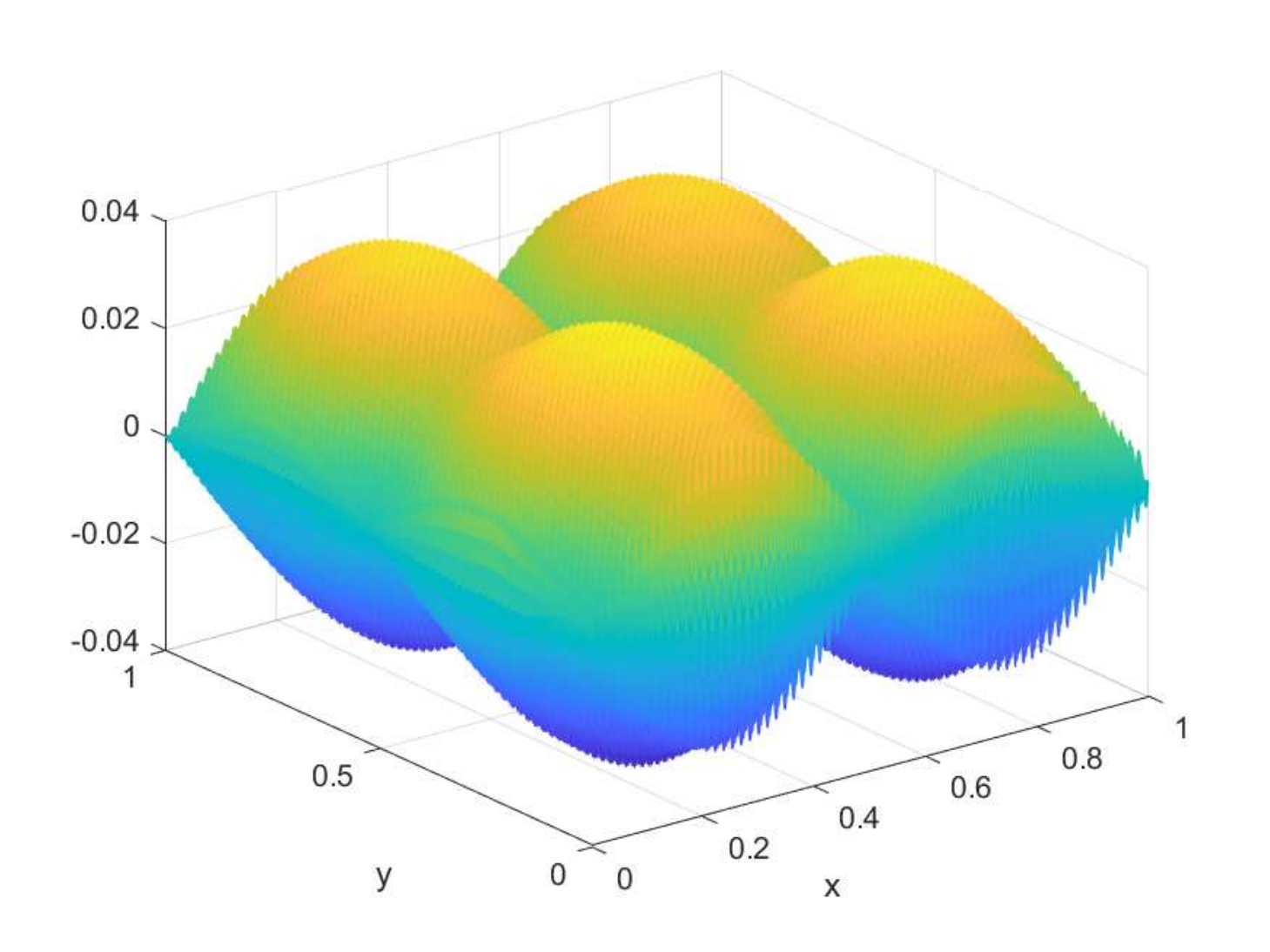}
	\end{subfigure}
	\begin{subfigure}[b]{0.3\textwidth}
		 \includegraphics[width=6cm,height=5cm]{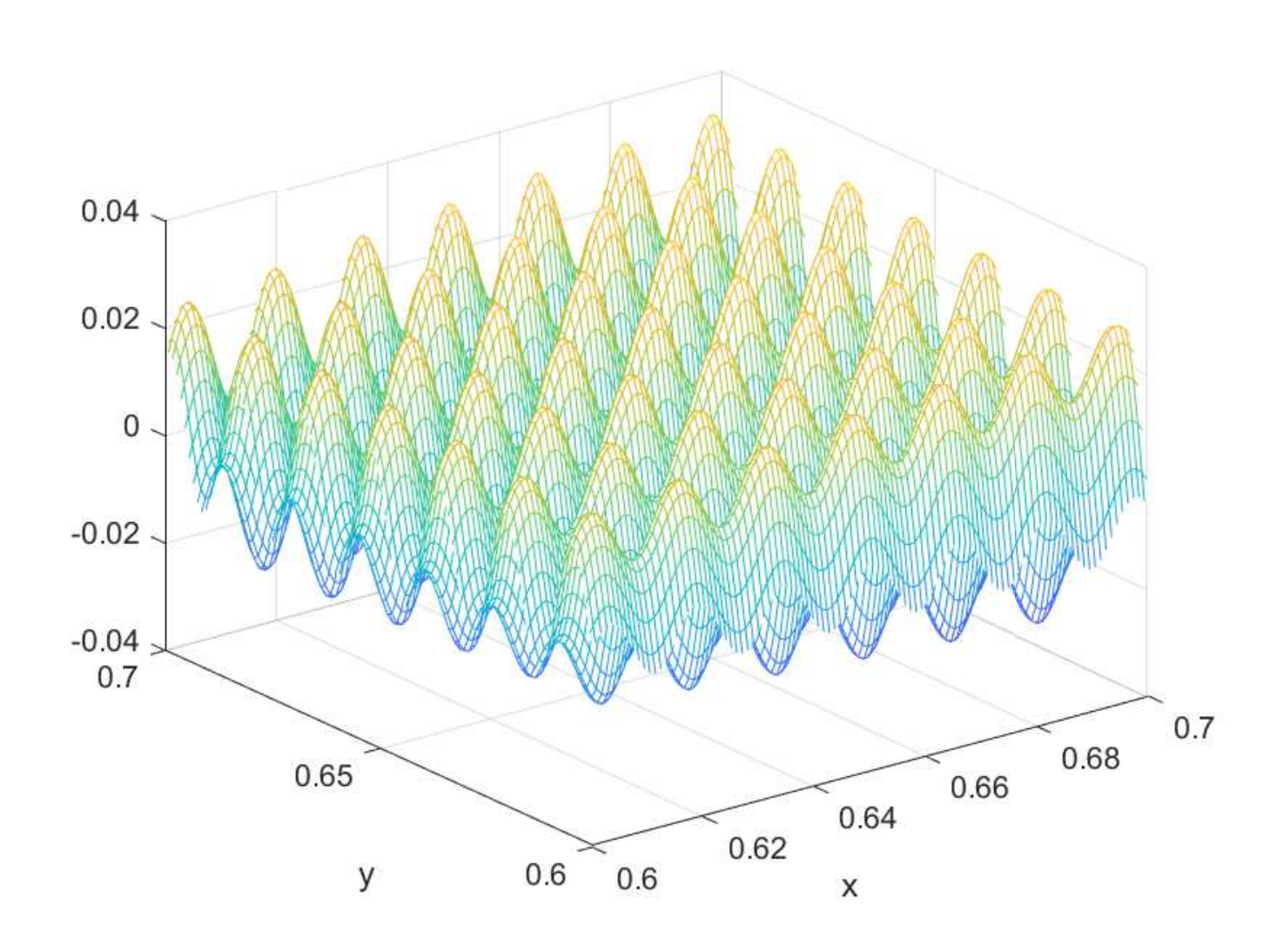}
	\end{subfigure}
	\begin{subfigure}[b]{0.3\textwidth}
		 \includegraphics[width=6cm,height=5cm]{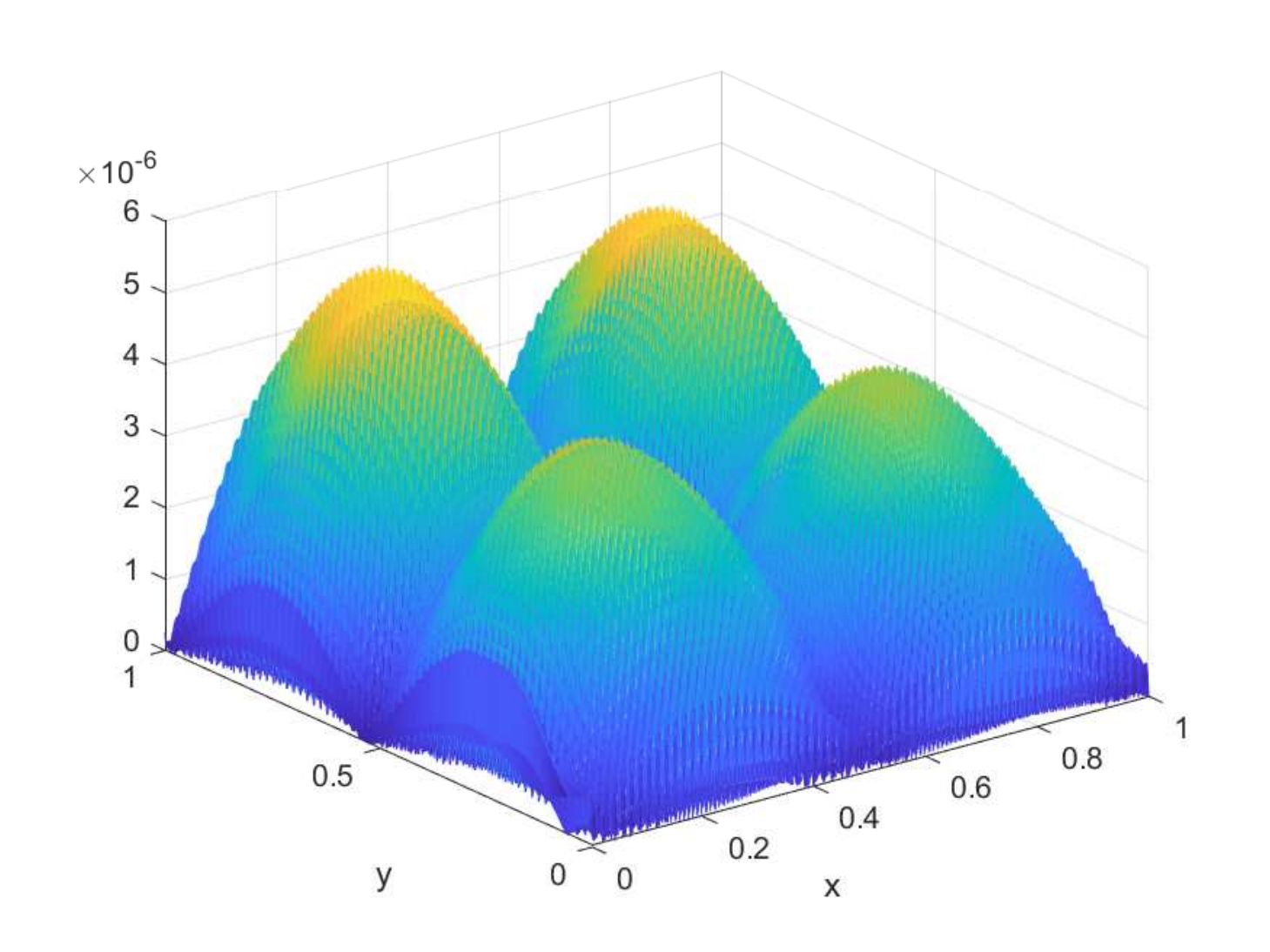}
	\end{subfigure}
	\begin{subfigure}[b]{0.3\textwidth}
		 \includegraphics[width=6cm,height=5cm]{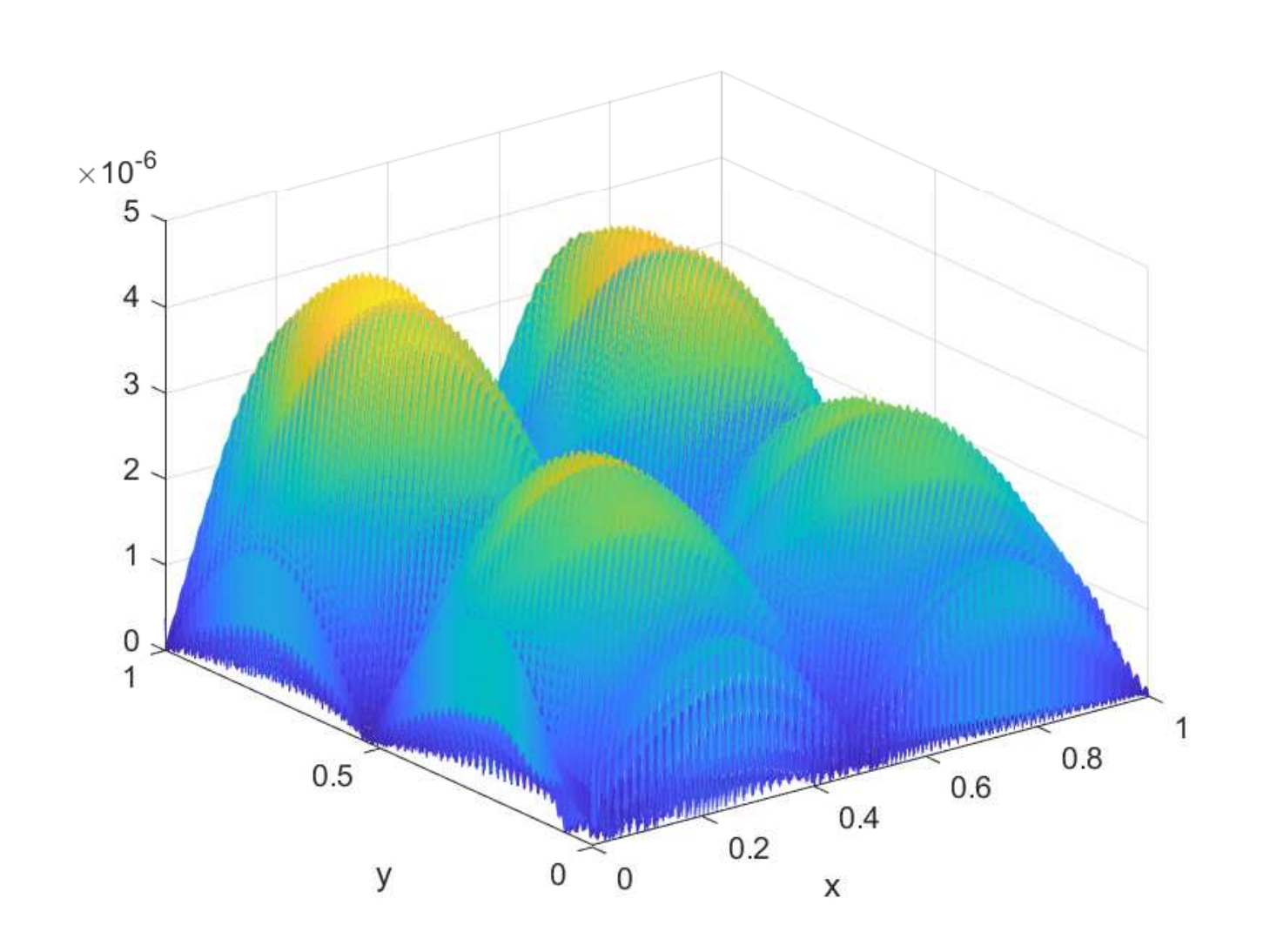}
	\end{subfigure}
	\caption
		{ \cref{ex4} with $\ka=400$ and $h=1/2^9$. First row: $\Re(u_{h/2})$ in $[0,1]^2$ (left), $\Re(u_{h/2})$ in $[0.6,0.7]^2$ (middle), $\Im(u_{h/2})$ in $[0,1]^2$ (right). Second row: $\Im(u_{h/2})$ in $[0.6,0.7]^2$ (left), $|\Re(u_h-u_{h/2})|$ in $[0,1]^2$ (middle), $|\Im(u_h-u_{h/2})|$ in $[0,1]^2$ (right).}
	\label{fig:U:graph:k400}
\end{figure}

\begin{figure}[htbp]
	\centering
	\begin{subfigure}[b]{0.3\textwidth}
		 \includegraphics[width=6cm,height=5cm]{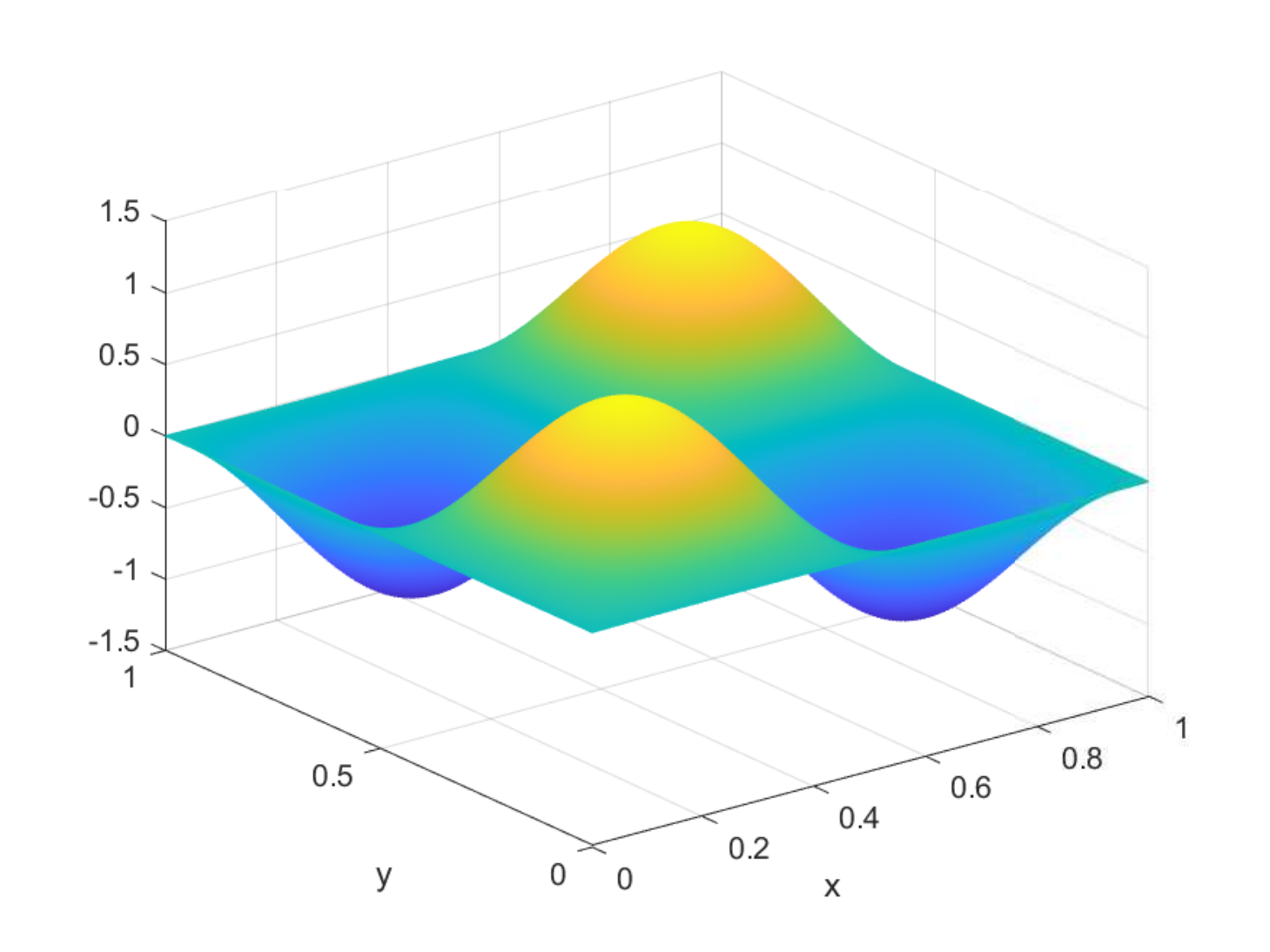}
	\end{subfigure}
	\begin{subfigure}[b]{0.3\textwidth}
		 \includegraphics[width=6cm,height=5cm]{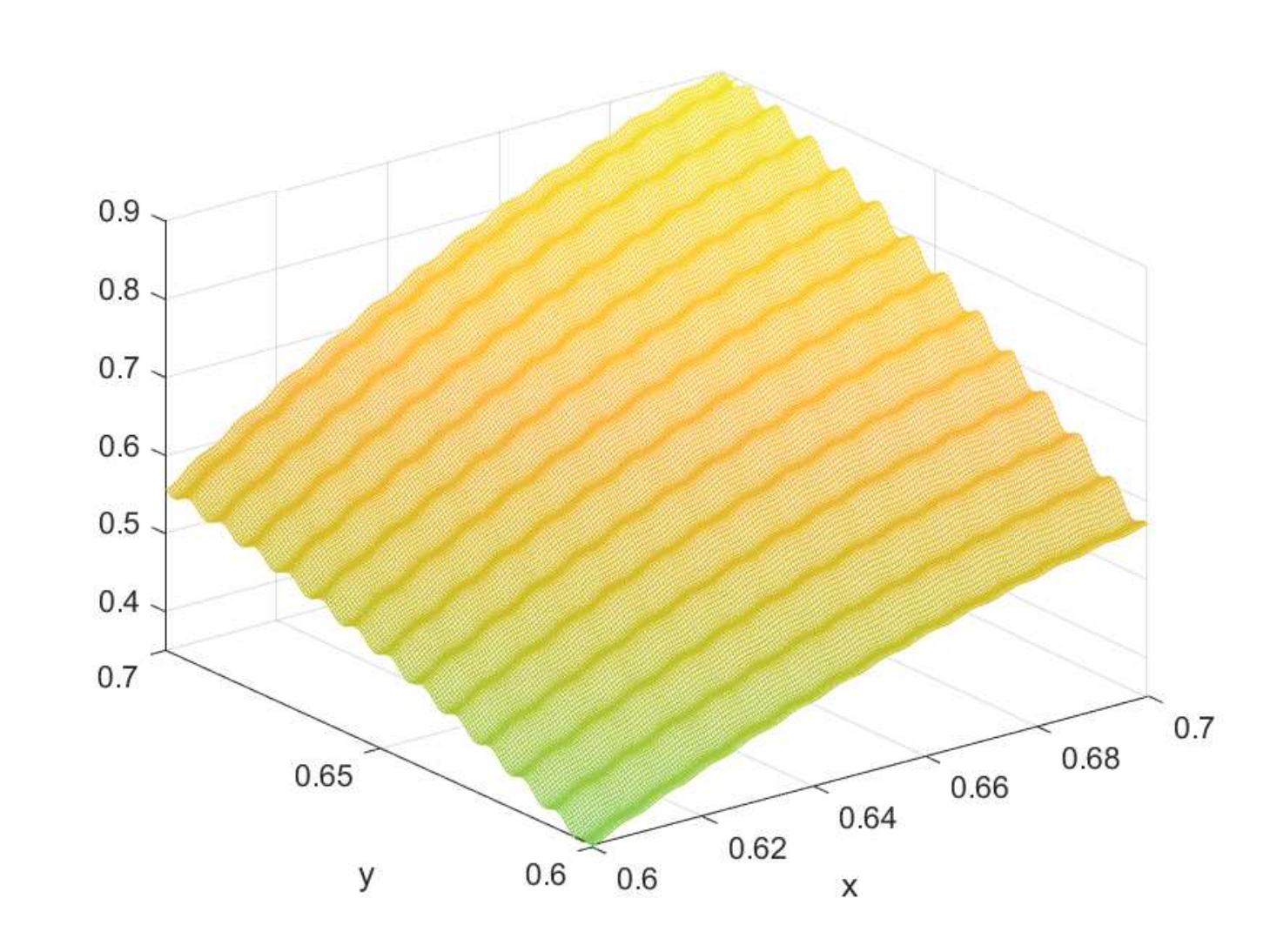}
	\end{subfigure}
	\begin{subfigure}[b]{0.3\textwidth}
		 \includegraphics[width=6cm,height=5cm]{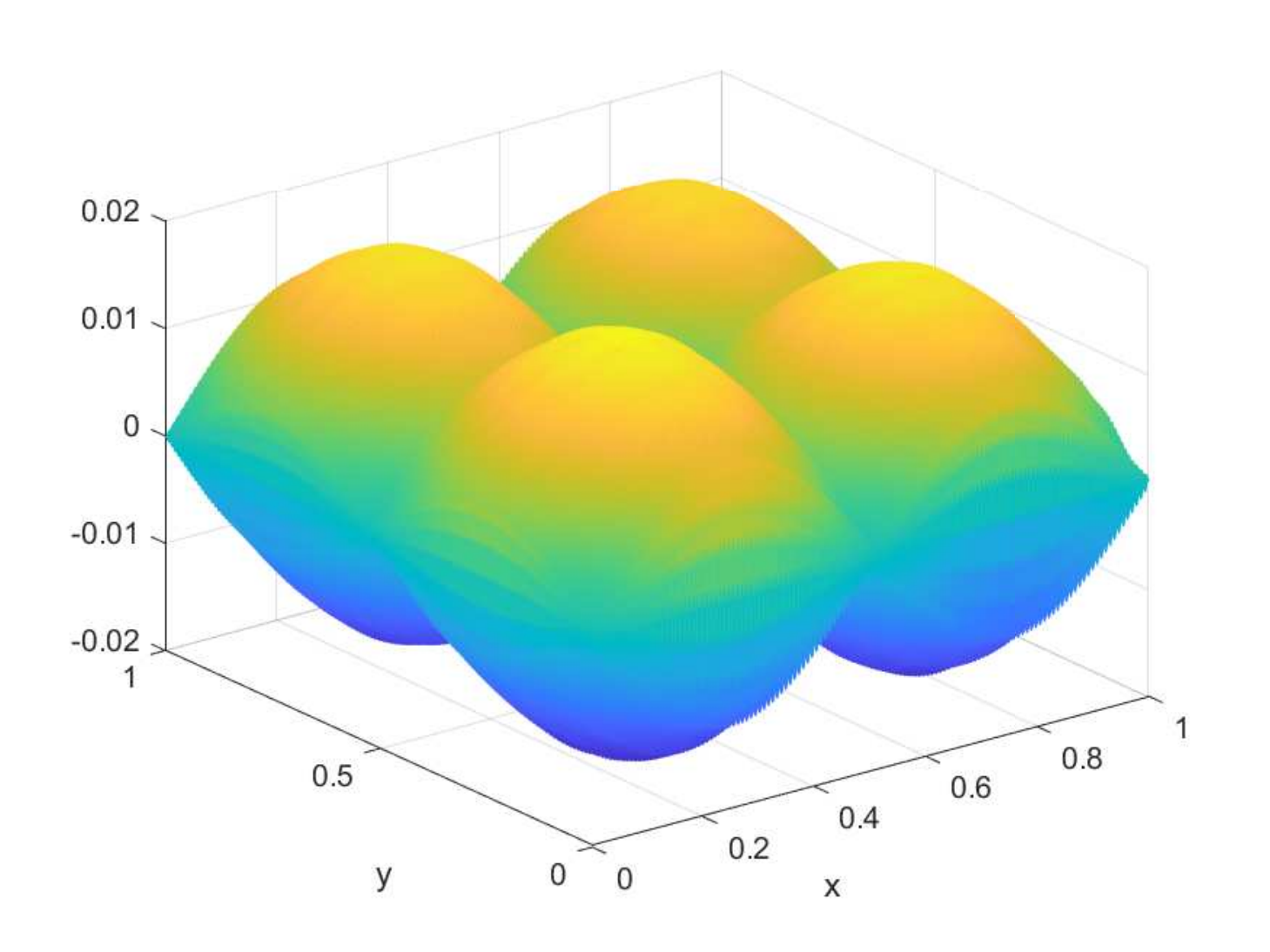}
	\end{subfigure}
	\begin{subfigure}[b]{0.3\textwidth}
		 \includegraphics[width=6cm,height=5cm]{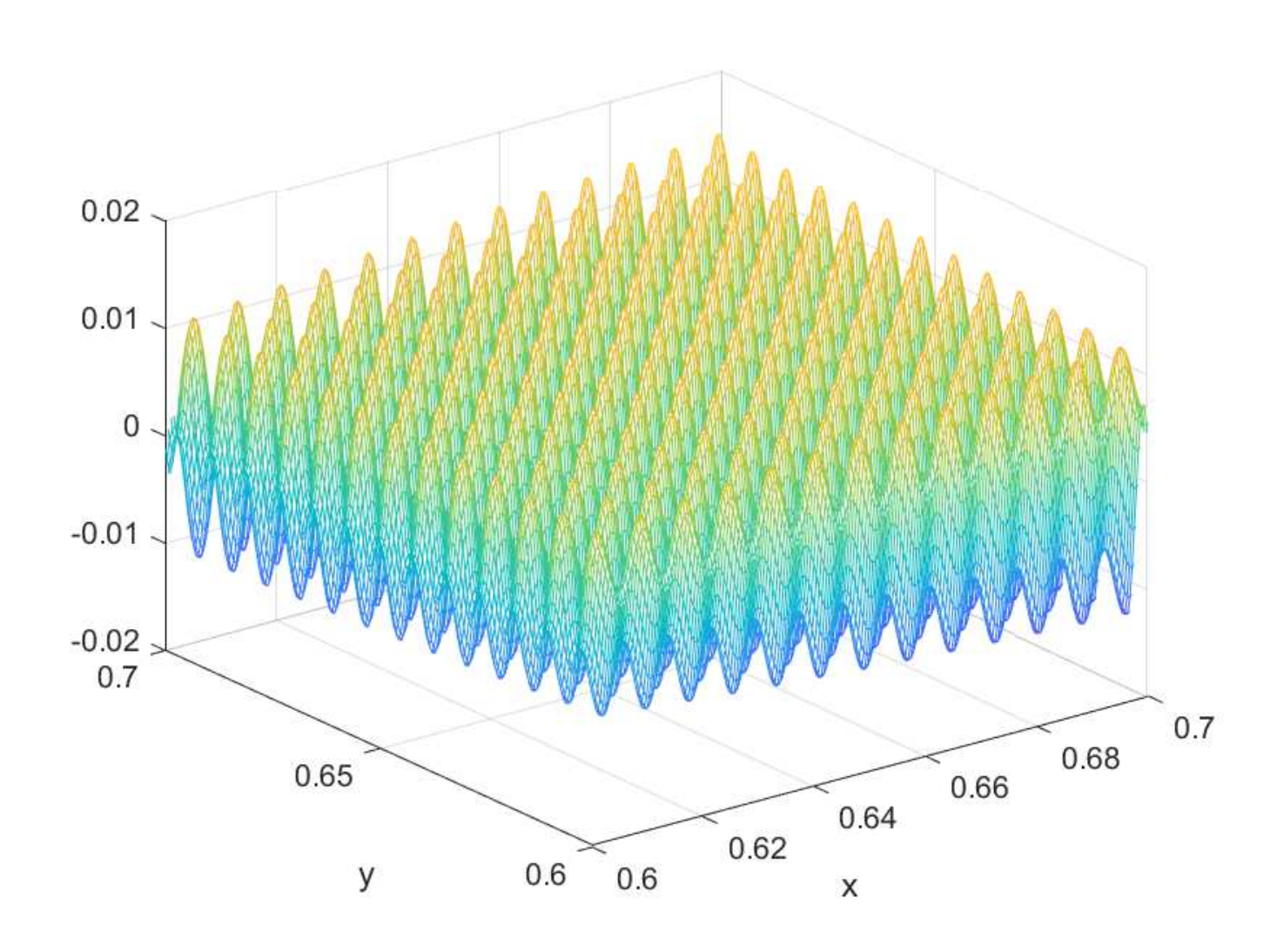}
	\end{subfigure}
	\begin{subfigure}[b]{0.3\textwidth}
		 \includegraphics[width=6cm,height=5cm]{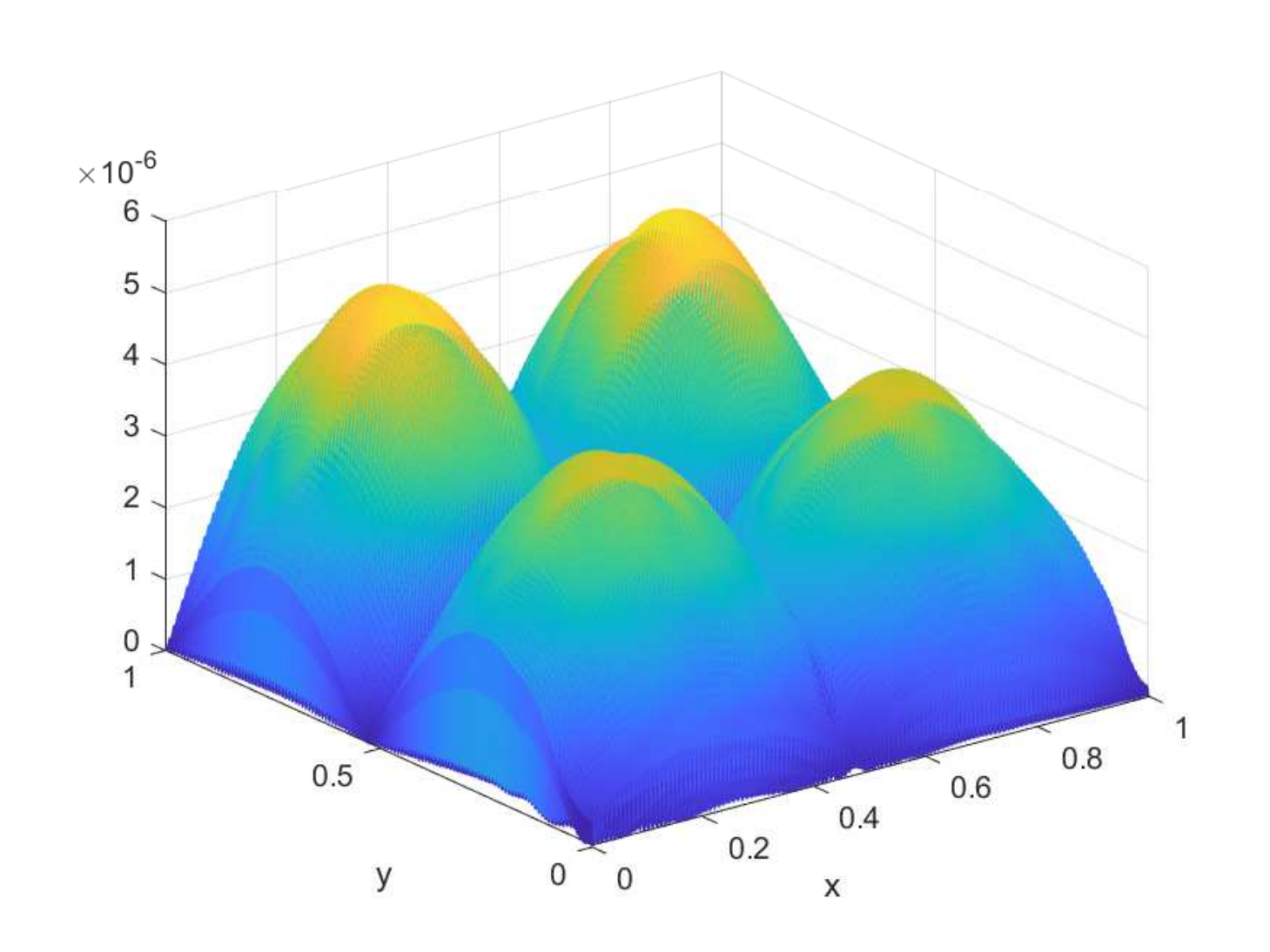}
	\end{subfigure}
	\begin{subfigure}[b]{0.3\textwidth}
		 \includegraphics[width=6cm,height=5cm]{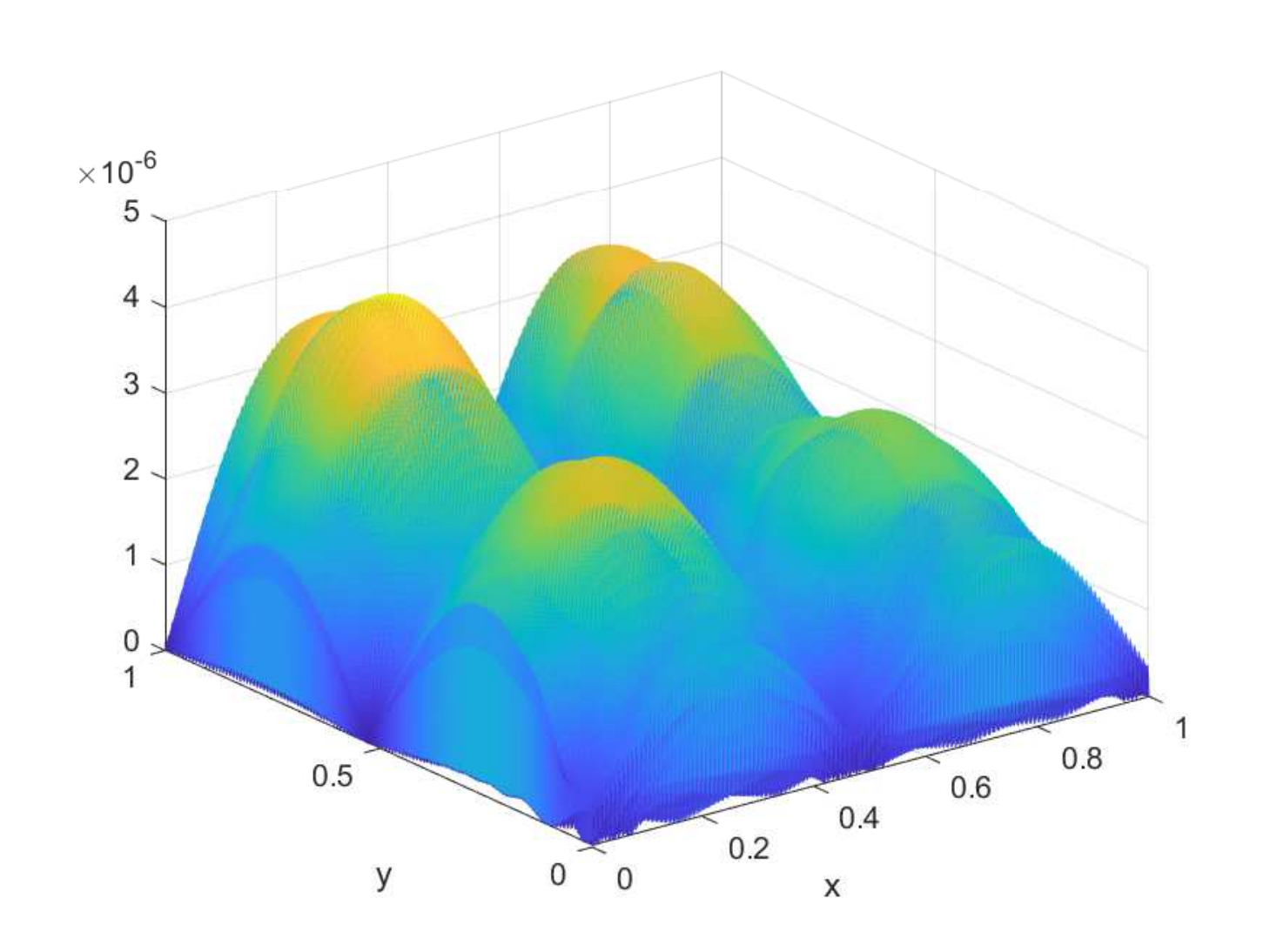}
	\end{subfigure}
	\caption
		{ \cref{ex4} with $\ka=800$ and $h=1/2^{10}$. First row: $\Re(u_{h/2})$ in $[0,1]^2$ (left), $\Re(u_{h/2})$ in $[0.6,0.7]^2$ (middle), $\Im(u_{h/2})$ in $[0,1]^2$ (right). Second row: $\Im(u_{h/2})$ in $[0.6,0.7]^2$ (left), $|\Re(u_h-u_{h/2})|$ in $[0,1]^2$ (middle), $|\Im(u_h-u_{h/2})|$ in $[0,1]^2$ (right).}
	\label{fig:U:graph:k800}
\end{figure}

	We provide four numerical experiments for this case (\cref{ex1,ex2,ex3,ex4}). The first two examples compare our method, denoted by `Proposed', with those proposed in \cite{CGX21, TGGT13, WX18}, denoted by `\cite{CGX21}', `\cite{TGGT13}' and `\cite{WX18}' respectively. Recall that $\frac{2\pi}{\ka h}$ corresponds to the number of points per wavelength. While the first example deals with Dirichlet boundary conditions, the other three examples deal with mixed boundary conditions. We know the true solutions of all examples in this section except for the last one. 
	\begin{example}\label{ex1}
		\normalfont
		Consider the problem \eqref{Qeques2} with no interface curve $\Gamma$ in $\Omega=(0,1)^2$, and
		\begin{align*}	
&  	\ka \in \{50, 150,450\},	\qquad f=0, \qquad u(x,y,\theta)=\exp({\ia}{\ka}(\cos(\theta)x+\sin(\theta)y)),\\			
& u(0,y)=g_1, \quad \mbox{and} \quad  u(1,y)=g_2 \quad \mbox{for} \quad y\in(0,1),\\
& u(x,0)=g_3, \quad \mbox{and} \quad  u(x,1)=g_4  \quad \mbox{for} \quad x\in(0,1),
		\end{align*}
	i.e., we consider
all Dirichlet boundary conditions and the exact solution is the plane wave with the angle $\theta$.
The boundary data $g_1,\ldots,g_4$ are obtained from the above data and the model problem.
We define the following average error for plane wave solutions along all different angles $\theta$ by
	\begin{align*}
		 &\frac{\|u_{h}-u\|_{2,\wa}}{\|u\|_{2,\wa}}:=
		\frac{1}{N_3}\sum_{k=0}^{N_3-1}
		 \sqrt{\frac{\sum_{i=0}^{N_1}\sum_{j=0}^{N_1} \left((u_h)_{i,j,k}-u(x_i,y_j,\theta_k)\right)^2}{\sum_{i=0}^{N_1}\sum_{j=0}^{N_1} \left(u(x_i,y_j,\theta_k)\right)^2}}, \quad \theta_k:=kh_{\theta},  \quad h_{\theta}:=\frac{2\pi}{N_3},\quad  N_3 \in \NN,
	\end{align*}
	where
 $(u_{h})_{i,j,k}$ is the value of the numerical solution $u_h$ at the grid point $(x_i, y_j)$ with a plane wave angle $\theta_k$.
See \cref{table:QSp1} for numerical results.
	\end{example}
	\begin{table}[htbp]
		\caption{Numerical results for \cref{ex1} with $h=1/2^J$. The ratio $\ra$ is equal to  $\frac{\|u_{h}-u\|_{2,\wa}}{\|u\|_{2,\wa}}$ of \cite{CGX21} divided by $\frac{\|u_{h}-u\|_{2,\wa}}{\|u\|_{2,\wa}}$ of our method. In other words, for the same mesh size $h$ with $h=2^{-J}$, the error of \cite{CGX21} is $\ra$ times larger than that of our method.}
		\centering
				\scalebox{0.85}{
		\setlength{\tabcolsep}{0.5mm}{
			 \begin{tabular}{c|c|c|c|c|c|c|c|c|c|c|c|c|c|c|cc}
				\hline
				\multicolumn{1}{c|}{} &
				\multicolumn{5}{c|}{$\ka=50, N_3=30$} &
				\multicolumn{5}{c|}{$\ka=150, N_3=30$} &
				\multicolumn{5}{c}{$\ka=450, N_3=30$} & \\
				\cline{1-17}
				\hline
				\multicolumn{1}{c|}{} &
				 \multicolumn{1}{c|}{\cite{CGX21}} &
				\multicolumn{2}{c|}{Proposed} &
				\multicolumn{1}{c}{} &
				\multicolumn{1}{c|}{} &
				 \multicolumn{1}{c|}{\cite{CGX21}} &
				\multicolumn{2}{c|}{Proposed} &
				\multicolumn{1}{c}{} &
				\multicolumn{1}{c|}{} &
				 \multicolumn{1}{c|}{\cite{CGX21}} &
				\multicolumn{2}{c|}{Proposed} &
				\multicolumn{1}{c}{} &
				\multicolumn{1}{c}{}  \\
				\cline{1-17}
				$J$
				& $\frac{\|u_{h}-u\|_{2,\wa}}{\|u\|_{2,\wa}}$
				&  $\frac{\|u_{h}-u\|_{2,\wa}}{\|u\|_{2,\wa}}$
				& order
				& $\frac{2\pi}{\ka h}$
				& $\ra$
				&  $\frac{\|u_{h}-u\|_{2,\wa}}{\|u\|_{2,\wa}}$
				& $\frac{\|u_{h}-u\|_{2,\wa}}{\|u\|_{2,\wa}}$
				&  order
				& $\frac{2\pi}{\ka h}$
				& $\ra$
				&  $\frac{\|u_{h}-u\|_{2,\wa}}{\|u\|_{2,\wa}}$
				& $\frac{\|u_{h}-u\|_{2,\wa}}{\|u\|_{2,\wa}}$
				&  order
				& $\frac{2\pi}{\ka h}$
				& $\ra$ \\
				\hline
3   &   &   &   &   &   &   &   &   &   &   &   &   &   &   &\\
4   &9.47E+00   &5.33E-01   &   &2.0   &\textbf{17.78}   &   &   &   &   &   &   &   &   &   &\\
5   &1.55E-02   &1.01E-03   &9.0   &4.0   &\textbf{15.35}   &   &   &   &   &   &   &   &   &   &\\
6   &4.97E-05   &1.20E-05   &6.4   &8.0   &\textbf{4.13}   &3.67E+00   &6.25E-02   &   &2.7   &\textbf{58.66}   &   &   &   &   &\\
7   &2.33E-07   &1.77E-07   &6.1   &16.1   &\textbf{1.32}   &6.04E-03   &6.71E-04   &6.5   &5.4   &\textbf{8.99}   &   &   &   &   &\\
8   &   &   &   &   &   &2.56E-05   &9.09E-06   &6.2   &10.7   &\textbf{2.81}   &1.26E+00   &5.40E-02   &   &3.6   &\textbf{23.24}\\
9   &   &   &   &   &   &1.78E-07   &1.37E-07   &6.0   &21.4   &\textbf{1.30}   &4.72E-03   &7.72E-04   &6.1   &7.1   &\textbf{6.11}\\
10   &   &   &   &   &   &2.27E-09   &2.13E-09   &6.0   &42.9   &\textbf{1.06}   &2.25E-05   &1.12E-05   &6.1   &14.3   &\textbf{2.02}\\
11   &   &   &   &   &   &   &   &   &   &   &1.85E-07   &1.71E-07   &6.0   &28.6   &\textbf{1.08}\\
				\hline	
		\end{tabular}}
		\label{table:QSp1}
	}
	\end{table}	
	\begin{example}\label{ex2}
		\normalfont
		Consider the problem \eqref{Qeques2}  with no interface curve $\Gamma$ in $\Omega=(0,1)^2$, and
		\begin{align*}
			& \ka=300,\qquad 	 u=(y-1)\cos(\alpha x)\sin(\beta(y-1)),  \quad \alpha,\beta\in \R,\\
			& u(0,y)=g_1, \quad \mbox{and} \quad  u(1,y)=g_2 \quad \mbox{for} \quad y\in(0,1),\\
			& u(x,0)=g_3, \quad \mbox{and} \quad  u_y(x,1)-\ia\ka u(x,1)=g_4  \quad \mbox{for} \quad x\in(0,1),	
		\end{align*}
		where the boundary data $g_1,\ldots,g_4$ and the source term $f$  are obtained from the above data and the model problem.
See \cref{table:QSp2,table:QSp3} for numerical results for various choices of $\alpha$ and $\beta$. Note that the errors of the methods in \cite{TGGT13, WX18} are at least twice as large as ours.
	\end{example}
	\begin{table}[htbp]
		\caption{Numerical results of \cref{ex2} with $h=1/2^J$ and $\ka=300$. The ratio $\ra_1$ is equal to  $\|u_{h}-u\|_{\infty}$ of \cite{TGGT13} divided by $\|u_{h}-u\|_{\infty}$ of our method and the ratio $\ra_2$ is equal to  $\|u_{h}-u\|_{\infty}$ of \cite{WX18} divided by $\|u_{h}-u\|_{\infty}$ of our method. In other words, for the same grid size $h$ with $h=2^{-J}$, the errors of \cite{TGGT13} and \cite{WX18}
are $\ra_1$ and $\ra_2$ times larger than those of our method, respectively.
}
	\scalebox{0.73}{
		\centering
		\setlength{\tabcolsep}{0.01mm}{
			 \begin{tabular}{c|c|c|c|c|c|c|c|c|c|c|c|c|c|c|c|c|c|c|c}
				\hline
				\multicolumn{1}{c}{} &
				\multicolumn{1}{c|}{} &
				 \multicolumn{6}{c|}{$\alpha=50$, $\beta=290$} &
				 \multicolumn{6}{c|}{$\alpha=100$, $\beta=275$} &
				 \multicolumn{6}{c}{$\alpha=150$, $\beta=255$} \\					
				\cline{1-17}
				\hline
				\multicolumn{1}{c}{} &
				\multicolumn{1}{c|}{} &
				 \multicolumn{1}{c|}{\cite{TGGT13}} &
				 \multicolumn{1}{c|}{\cite{WX18}} &
				 \multicolumn{2}{c|}{\tiny{Proposed}} &
				\multicolumn{2}{c|}{} &
				 \multicolumn{1}{c|}{\cite{TGGT13}} &
				 \multicolumn{1}{c|}{\cite{WX18}} &
				 \multicolumn{2}{c|}{\tiny{Proposed}} &
				\multicolumn{2}{c|}{} &		
				 \multicolumn{1}{c|}{\cite{TGGT13}} &
				 \multicolumn{1}{c|}{\cite{WX18}} &
				 \multicolumn{2}{c|}{\tiny{Proposed}} &
				\multicolumn{2}{c}{}\\		
				\cline{1-17}
				\hline
				\multicolumn{1}{c|}{$J$} &
				 \multicolumn{1}{c|}{$\frac{2\pi}{\ka h}$} &
				 \multicolumn{1}{c|}{$\|u_{h}-u\|_{\infty}$} &
				 \multicolumn{1}{c|}{$\|u_{h}-u\|_{\infty}$} &
				 \multicolumn{1}{c|}{$\|u_{h}-u\|_{\infty}$} &
				 \multicolumn{1}{c|}{\tiny{order}} &
				\multicolumn{1}{c|}{$\ra_1$} &
				\multicolumn{1}{c|}{$\ra_2$} &
				 \multicolumn{1}{c|}{$\|u_{h}-u\|_{\infty}$} &
				 \multicolumn{1}{c|}{$\|u_{h}-u\|_{\infty}$} &
				 \multicolumn{1}{c|}{$\|u_{h}-u\|_{\infty}$} &
				  \multicolumn{1}{c|}{\tiny{order}} &
				\multicolumn{1}{c|}{$\ra_1$} &
				\multicolumn{1}{c|}{$\ra_2$} &
				 \multicolumn{1}{c|}{$\|u_{h}-u\|_{\infty}$} &
				 \multicolumn{1}{c|}{$\|u_{h}-u\|_{\infty}$} &
				 \multicolumn{1}{c|}{$\|u_{h}-u\|_{\infty}$} &
				  \multicolumn{1}{c|}{\tiny{order}} &
				\multicolumn{1}{c|}{$\ra_1$} &
				\multicolumn{1}{c}{$\ra_2$} \\	
				\cline{1-17}
				\hline
7   &2.7   &1.17E+00   &7.13E-02   &2.71E-02   &   &\textbf{43}   &\textbf{2.6}   &1.37E+00   &1.43E-01   &3.02E-02   &   &\textbf{45}   &\textbf{4.7}   &2.71E+00   &1.07E-01   &4.93E-02   &   &\textbf{55}   &\textbf{2.2}\\
8   &5.4   &6.09E-03   &4.39E-04   &8.81E-05   &8.3   &\textbf{69}   &\textbf{5.0}   &8.72E-03   &9.37E-04   &2.17E-04   &7.1   &\textbf{40}   &\textbf{4.3}   &1.51E-02   &7.39E-04   &8.39E-05   &9.2   &\textbf{180}   &\textbf{8.8}\\
9   &10.7   &8.69E-05   &5.99E-06   &1.90E-06   &5.5   &\textbf{46}   &\textbf{3.2}   &1.24E-04   &1.28E-05   &4.04E-06   &5.7   &\textbf{31}   &\textbf{3.2}   &2.22E-04   &1.08E-05   &1.42E-06   &5.9   &\textbf{156}   &\textbf{7.6}\\
10   &21.4   &1.32E-06   &8.59E-08   &3.12E-08   &5.9   &\textbf{42}   &\textbf{2.8}   &1.89E-06   &1.88E-07   &6.50E-08   &6.0   &\textbf{29}   &\textbf{2.9}   &3.39E-06   &1.57E-07   &2.42E-08   &5.9   &\textbf{140}   &\textbf{6.5}\\
11   &42.9   &2.07E-08   &1.32E-09   &4.96E-10   &6.0   &\textbf{42}   &\textbf{2.7}   &2.94E-08   &2.89E-09   &1.03E-09   &6.0   &\textbf{29}   &\textbf{2.8}   &5.27E-08   &2.41E-09   &3.89E-10   &6.0   &\textbf{135}   &\textbf{6.2}\\
				\hline
		\end{tabular}}
		\label{table:QSp2}}
	\end{table}	
	%
	\begin{table}[htbp]
		\caption{Numerical results of \cref{ex2} with $h=1/2^J$ and $\ka=300$. The ratio $\ra_1$ is equal to  $\|u_{h}-u\|_{\infty}$ of \cite{TGGT13} divided by $\|u_{h}-u\|_{\infty}$ of our method and the ratio $\ra_2$ is equal to  $\|u_{h}-u\|_{\infty}$ of \cite{WX18} divided by $\|u_{h}-u\|_{\infty}$ of our method. I.e., for the same grid size $h$ with $h=2^{-J}$, the errors of \cite{TGGT13} and \cite{WX18}
			are $\ra_1$ and $\ra_2$ times larger than those of our method, respectively.}
		\centering
			\scalebox{0.7}{
		\setlength{\tabcolsep}{0.01mm}{
			 \begin{tabular}{c|c|c|c|c|c|c|c|c|c|c|c|c|c|c|c|c|c|c|c}
				\hline
				\multicolumn{1}{c}{} &
				\multicolumn{1}{c|}{} &
				 \multicolumn{6}{c|}{$\alpha=200$, $\beta=200$} &
				 \multicolumn{6}{c|}{$\alpha=250$, $\beta=160$} &
				 \multicolumn{6}{c}{$\alpha=290$, $\beta=50$} \\					
				\cline{1-17}
			\hline
			\multicolumn{1}{c}{} &
			\multicolumn{1}{c|}{} &
			 \multicolumn{1}{c|}{\cite{TGGT13}} &
			\multicolumn{1}{c|}{\cite{WX18}} &
			 \multicolumn{2}{c|}{\tiny{Proposed}} &
			\multicolumn{2}{c|}{} &
			 \multicolumn{1}{c|}{\cite{TGGT13}} &
			\multicolumn{1}{c|}{\cite{WX18}} &
			 \multicolumn{2}{c|}{\tiny{Proposed}} &
			\multicolumn{2}{c|}{} &		
			 \multicolumn{1}{c|}{\cite{TGGT13}} &
			\multicolumn{1}{c|}{\cite{WX18}} &
			 \multicolumn{2}{c|}{\tiny{Proposed}} &
			\multicolumn{2}{c}{}\\		
			\cline{1-17}
			\hline
			\multicolumn{1}{c|}{$J$} &
			 \multicolumn{1}{c|}{$\frac{2\pi}{\ka h}$} &
			 \multicolumn{1}{c|}{$\|u_{h}-u\|_{\infty}$} &
			 \multicolumn{1}{c|}{$\|u_{h}-u\|_{\infty}$} &
			 \multicolumn{1}{c|}{$\|u_{h}-u\|_{\infty}$} &
			\multicolumn{1}{c|}{\tiny{order}} &
			\multicolumn{1}{c|}{$\ra_1$} &
			\multicolumn{1}{c|}{$\ra_2$} &
			 \multicolumn{1}{c|}{$\|u_{h}-u\|_{\infty}$} &
			 \multicolumn{1}{c|}{$\|u_{h}-u\|_{\infty}$} &
			 \multicolumn{1}{c|}{$\|u_{h}-u\|_{\infty}$} &
			\multicolumn{1}{c|}{\tiny{order}} &
			\multicolumn{1}{c|}{$\ra_1$} &
			\multicolumn{1}{c|}{$\ra_2$} &
			 \multicolumn{1}{c|}{$\|u_{h}-u\|_{\infty}$} &
			 \multicolumn{1}{c|}{$\|u_{h}-u\|_{\infty}$} &
			 \multicolumn{1}{c|}{$\|u_{h}-u\|_{\infty}$} &
			\multicolumn{1}{c|}{\tiny{order}} &
			\multicolumn{1}{c|}{$\ra_1$} &
			\multicolumn{1}{c}{$\ra_2$} \\	
			\cline{1-17}
			\hline
7   &2.7   &8.19E-01   &1.03E-01   &1.18E-01   &   &\textbf{7}   &\textbf{0.9}   &3.32E+00   &1.21E-01   &4.45E-02   &   &\textbf{74}   &\textbf{2.7}   &5.12E+00   &1.00E-01   &4.31E-02   &   &\textbf{119}   &\textbf{2.3}\\
8   &5.4   &6.05E-03   &7.81E-04   &3.08E-04   &8.6   &\textbf{20}   &\textbf{2.5}   &3.26E-02   &9.33E-04   &4.78E-05   &9.9   &\textbf{682}   &\textbf{19.5}   &8.03E-03   &5.73E-04   &1.33E-04   &8.3   &\textbf{61}   &\textbf{4.3}\\
9   &10.7   &9.31E-05   &1.10E-05   &2.90E-06   &6.7   &\textbf{32}   &\textbf{3.8}   &4.84E-04   &1.42E-05   &1.58E-06   &4.9   &\textbf{306}   &\textbf{9.0}   &1.17E-04   &7.98E-06   &1.49E-06   &6.5   &\textbf{79}   &\textbf{5.4}\\
10   &21.4   &1.46E-06   &1.66E-07   &3.92E-08   &6.2   &\textbf{37}   &\textbf{4.2}   &7.51E-06   &2.04E-07   &2.75E-08   &5.8   &\textbf{273}   &\textbf{7.4}   &1.79E-06   &1.14E-07   &2.12E-08   &6.1   &\textbf{85}   &\textbf{5.4}\\
11   &42.9   &2.27E-08   &2.57E-09   &5.86E-10   &6.1   &\textbf{39}   &\textbf{4.4}   &1.18E-07   &3.14E-09   &4.48E-10   &5.9   &\textbf{262}   &\textbf{7.0}   &2.81E-08   &1.76E-09   &3.23E-10   &6.0   &\textbf{87}   &\textbf{5.4} \\
				\hline				
		\end{tabular}}
		\label{table:QSp3}}
	\end{table}		
	\begin{example}\label{ex3}
		\normalfont
		Consider the problem \eqref{Qeques2}  with no interface curve $\Gamma$ in $\Omega=(0,1)^2$, and
	\[
	\begin{aligned}
		& \ka \in \{450, 650\},\qquad 	 u=\sin(\alpha x+\beta y), \quad \alpha,\ \beta\in \R,\\
			& -u_x(0,y)-\ia\ka u(0,y)= g_1 \quad \mbox{and}
\quad u(1,y)= g_2 \quad \mbox{for} \quad y\in(0,1),\\
& -u_y(x,0)= g_3 \quad \mbox{and}
\quad u_y(x,1)-\ia\ka u(x,1)= g_4 \quad \mbox{for} \quad x\in(0,1),
\end{aligned}
	\]
		where the boundary data $g_1,\ldots,g_4$ and the source term $f$ are obtained from the above data and the model problem.
		See \cref{table:QSp4} for numerical results for various choices of $\alpha$ and $\beta$.	
	\end{example}
	\begin{table}[htbp]
		\caption{Numerical results of \cref{ex3} with $h=1/2^J$ using our method.}
		\centering
		\setlength{\tabcolsep}{0.5mm}{
			 \begin{tabular}{c|c|c|c|c|c|c|c|c|c|c|c}
				\hline
				 \multicolumn{6}{c|}{$\ka=450$, $\alpha=400$, $\beta=200$} &
				\multicolumn{6}{c}{$\ka=650$, $\alpha=250$, $\beta=600$} \\				
				\cline{1-12}
				\hline
				$J$&    $\frac{2\pi}{\ka h}$  &  $\frac{\|u_{h}-u\|_{2}}{\|u\|_{2}}$    &order &   ${\|u_{h}-u\|_{\infty}}$    &order &   $J$ & $\frac{2\pi}{\ka h}$ &  $\frac{\|u_{h}-u\|_{2}}{\|u\|_{2}}$    &order &   ${\|u_{h}-u\|_{\infty}}$     &order  \\
				\hline
8   &3.57   &1.6912E-02   &   &2.9616E-02   &   &8   &2.47   &6.0301E-01   &   &9.5806E-01   &\\
9   &7.15   &1.6013E-04   &6.7   &2.4755E-04   &6.9   &9   &4.95   &3.9578E-03   &7.3   &6.8610E-03   &7.1\\
10   &14.30   &2.3644E-06   &6.1   &3.8461E-06   &6.0   &10   &9.90   &4.9900E-05   &6.3   &8.6360E-05   &6.3\\
11   &28.60   &3.7478E-08   &6.0   &6.3435E-08   &5.9   &11   &19.80   &7.3859E-07   &6.1   &1.2928E-06   &6.1\\
				\hline
		\end{tabular}}
		\label{table:QSp4}
	\end{table}	
	\begin{example}\label{ex4}
		\normalfont
		Consider the problem \eqref{Qeques2}  with no interface curve $\Gamma$ in $\Omega=(0,1)^2$, and
		\begin{align*}
			& \ka \in \{200,400,800\}, \qquad f(x,y)=\ka^2\sin(2\pi x)\sin(2\pi y), \\
			& -u_x(0,y)-\ia\ka u(0,y)= \sin(\pi y) \quad \mbox{and}
			\quad u(1,y)= 0 \quad \mbox{for} \quad y\in(0,1),\\
			& -u_y(x,0)= \sin(\pi x) \quad \mbox{and}
			\quad u_y(x,1)-\ia\ka u(x,1)= \sin(\pi x) \quad \mbox{for} \quad x\in(0,1).
		\end{align*}
		The exact solution $u$ is unknown in this example. See \cref{table:QSp5} and \cref{fig:U:graph:k200,fig:U:graph:k400,fig:U:graph:k800} for numerical results.	
	\end{example}	
	\begin{table}[htbp]
		\caption{Numerical results of \cref{ex4} with $h=1/2^J$ using our method.}
		\centering
			\scalebox{0.77}{
		\setlength{\tabcolsep}{0.0001mm}{
			 \begin{tabular}{c|c|c|c|c|c|c|c|c|c|c|c|c|c|c|c|c|c}
				\hline
				 \multicolumn{6}{c|}{$\ka=200$} &
				 \multicolumn{6}{c|}{$\ka=400$} &
				\multicolumn{6}{c}{ $\ka=800$ }
				\\				
				\cline{1-16}
				\hline
				$J$&    $\frac{2\pi}{\ka h}$  &  $\|u_{h}-u_{h/2}\|_{2}$ &order &   $\|u_{h}-u_{h/2}\|_{\infty}$ &order
				&  $J$&  $\frac{2\pi}{\ka h}$  &  $\|u_{h}-u_{h/2}\|_{2}$    &order &   $\|u_{h}-u_{h/2}\|_{\infty}$ &order
				&  $J$&  $\frac{2\pi}{\ka h}$  &  $\|u_{h}-u_{h/2}\|_{2}$    &order &   $\|u_{h}-u_{h/2}\|_{\infty}$ &order  \\
				\hline
4   &0.50   &8.121E+01   &   &1.616E+02   &   &5   &0.50   &8.307E+01   &   &1.661E+02   &   &6   &0.50   &8.360E+01   &   &1.672E+02   & \\
5   &1.01   &1.955E+00   &5.4   &3.899E+00   &5.4   &6   &1.01   &1.874E+00   &5.5   &3.746E+00   &5.5   &7   &1.01   &1.855E+00   &5.5   &3.709E+00   &5.5 \\
6   &2.01   &2.653E-02   &6.2   &6.984E-02   &5.8   &7   &2.01   &1.935E-02   &6.6   &4.422E-02   &6.4   &8   &2.01   &1.239E-02   &7.2   &3.033E-02   &6.9 \\
7   &4.02   &1.449E-04   &7.5   &3.333E-04   &7.7   &8   &4.02   &1.805E-04   &6.7   &4.443E-04   &6.6   &9   &4.02   &1.793E-04   &6.1   &4.328E-04   &6.1 \\
8   &8.04   &1.731E-06   &6.4   &4.034E-06   &6.4   &9   &8.04   &2.153E-06   &6.4   &5.468E-06   &6.3   &10   &8.04   &2.055E-06   &6.4   &5.190E-06   &6.4 \\
				\hline
		\end{tabular}}
		\label{table:QSp5}}
	\end{table}	
%
%
%
%
%
%
	%
	%
	%
	\subsection{Numerical examples with interfaces}
\begin{figure}[htbp]
	\centering
	\begin{subfigure}[b]{0.3\textwidth}
		 \includegraphics[width=5.8cm,height=5cm]{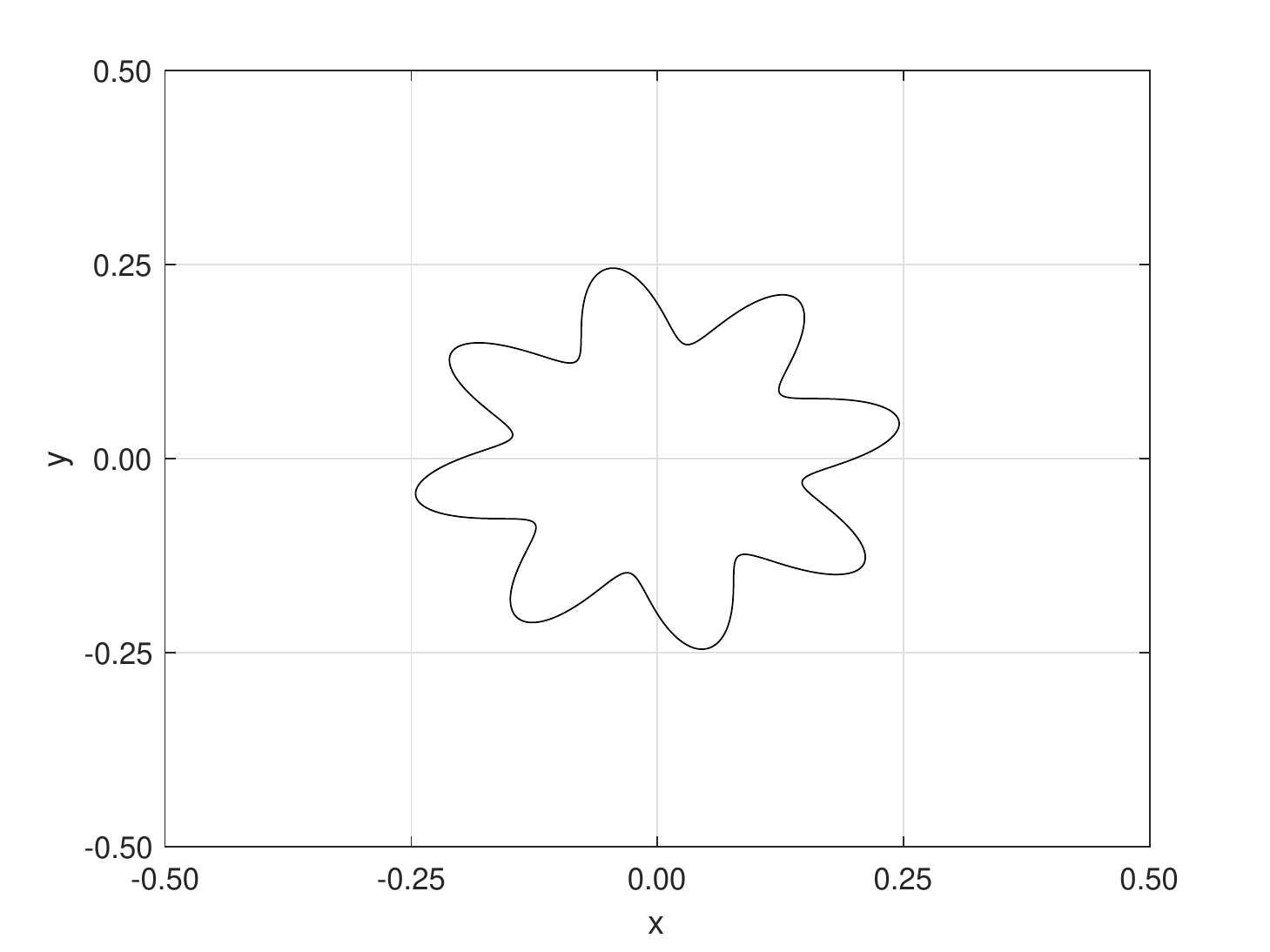}
	\end{subfigure}
	\begin{subfigure}[b]{0.3\textwidth}
		 \includegraphics[width=6cm,height=5cm]{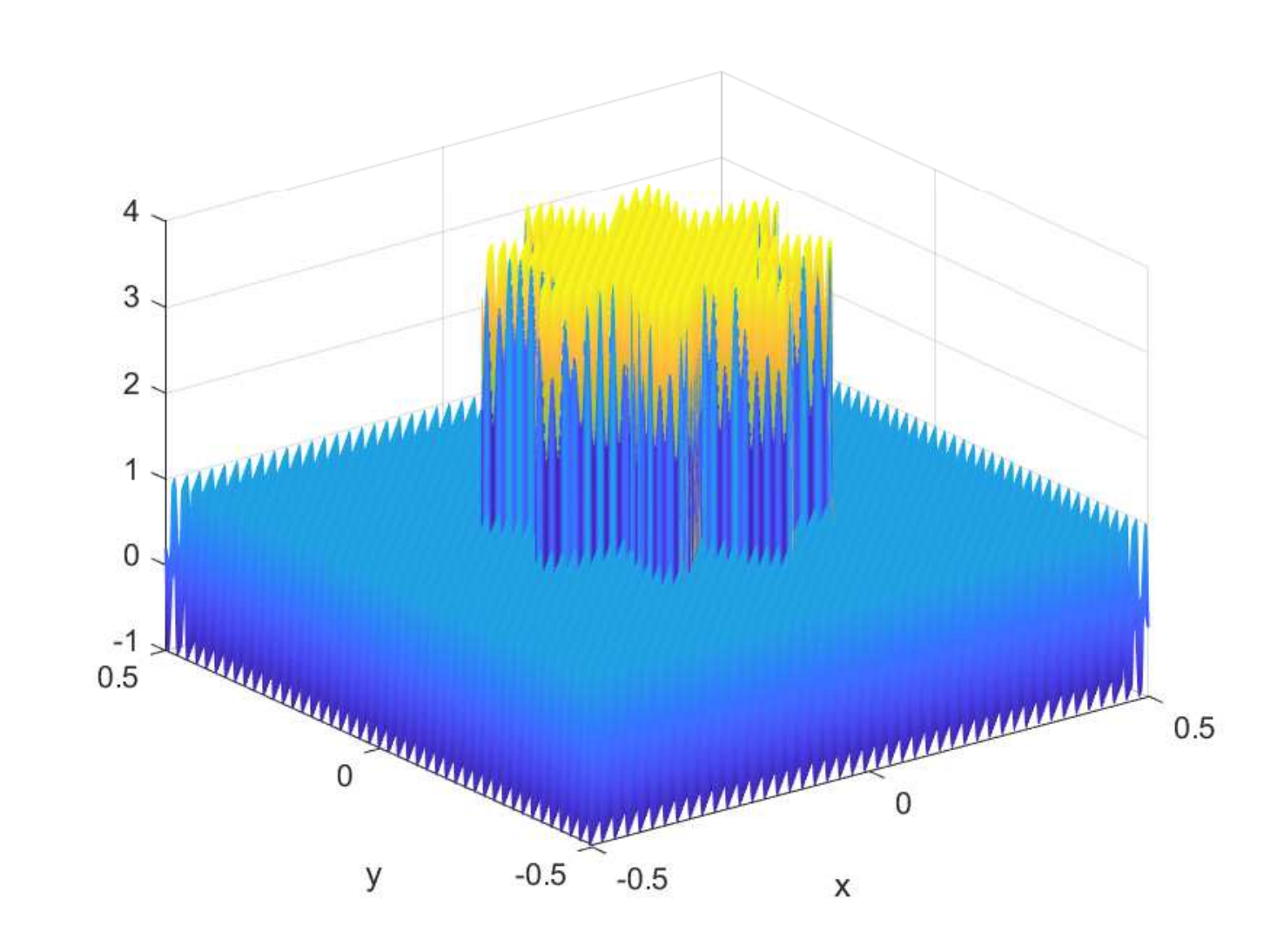}
	\end{subfigure}
	\begin{subfigure}[b]{0.3\textwidth}
		 \includegraphics[width=6cm,height=5cm]{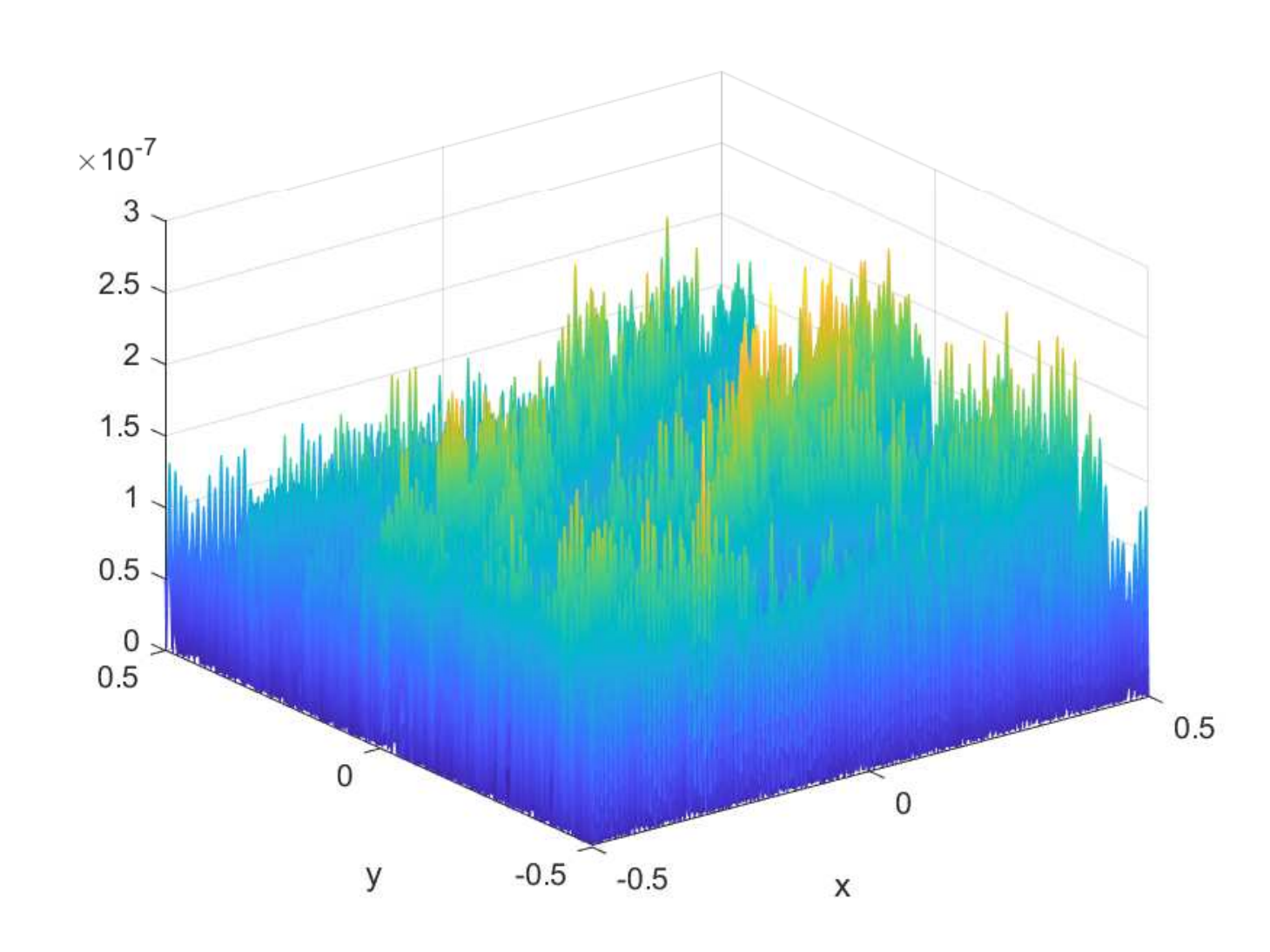}
	\end{subfigure}
	\caption
	{ \cref{interface:ex1}: The interface curve $\Gamma$ (left), $\Re(u_h)$ with $h=1/2^{11}$ (middle), and $|\Re(u_h-u)|$ with $h=1/2^{11}$ (right).}
	\label{fig:interface:1}
\end{figure}

\begin{figure}[htbp]
		\centering
		\begin{subfigure}[b]{0.3\textwidth}
			 \includegraphics[width=5.8cm,height=5cm]{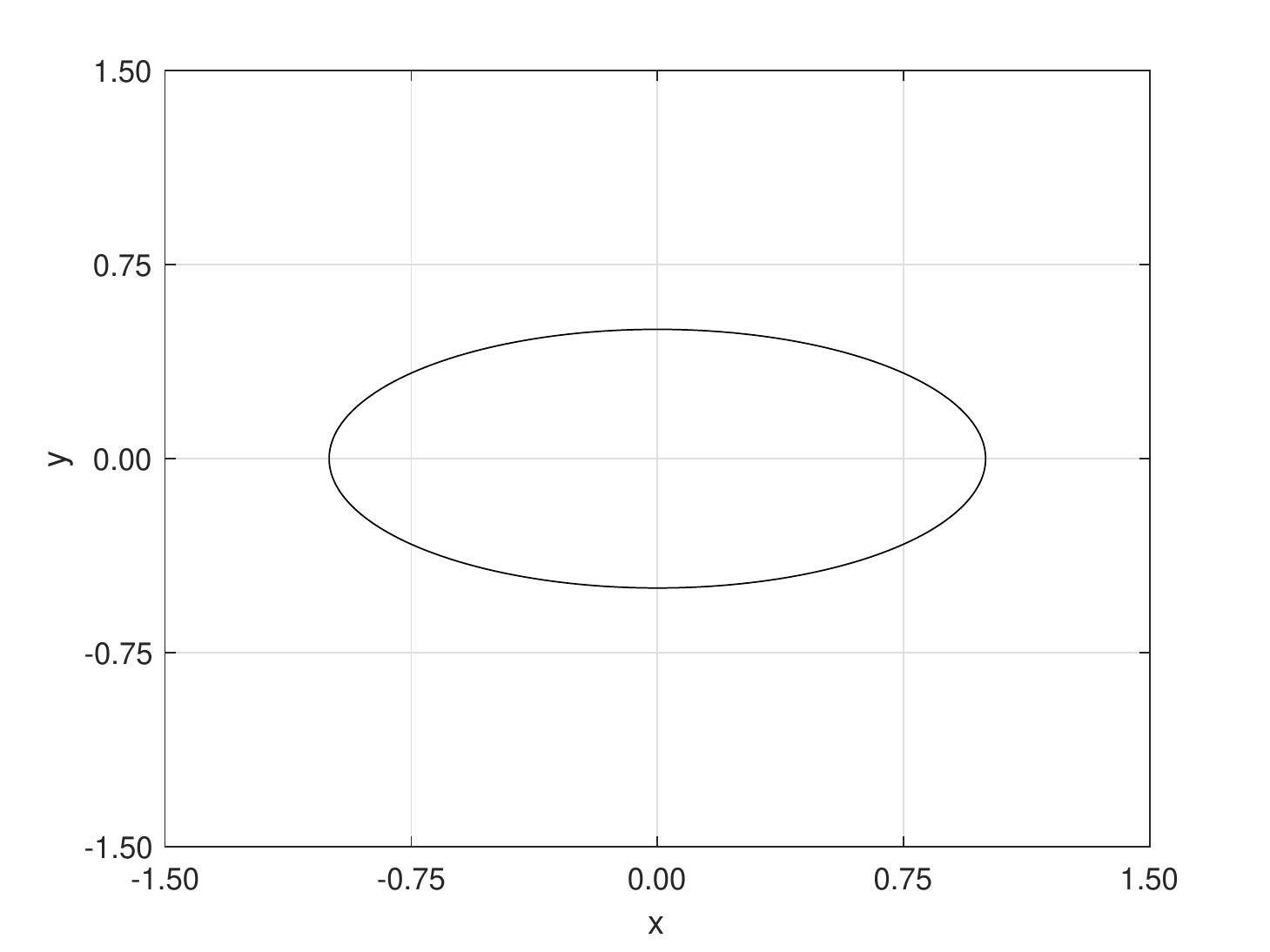}
		\end{subfigure}
		\begin{subfigure}[b]{0.3\textwidth}
			 \includegraphics[width=6cm,height=5cm]{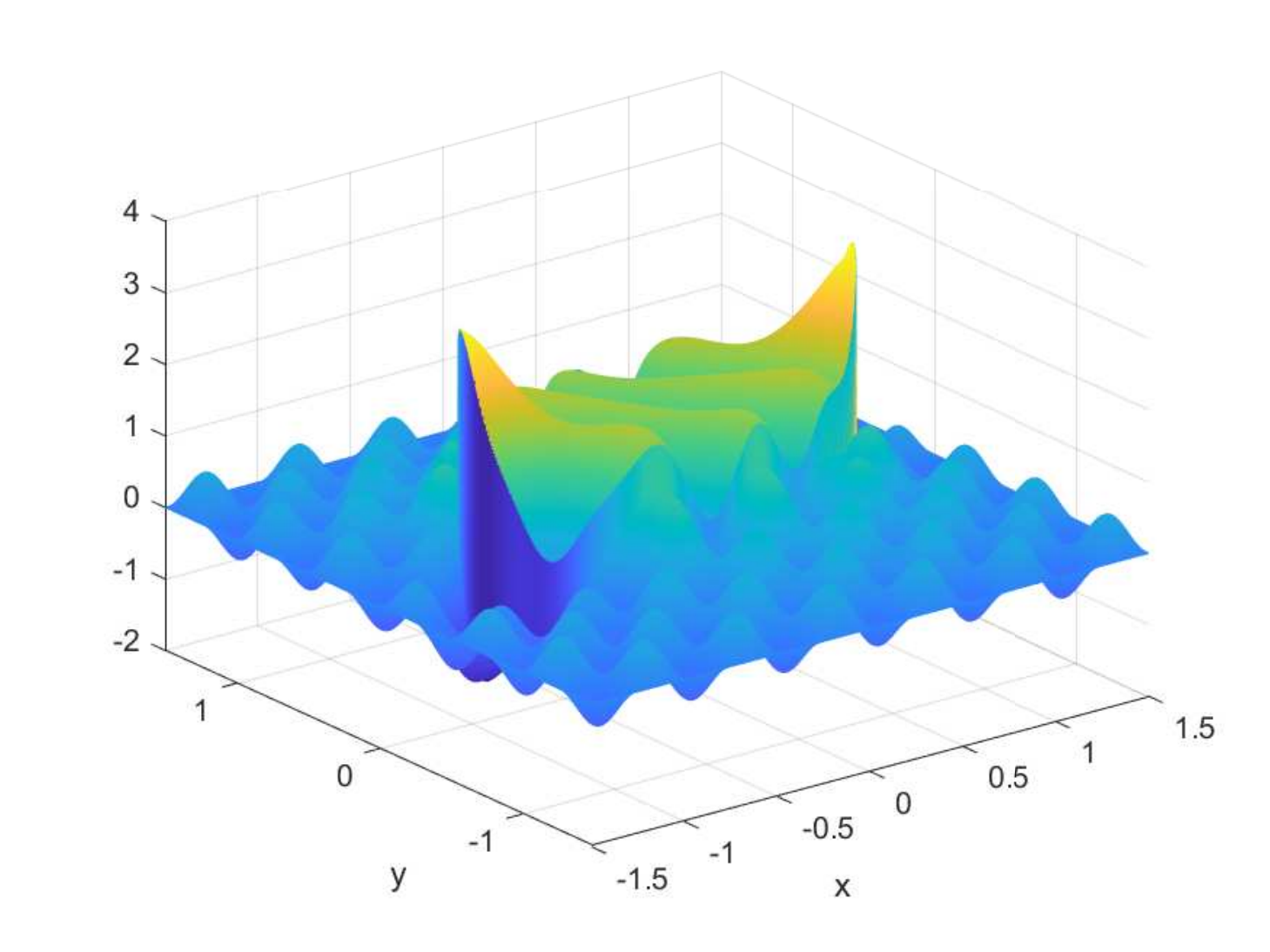}
		\end{subfigure}
		\begin{subfigure}[b]{0.3\textwidth}
			 \includegraphics[width=6cm,height=5cm]{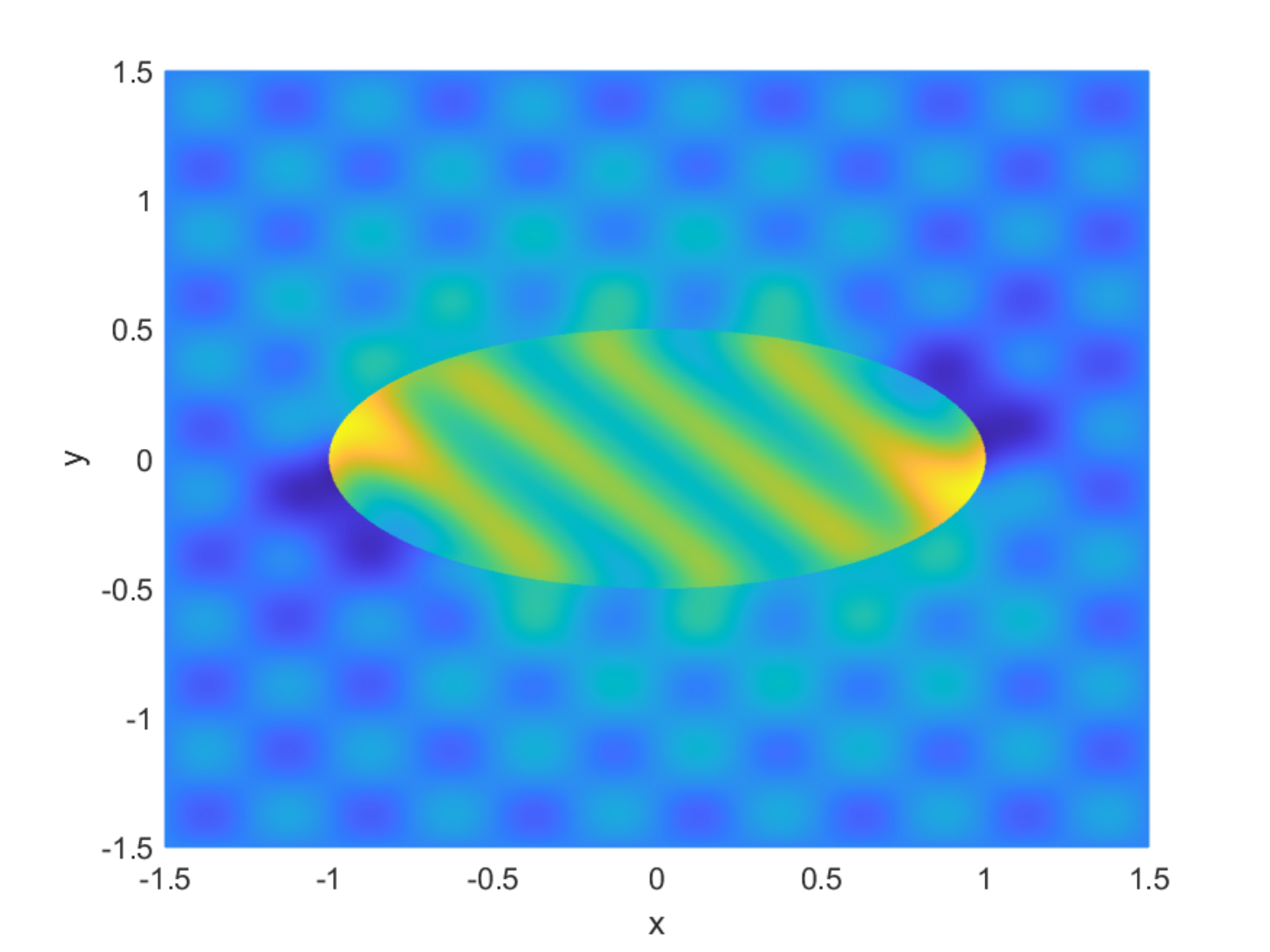}
		\end{subfigure}
	  \begin{subfigure}[b]{0.35\textwidth}
		\includegraphics[width=6cm,height=5cm]{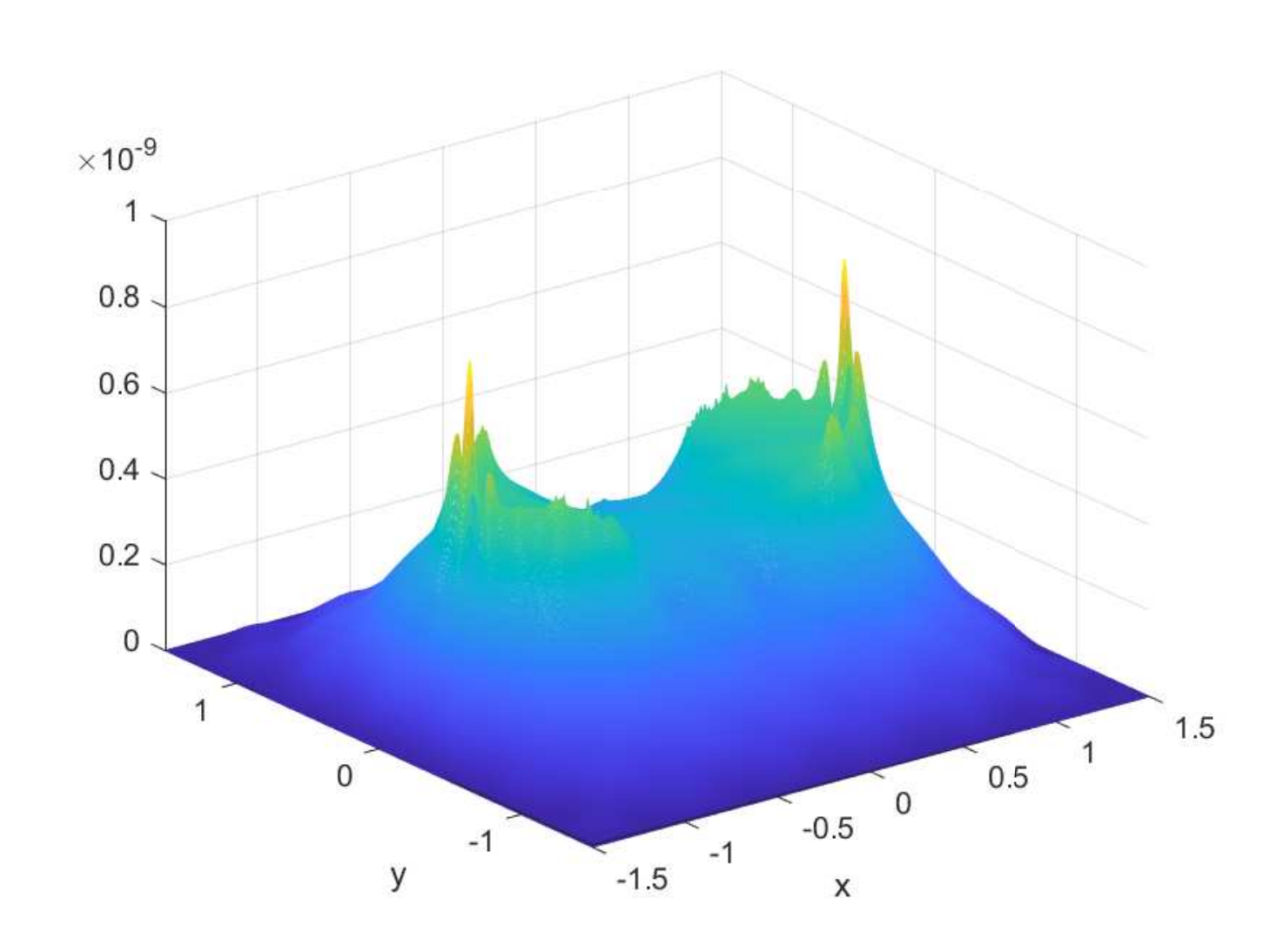}
	\end{subfigure}
	\begin{subfigure}[b]{0.35\textwidth}
		\includegraphics[width=6cm,height=5cm]{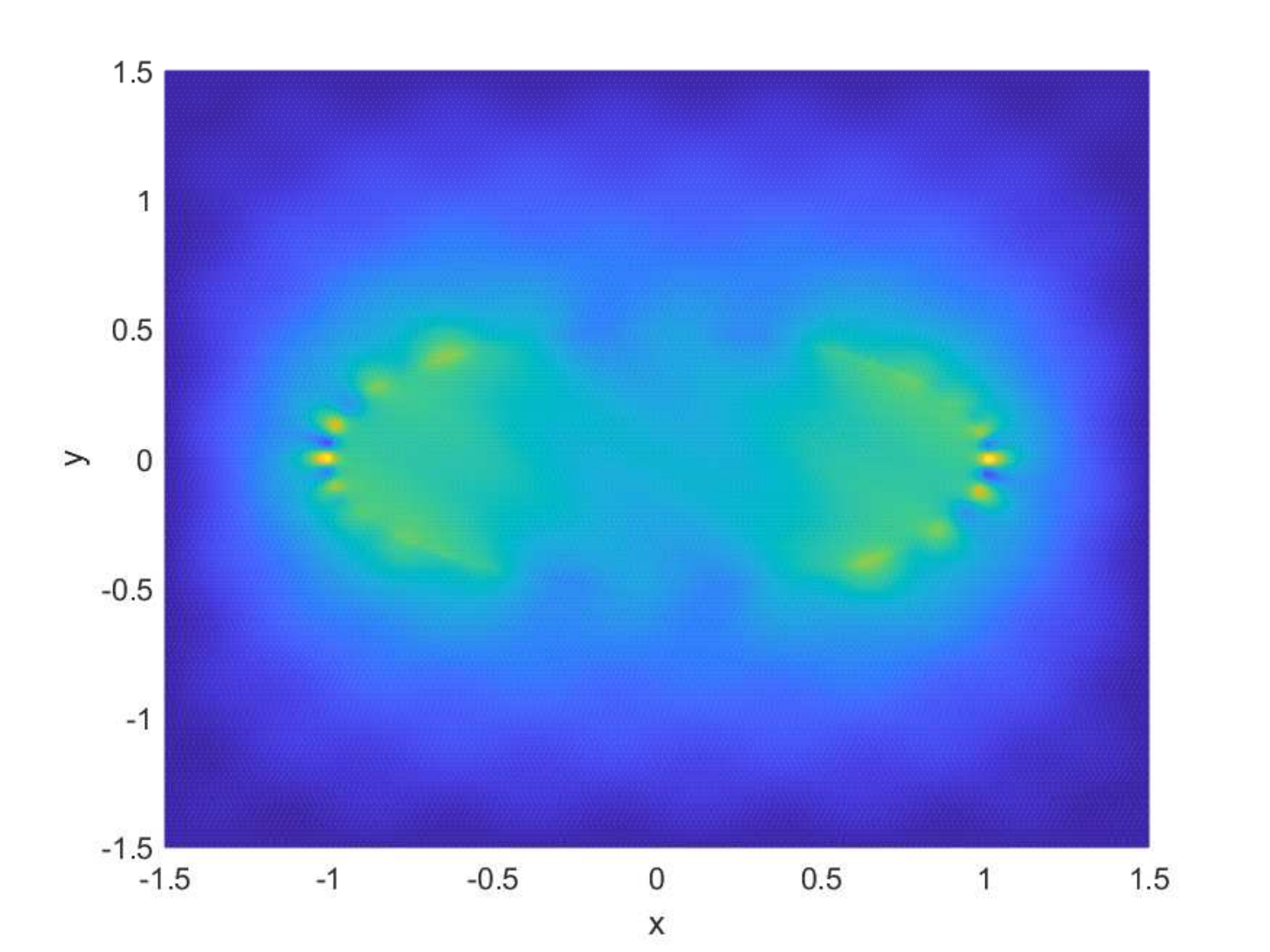}
	\end{subfigure}
		\caption
		{ \cref{interface:ex2} with $\ka_{+}=\ka_{-}=0$. First row: The interface curve $\Gamma$ (left), $u_h$ with  $h=3/2^{10}$ (middle and right).  Second row: $|u_h-u_{h/2}|$ with  $h=3/2^{9}$ (left and right).}
		\label{fig:curvature}
	\end{figure}

\begin{figure}[htbp]
	\centering
	\begin{subfigure}[b]{0.3\textwidth}
		 \includegraphics[width=5.8cm,height=5cm]{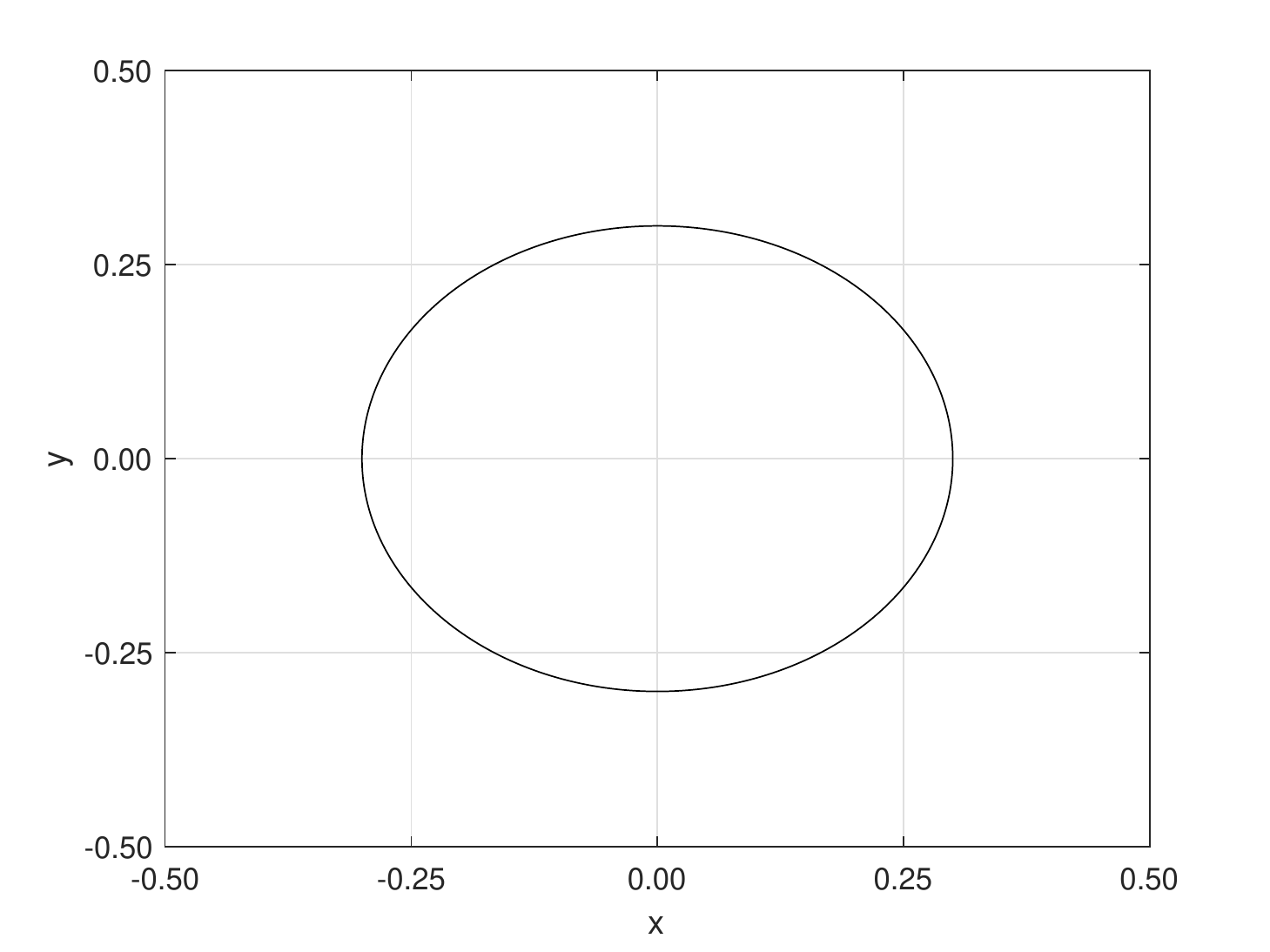}
	\end{subfigure}
	\begin{subfigure}[b]{0.3\textwidth}
		 \includegraphics[width=6cm,height=5cm]{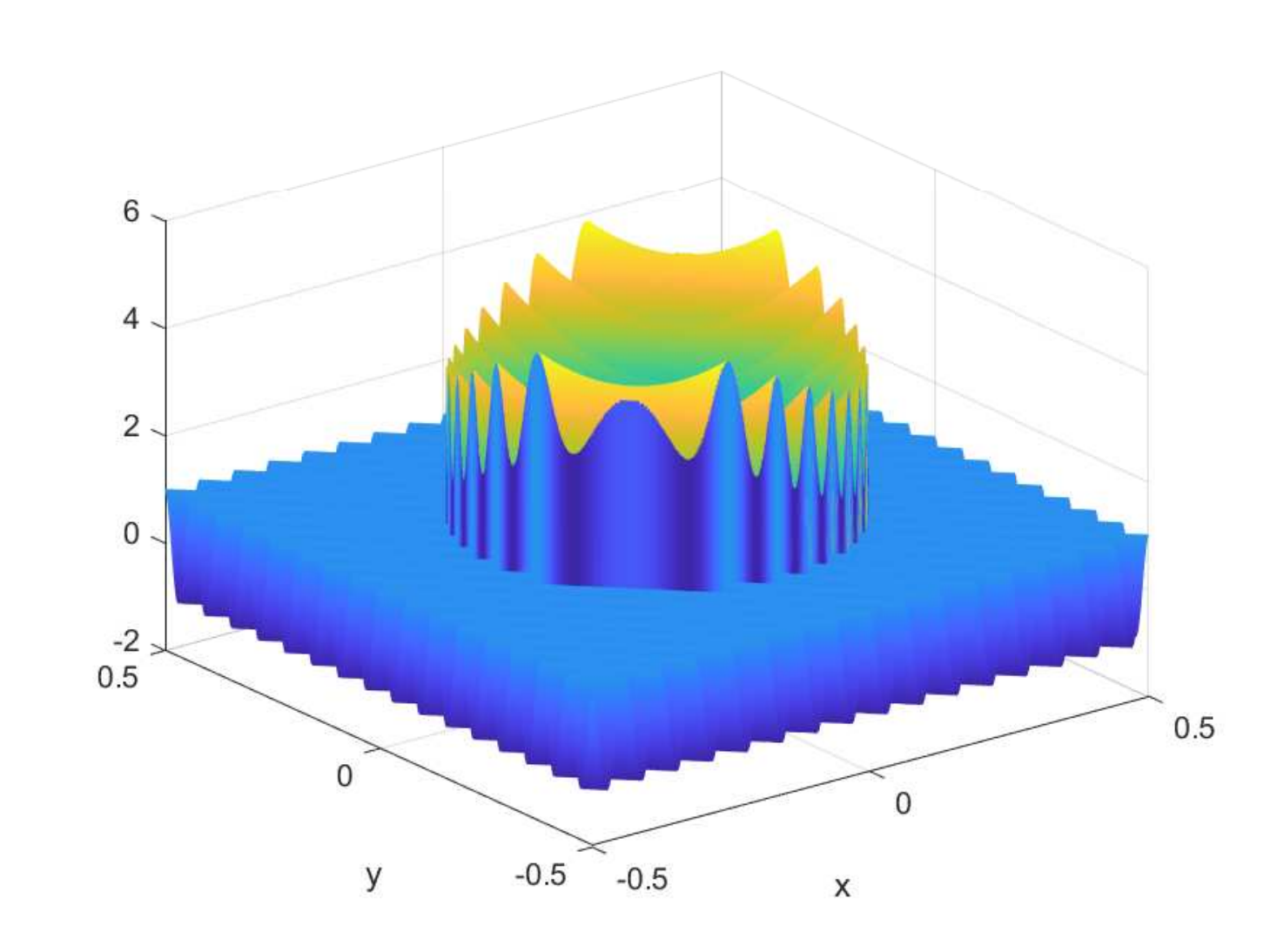}
	\end{subfigure}
	\begin{subfigure}[b]{0.3\textwidth}
		 \includegraphics[width=6cm,height=5cm]{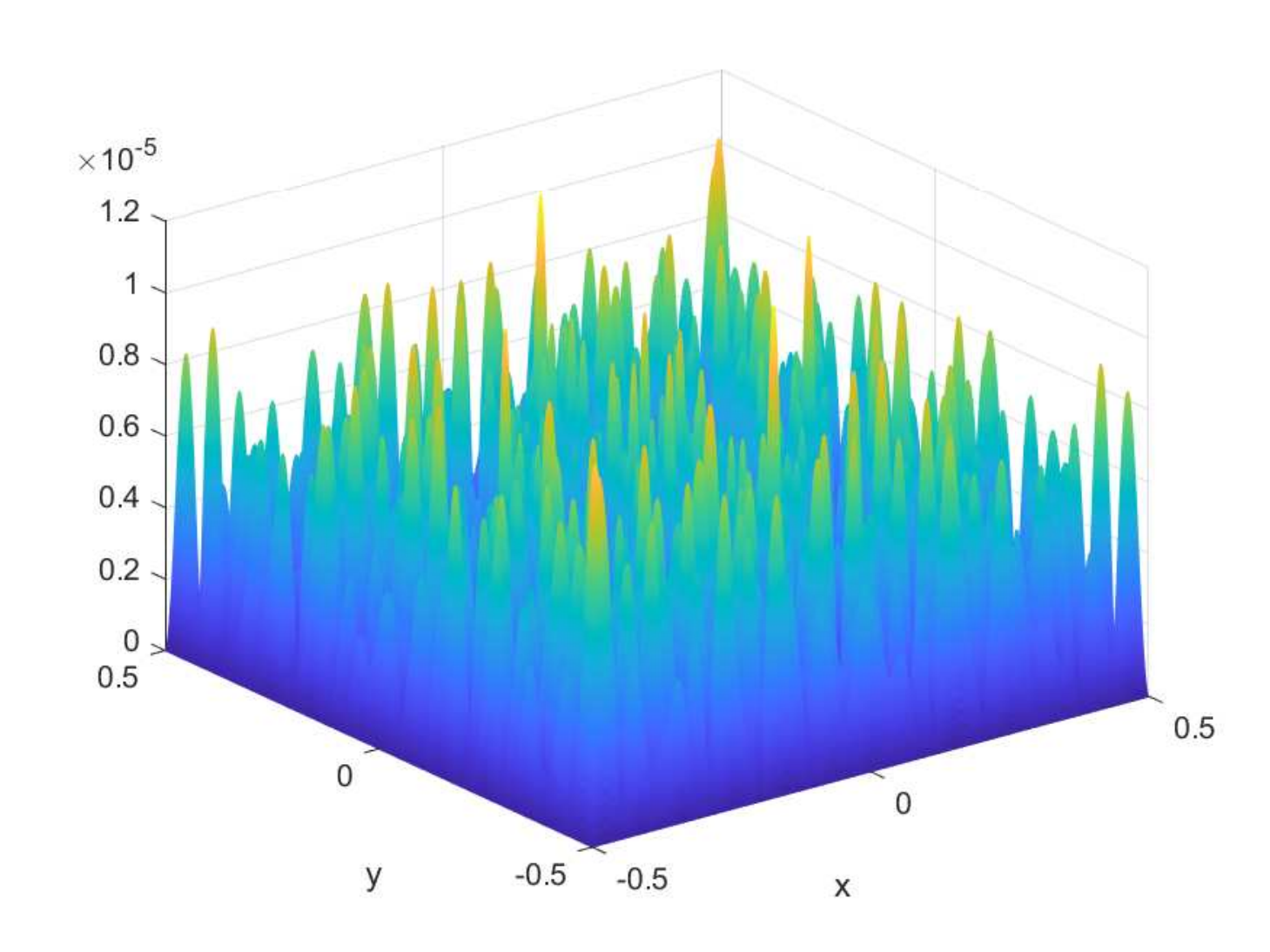}
	\end{subfigure}
	\caption
	{ \cref{interface:ex3} with $\ka_{+}=100$, $\ka_{-}=150$, $K=100$:  The interface curve $\Gamma$ (left), $u_h$ with $h=1/2^{11}$ (middle), $|u_h-u|$ with $h=1/2^{11}$ (right). }
	\label{fig:diskapa:1}
\end{figure}

\begin{figure}[htbp]
	\centering
	\begin{subfigure}[b]{0.3\textwidth}
		 \includegraphics[width=5.8cm,height=5cm]{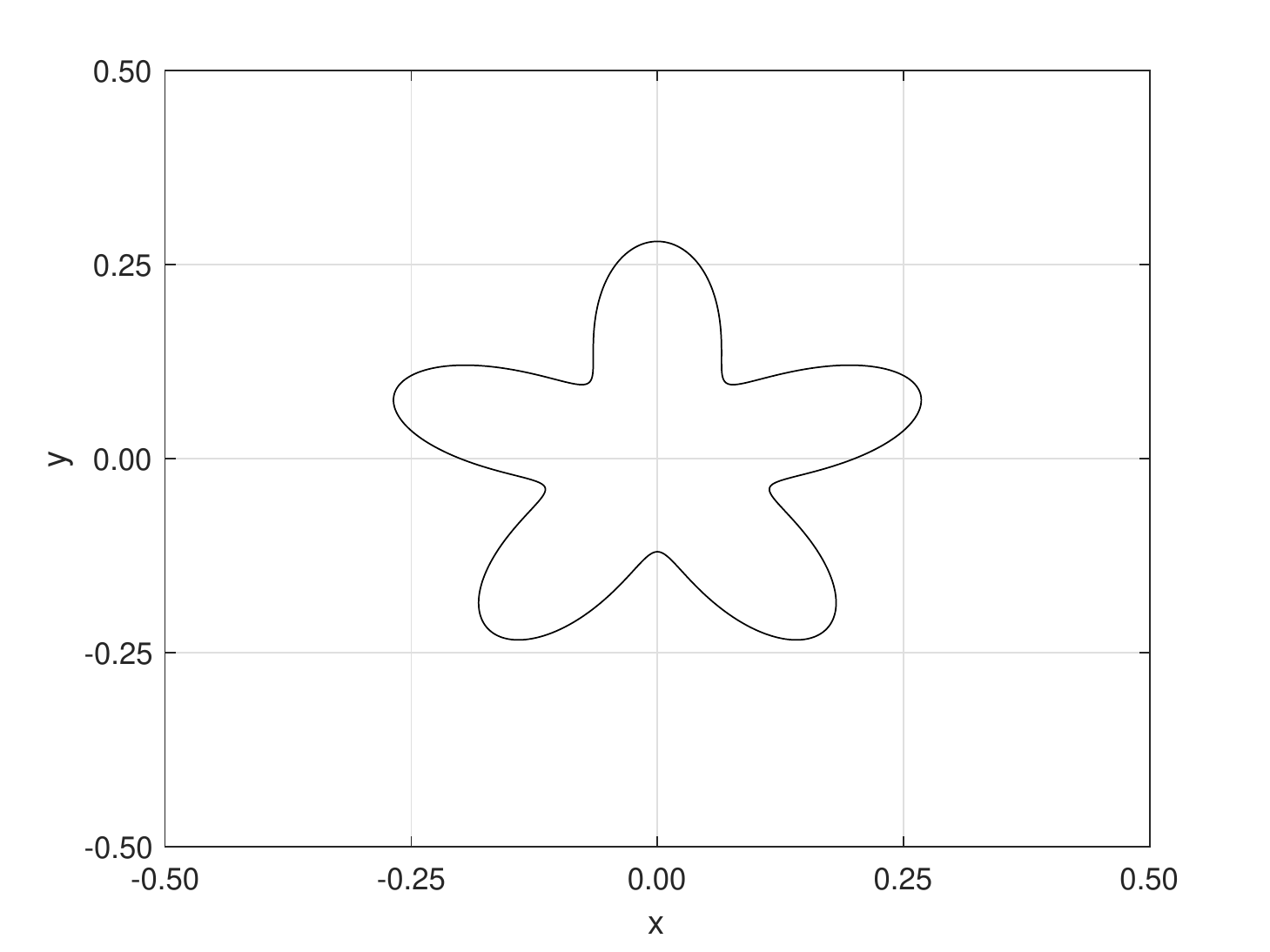}
	\end{subfigure}
	\begin{subfigure}[b]{0.3\textwidth}
	 \includegraphics[width=6cm,height=5cm]{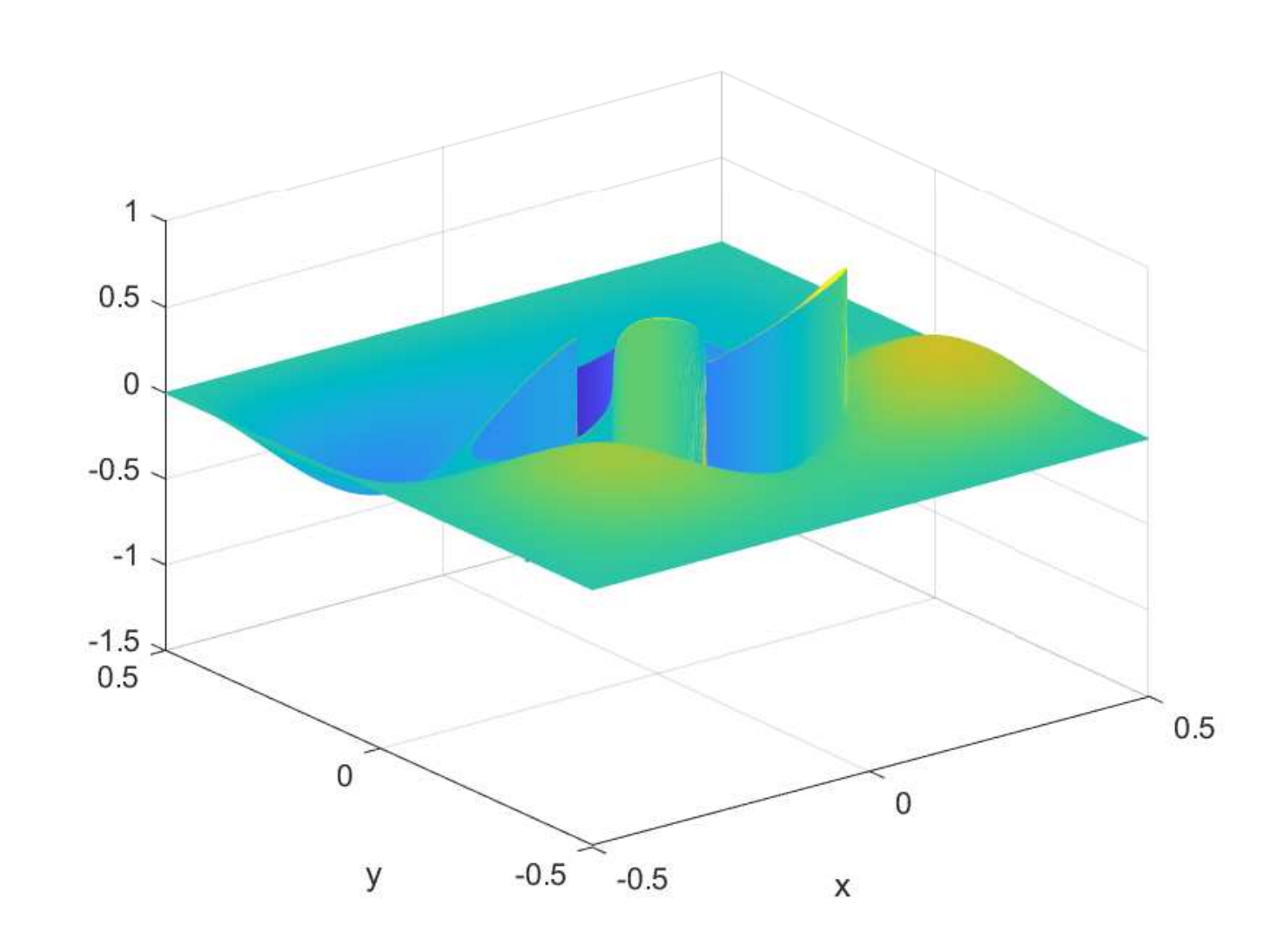}
\end{subfigure}
	\begin{subfigure}[b]{0.3\textwidth}
	 \includegraphics[width=6cm,height=5cm]{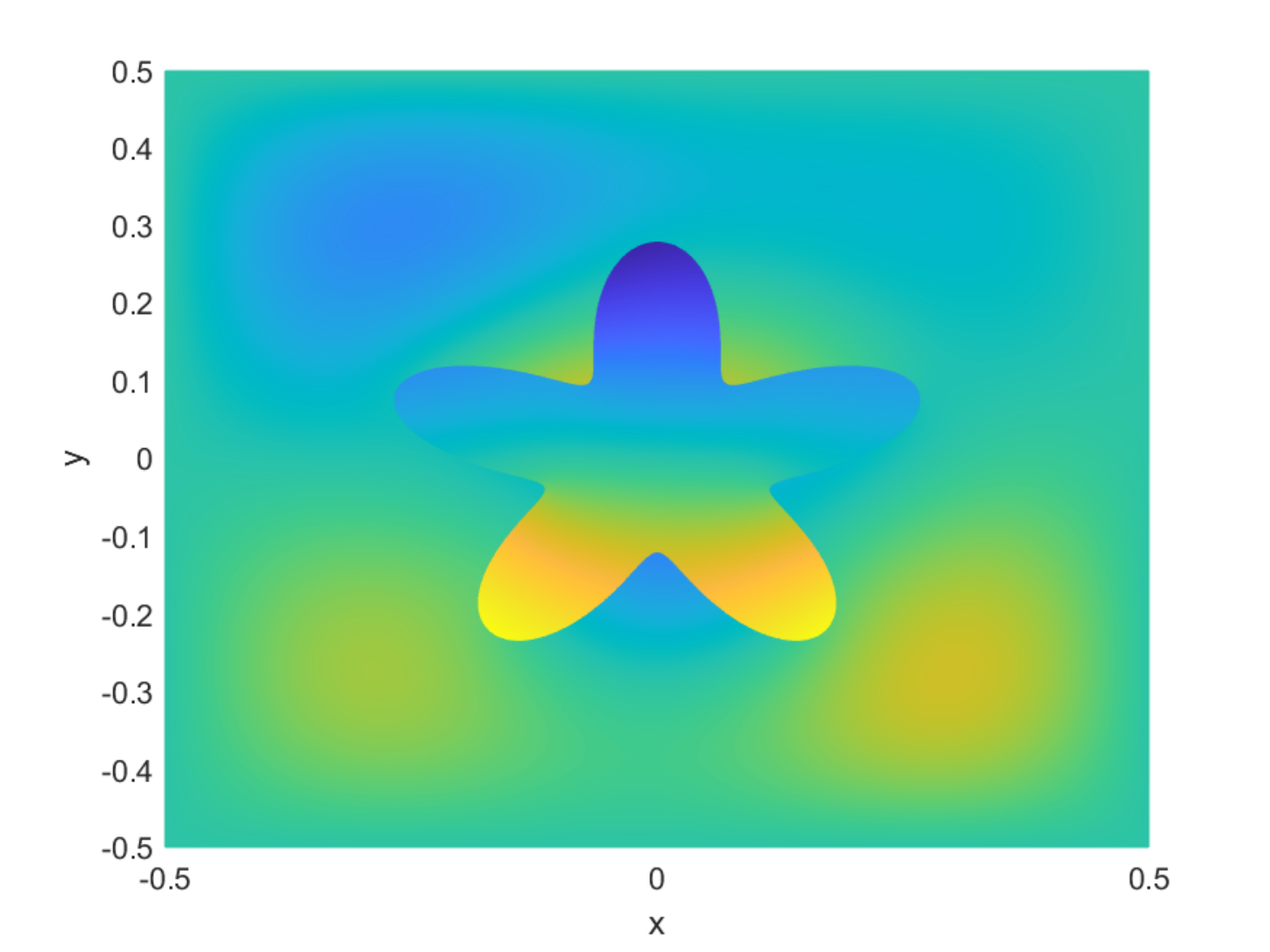}
\end{subfigure}
\begin{subfigure}[b]{0.3\textwidth}
	 \includegraphics[width=6cm,height=5cm]{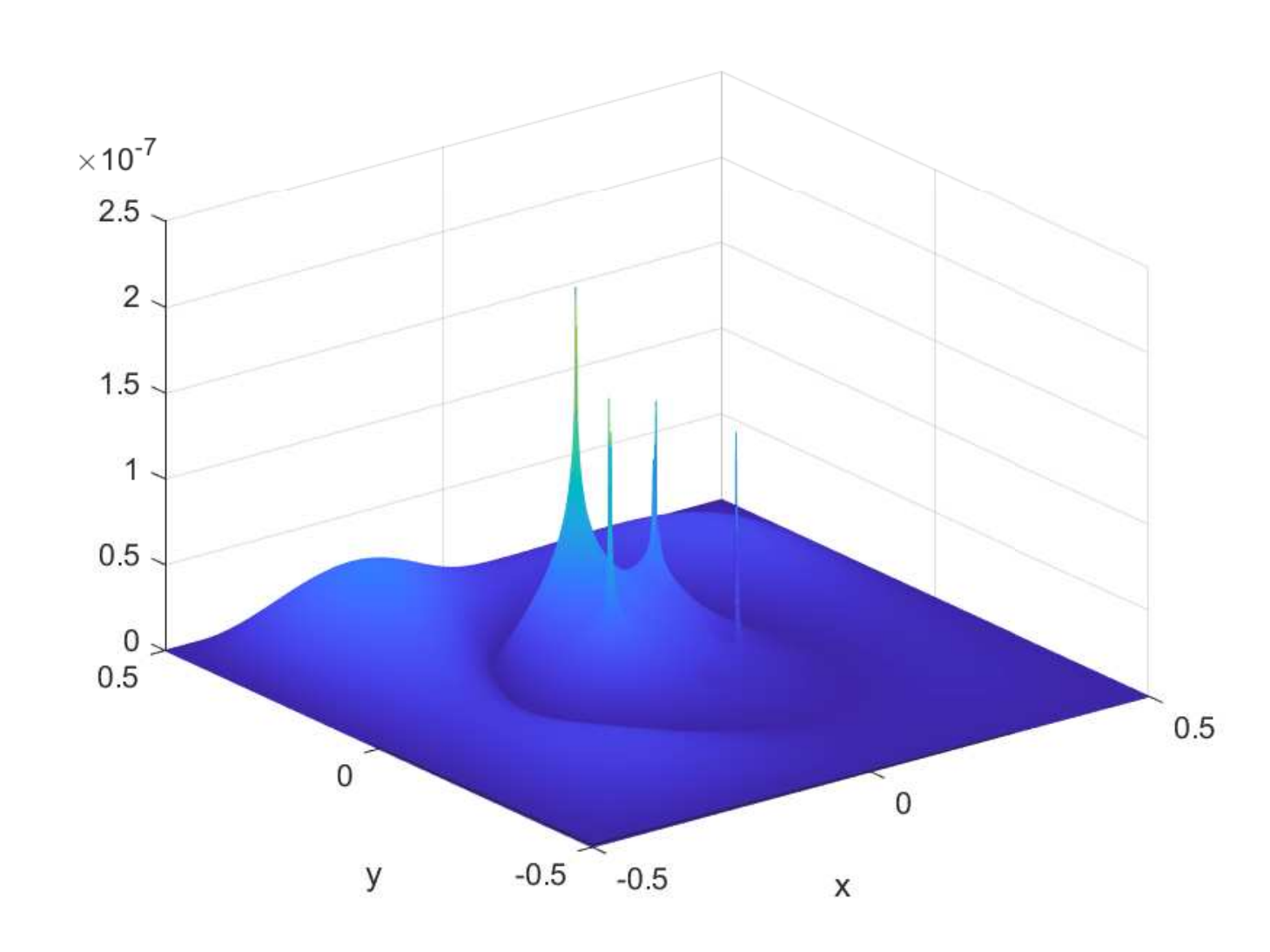}
\end{subfigure}
	\begin{subfigure}[b]{0.3\textwidth}
	 \includegraphics[width=6cm,height=5cm]{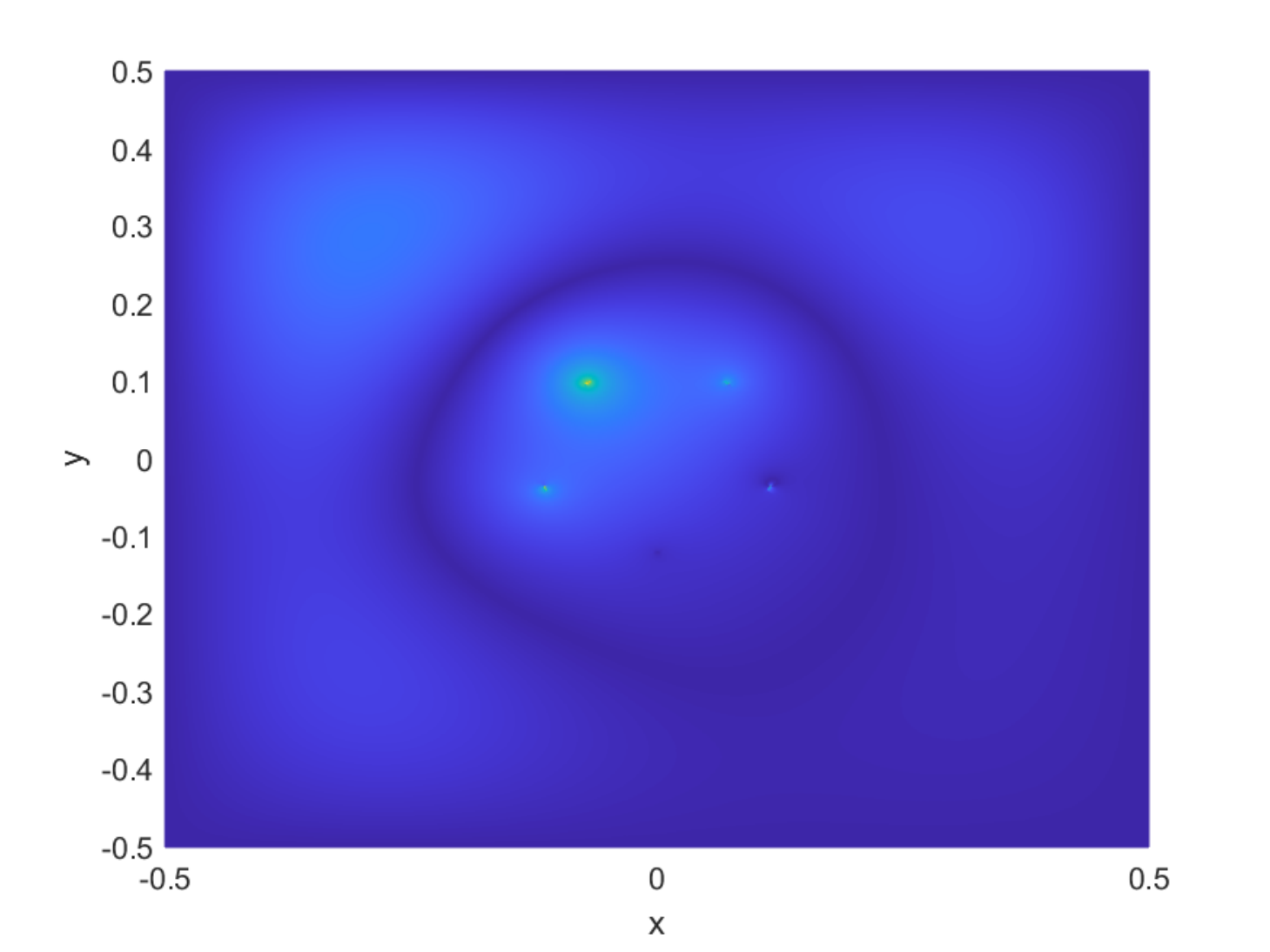}
\end{subfigure}
	\begin{subfigure}[b]{0.3\textwidth}
		 \includegraphics[width=6cm,height=5cm]{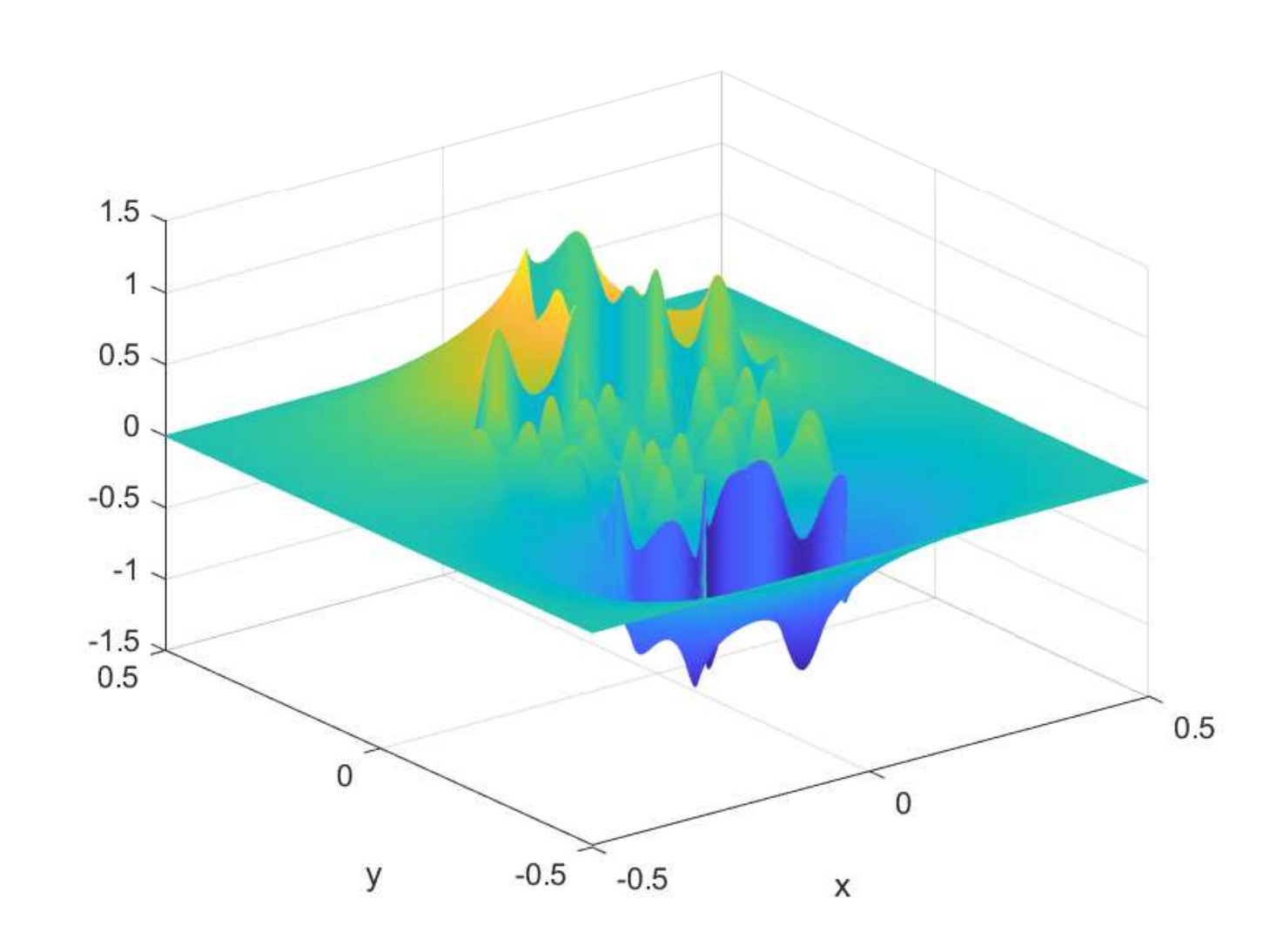}
	\end{subfigure}
	\caption
	{ \cref{interface:ex4}: The first row: the interface curve $\Gamma$ (left), $u_h$ with  $\ka_{+}=10$, $\ka_{-}=1$ and $h=1/2^{11}$ (middle and right). The second row: $|u_h-u_{h/2}|$ with  $\ka_{+}=10$, $\ka_{-}=1$ and $h=1/2^{10}$ (left and middle), $u_h$ with  $\ka_{+}=1$,  $\ka_{-}=100$ and $h=1/2^{11}$ (right).}
	\label{fig:diskapa:2}
\end{figure}

\begin{figure}[htbp]
	\centering
	\begin{subfigure}[b]{0.3\textwidth}
		 \includegraphics[width=6cm,height=5cm]{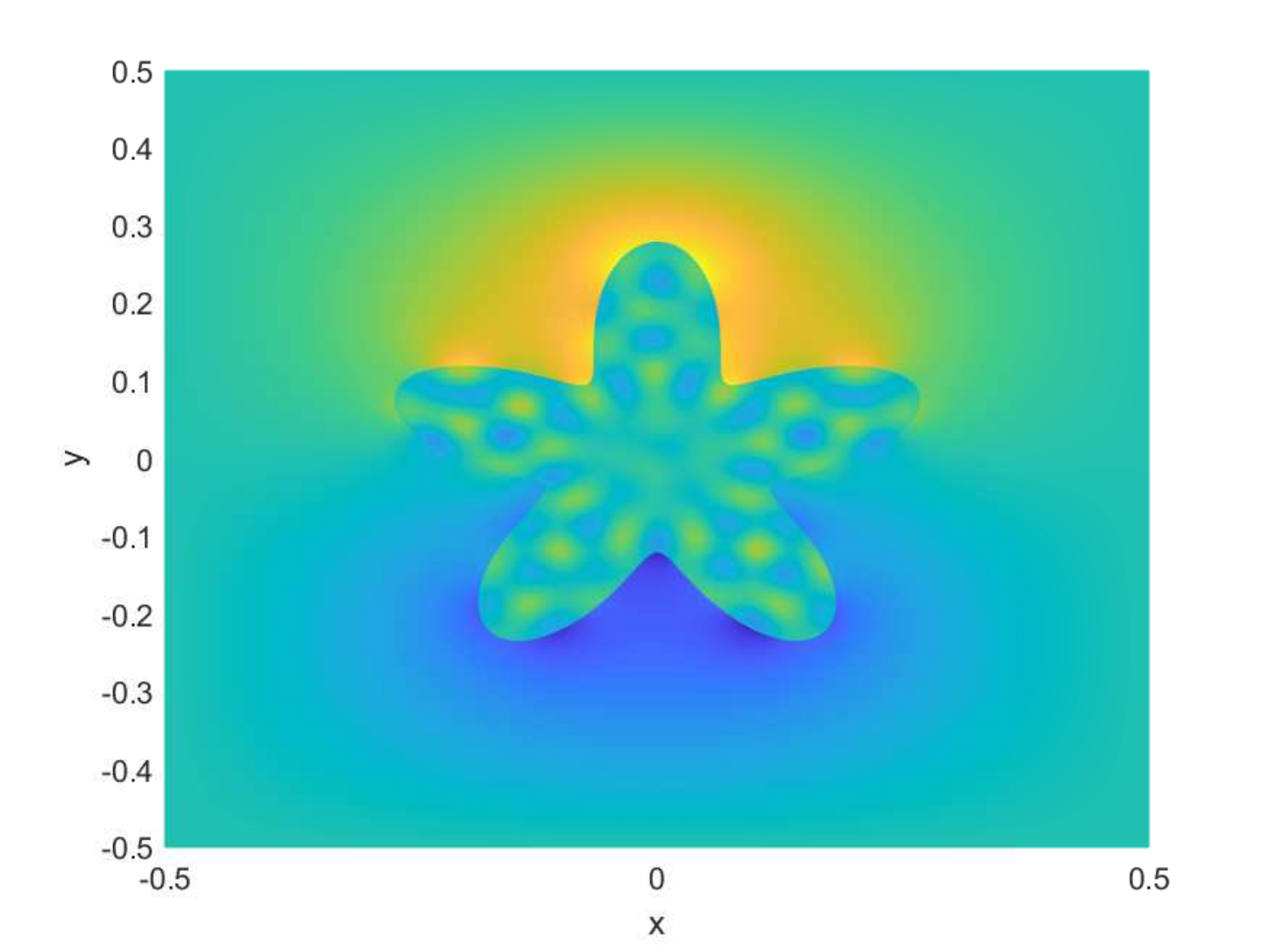}
	\end{subfigure}
	\begin{subfigure}[b]{0.3\textwidth}
		 \includegraphics[width=6cm,height=5cm]{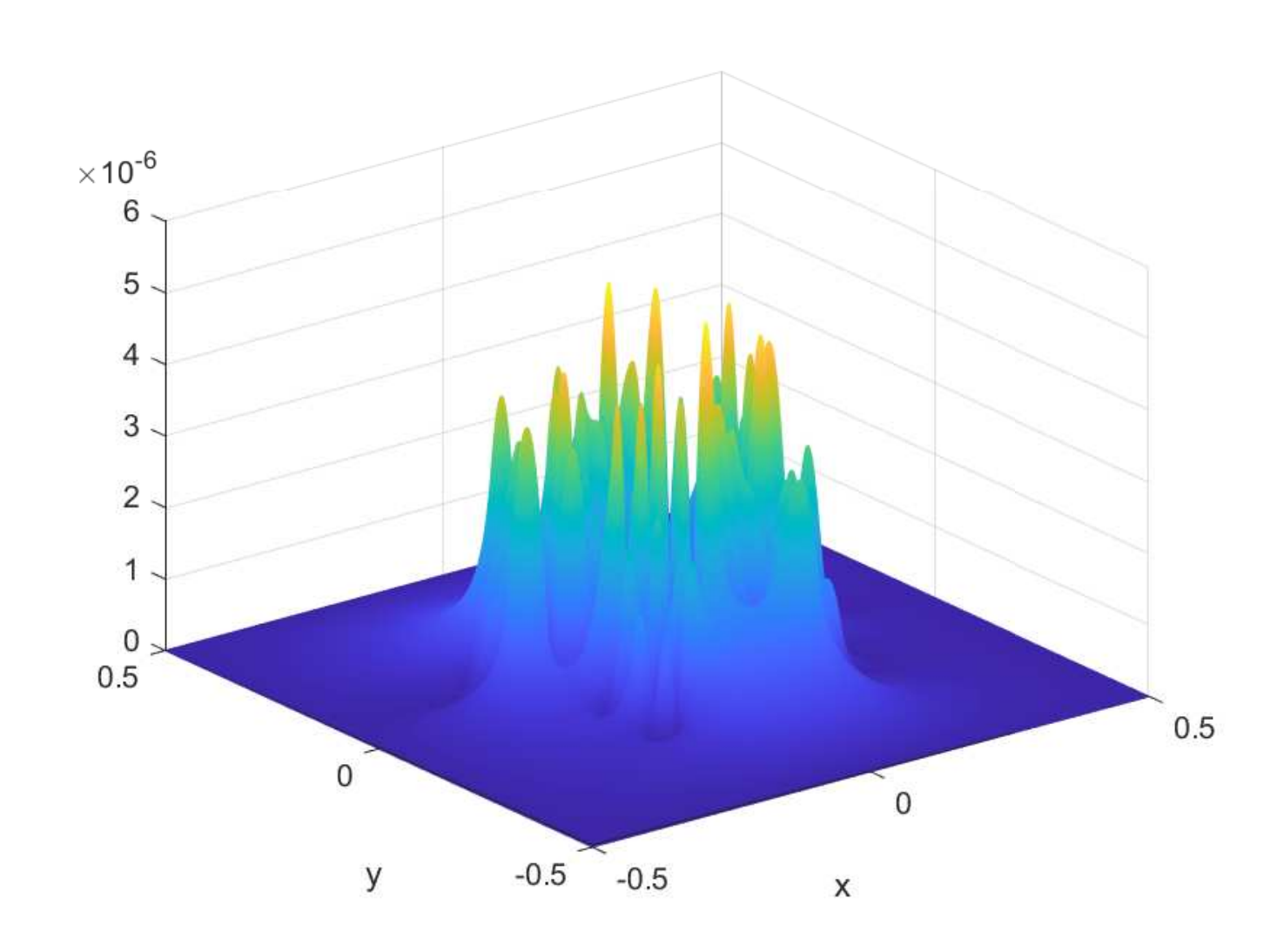}
	\end{subfigure}
	\begin{subfigure}[b]{0.3\textwidth}
		 \includegraphics[width=6cm,height=5cm]{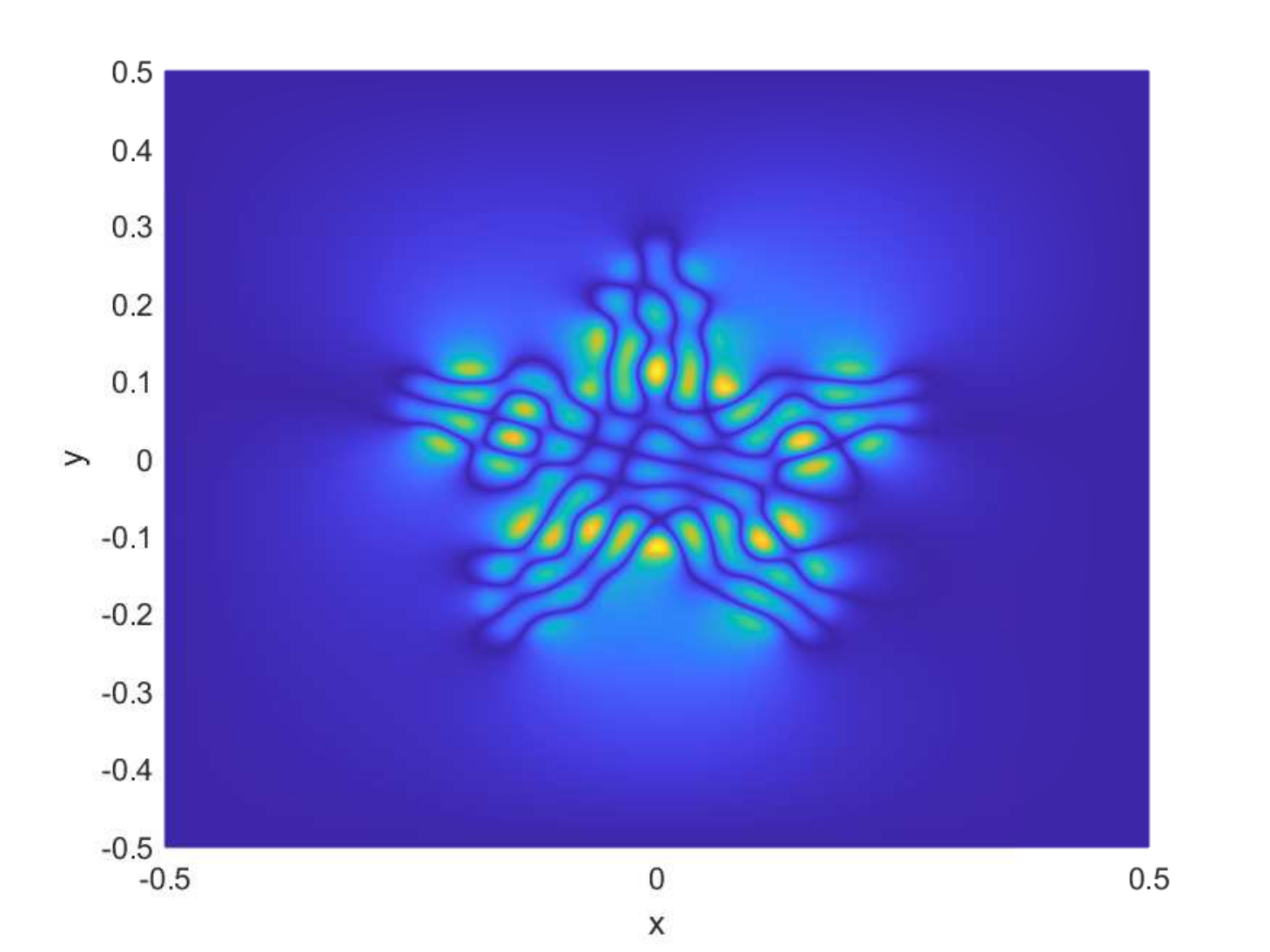}
	\end{subfigure}
	\caption
	{ \cref{interface:ex4}: $u_h$ with $\ka_{+}=1$, $\ka_{-}=100$ and  $h=1/2^{11}$ (left), $|u_h-u_{h/2}|$ with  $\ka_{+}=1$, $\ka_{-}=100$ and $h=1/2^{10}$ (middle and right). }
	\label{fig:diskapa:3}
\end{figure}

	We provide four numerical experiments here for this case (\cref{interface:ex1,interface:ex2,interface:ex3,interface:ex4}). The interfaces we consider are a five star interface, an eight star interface, an ellipse, and a circle. The first two examples consider continuous wavenumbers, while the last two examples consider discontinuous wavenumbers. We choose $(x_i^*,y_j^*)\in \Gamma$ to be the orthogonal projection of $(x_i,y_j)$ in this subsection.
	\begin{example}\label{interface:ex1}
		\normalfont
		Consider the problem \eqref{Qeques2} in $\Omega=(-1/2,1/2)^2$ with
		\begin{align*}
		& \Gamma=\{ (x,y) : \ x(\theta):=(0.2+0.05\sin(8\theta))\cos(\theta),\ y(\theta)=(0.2+0.05\sin(8\theta))\sin(\theta)\}, \\
		& \Omega_{+}=\{(x,y)\in \Omega : x^2(\theta)+y^2(\theta)> (0.2+0.05\sin(8\theta))^2\},\\
		& \Omega_{-}=\{(x,y)\in \Omega : x^2(\theta)+y^2(\theta)< (0.2+0.05\sin(8\theta))^2\},\\
		& \ka_{+}=\ka_{-}=400,\qquad g=-3, \qquad g_\Gamma=0,\\
		&  u_{+}=\sin(280 x)\cos(280 y), \quad u_{-}=\sin(280 x)\cos(280 y)+3,\\
		& -u_x(-1/2,y)-\ia\ka_{+} u(-1/2,y)= g_1 \quad \mbox{and}
			\quad u(1/2,y)= g_2 \quad \mbox{for} \quad y\in(-1/2,1/2),\\
		& -u_y(x,-1/2)= g_3 \quad \mbox{and}
			\quad u_y(x,1/2)-\ia\ka_{+} u(x,1/2)= g_4 \quad \mbox{for} \quad x\in(-1/2,1/2),
		\end{align*}
	where	the boundary data $g_1, \ldots,g_4$ and the source term  $f_{\pm}$ are obtained from the above data and the model problem.	
 See \cref{table:QSp6} and \cref{fig:interface:1} for numerical results.
	\end{example}
		\begin{table}[htbp]
		\caption{Numerical results of \cref{interface:ex1} with $h=1/2^J$ using our method.}
		\centering
		\setlength{\tabcolsep}{3mm}{
			\begin{tabular}{c|c|c|c|c|c}
				\hline
				 \multicolumn{6}{c}{\cref{interface:ex1} with $\ka_{+}=\ka_{-}=400$} \\
				\cline{1-6}
				$J$  	&  $\frac{2\pi}{ h\ka }$  &  $\frac{\|u_{h}-u\|_{2}}{\|u\|_2}$		 
				&order &  $\|u_{h}-u\|_{\infty}$
				&order   \\
				\hline
	8     &4.02     &1.99770E-01     &     &9.95173E-01     &\\
	9     &8.04     &1.48476E-03     &7.072     &6.98903E-03     &7.154\\
	10     &16.08     &1.09459E-05     &7.084     &5.38930E-05     &7.019\\
	11     &32.17     &7.51367E-08     &7.187     &3.76922E-07     &7.160\\	
				\hline
		\end{tabular}}
		\label{table:QSp6}
	\end{table}
	%
	%
%
%
	\begin{example}\label{interface:ex2}
		\normalfont
		Consider the problem \eqref{Qeques2} with $\Omega=(-3/2,3/2)^2$ with
			\begin{align*}
			& \Gamma =\{(x,y)\in \Omega \; :\; x^2+4y^2-1=0\},\\
			& \Omega_{+} =\{(x,y)\in \Omega : x^2+4y^2> 1\},\qquad  \Omega_{-} =\{(x,y)\in \Omega : x^2+4y^2< 1\},\\	
			& (\ka_{+},\ka_{-}) \in \{(0,0), (100,100)\}, \qquad x(\theta)=\cos(\theta),\quad y(\theta)={1}/{2}\sin(\theta),\\	
			& f_{+}=(4\pi)^2\sin(4\pi x)\sin(4\pi y),
			\qquad f_{-}=(4\pi)^2\cos(4\pi(x+y)), \\
			 &g=-\frac{|x'(\theta)y''(\theta)-x''(\theta)y'(\theta)|}{((x'(\theta))^2+(y'(\theta))^2)^{3/2}}, \qquad g_\Gamma=-\frac{|x'(\theta)y''(\theta)-x''(\theta)y'(\theta)|}{((x'(\theta))^2+(y'(\theta))^2)^{3/2}}, \quad \mbox{for} \quad \theta\in [0,2\pi),\\
			& u(-{3}/{2},y)=0, \quad \mbox{and} \quad  u({3}/{2},y)=0 \quad \mbox{for} \quad y\in(-{3}/{2},{3}/{2}),\\
			& u(x,-{3}/{2})=0, \quad \mbox{and} \quad   u(x,{3}/{2})=0  \quad \mbox{for} \quad x\in(-{3}/{2},{3}/{2}),
		\end{align*}
		i.e., the two jump functions $g$ and $g_\Gamma$ are curvatures of the interface curve $\Gamma$ and  we consider
		4 zero-Dirichlet boundary conditions.
		Note that the exact solution $u$ is unknown in this example.	See \cref{table:QSp7} and  \cref{fig:curvature} for numerical results.
	\end{example}
	\begin{table}[htbp]
		Numerical results of \cref{interface:ex2} with $h=3/2^J$ using our method.
		\centering
		\setlength{\tabcolsep}{1mm}{
			 \begin{tabular}{c|c|c|c|c|c|c|c|c|c|c}
				\hline
				 \multicolumn{5}{c|}{\cref{interface:ex2} with $\ka_{+}=\ka_{-}=0$} &
				 \multicolumn{6}{c}{\cref{interface:ex2} with $\ka_{+}=\ka_{-}=100$}\\
				\cline{1-11}
				$J$    &  $\|u_{h}-u_{h/2}\|_{2}$		
				&order &  $\|u_{h}-u_{h/2}\|_{\infty}$
				&order & 	$J$  & $\frac{2\pi}{h\ka}$  &   ${\|u_{h}-u_{h/2}\|_2}$    &order &  $\|u_{h}-u_{h/2}\|_{\infty}$
				&order  \\
				\hline
2   &4.1967E+05   &   &3.2562E+05   &   &   &   &   &   &   & \\
3   &3.5919E+03   &6.87   &3.5406E+03   &6.52   &   &   &   &   &   & \\
4   &3.8052E+01   &6.56   &4.0838E+01   &6.44   &   &   &   &   &   & \\
5   &2.9412E-01   &7.02   &3.8445E-01   &6.73   &6   &1.34   &1.0979E+03   &   &9.8002E+02   & \\
6   &1.9725E-03   &7.22   &1.9593E-03   &7.62   &7   &2.68   &1.3867E+01   &6.31   &1.3455E+01   &6.19 \\
7   &1.3459E-05   &7.20   &1.2578E-05   &7.28   &8   &5.36   &3.4798E-01   &5.32   &3.0775E-01   &5.45 \\
8   &8.9389E-08   &7.23   &8.0276E-08   &7.29   &9   &10.72   &4.7286E-03   &6.20   &4.2218E-03   &6.19 \\
9   &7.2057E-10   &6.95   &8.4663E-10   &6.57   &10   &21.45   &7.1356E-05   &6.05   &6.3680E-05   &6.05 \\
				\hline
		\end{tabular}}
		\label{table:QSp7}
	\end{table}
	%
	%
	%
%
%
%
	\begin{example}\label{interface:ex3}
	\normalfont
	Consider the problem \eqref{Qeques2} in $\Omega=(-1/2,1/2)^2$ with
	\begin{align*}
		& \Gamma=\{ (x,y) : \ x(\theta)={3}/{10}\cos(\theta),\ y(\theta)={3}/{10}\sin(\theta)\}, \\
		& \Omega_{+}=\{(x,y)\in \Omega : (10x/3)^2+(10y/3)^2> 1\},\\
		&\Omega_{-}=\{(x,y)\in \Omega : (10x/3)^2+(10y/3)^2< 1\},\\	
		& (\ka_{+},\ka_{-},K) \in \{(90,100,70), (100,150,100)\}, \\
		& u_{+}=\cos(K(x+y)),
		\qquad u_{-}=\cos(K(x+y))+40(x^2+y^2)+20xy, \\
		 &g=-{18}/{5}-{9}/{5}\cos(\theta)\sin(\theta), \qquad g_\Gamma=-12\cos(\theta)\sin(\theta)-24, \quad \mbox{for} \quad \theta\in [0,2\pi),\\
		& u(-1/2,y)=g_1, \quad \mbox{and} \quad  u(1/2,y)=g_2 \quad \mbox{for} \quad y\in(-1/2,1/2),\\
		& u(x,-1/2)=g_3, \quad \mbox{and} \quad  u(x,1/2)=g_4  \quad \mbox{for} \quad x\in(-1/2,1/2),
	\end{align*}
where the boundary data $g_1,\ldots,g_4$ and the source term $f_{\pm}$ are obtained from the above data and the model problem.	
 See \cref{table:QSp8} and \cref{fig:diskapa:1} for numerical results.	 
\end{example}	
		\begin{table}[htbp]
	\caption{Numerical results of \cref{interface:ex3} with $h=1/2^J$ using our method.}
	\centering
	\setlength{\tabcolsep}{1mm}{
		 \begin{tabular}{c|c|c|c|c|c|c|c|c|c|c|c}
			\hline
			 \multicolumn{6}{c|}{\cref{interface:ex3} with $\ka_{+}=90$, $\ka_{-}=100$, $K=70$} &
			 \multicolumn{6}{c}{\cref{interface:ex3} with $\ka_{+}=100$, $\ka_{-}=150$, $K=100$}\\
			\cline{1-12}
			$J$  	&  $\frac{2\pi}{ h\ka_{-} }$  &  $\frac{\|u_{h}-u\|_{2}}{\|u\|_2}$		 
			&order &  $\|u_{h}-u\|_{\infty}$
			&order  & 	$J$  	&  $\frac{2\pi}{ h\ka_{-} }$  &  $\frac{\|u_{h}-u\|_{2}}{\|u\|_2}$		
			&order &  $\|u_{h}-u\|_{\infty}$
			&order  \\
			\hline
7   &8.0   &1.8683E+00   &   &9.3194E+00   &   &7   &5.4   &1.2698E+00   &   &7.2414E+00   & \\
8   &16.1   &1.1556E-02   &7.3   &5.6877E-02   &7.4   &8   &10.7   &5.7245E-02   &4.5   &2.5975E-01   &4.8 \\
9   &32.2   &3.5860E-04   &5.0   &1.9017E-03   &4.9   &9   &21.4   &2.3353E-03   &4.6   &1.2106E-02   &4.4 \\
10   &64.3   &1.0785E-05   &5.1   &5.8872E-05   &5.0   &10   &42.9   &8.4024E-05   &4.8   &4.1842E-04   &4.9 \\
11   &128.7   &3.4121E-07   &5.0   &1.8572E-06   &5.0   &11   &85.8   &2.3915E-06   &5.1   &1.1928E-05   &5.1 \\
			\hline
	\end{tabular}}
	\label{table:QSp8}
\end{table}
%
%
%
%
%
%
%
%
	\begin{example}\label{interface:ex4}
		\normalfont
		Consider the problem \eqref{Qeques2} in $\Omega=(-1/2,1/2)^2$ with
		\begin{align*}
			& \Gamma=\{ (x,y) : \ x(\theta)=(0.2+0.08\sin(5\theta))\cos(\theta),\ y(\theta)=(0.2+0.08\sin(5\theta))\sin(\theta)\}, \\
			& \Omega_{+}=\{(x,y)\in \Omega : x^2(\theta)+y^2(\theta)> (0.2+0.08\sin(5\theta))^2\},\\
			& \Omega_{-}=\{(x,y)\in \Omega : x^2(\theta)+y^2(\theta)< (0.2+0.08\sin(5\theta))^2\},\\
			& (\ka_{+},\ka_{-})\in\{(10,1), (1,100)\}, \qquad  f_{+}=\sin(2\pi x)\sin(2\pi y),
			\qquad f_{-}=\cos(2\pi x)\cos(2\pi y), \\
			&g=\sin(\theta), \qquad g_\Gamma=\cos(\theta), \quad \mbox{for} \quad \theta\in [0,2\pi),\\
			& u(-1/2,y)=0, \quad \mbox{and} \quad  u(1/2,y)=0 \quad \mbox{for} \quad y\in(-1/2,1/2),\\
			& u(x,-1/2)=0, \quad \mbox{and} \quad  u(x,1/2)=0  \quad \mbox{for} \quad x\in(-1/2,1/2).
		\end{align*} 	
		Note that the exact solution $u$ is unknown in this example. See \cref{table:QSp9} and \cref{fig:diskapa:2,fig:diskapa:3} for numerical results.	
	\end{example}	
		\begin{table}[htbp]
	\caption{Numerical results of \cref{interface:ex4} with $h=1/2^J$ using our method.}
	\centering
	\setlength{\tabcolsep}{1mm}{
		 \begin{tabular}{c|c|c|c|c|c|c|c|c|c|c|c}
			\hline
			 \multicolumn{6}{c|}{\cref{interface:ex4} with  $\ka_{+}=10$ and $\ka_{-}=1$ } &
			 \multicolumn{6}{c}{\cref{interface:ex4} with  $\ka_{+}=1$ and $\ka_{-}=100$ }  \\
			\cline{1-12}
			$J$  	&  $\frac{2\pi}{ h\ka_{+} }$  &  $\|u_{h}-u_{h/2}\|_{2}$		
			&order &  $\|u_{h}-u_{h/2}\|_{\infty}$
			&order  & 	$J$	&  $\frac{2\pi}{ h\ka_{-} }$  &  $\|u_{h}-u_{h/2}\|_{2}$		
			&order &  $\|u_{h}-u_{h/2}\|_{\infty}$
			&order \\
			\hline
6   &40.2   &8.2461E-01   &   &3.6325E+00   &   &   &   &   &   &   & \\
7   &80.4   &8.0665E-03   &6.7   &4.4579E-02   &6.3   &7   &8.0   &3.5212E-02   &   &2.6626E-01   & \\
8   &160.8   &9.3400E-05   &6.4   &6.5308E-04   &6.1   &8   &16.1   &9.6191E-04   &5.2   &6.9248E-03   &5.3 \\
9   &321.7   &3.0871E-06   &4.9   &2.2701E-05   &4.8   &9   &32.2   &2.4508E-05   &5.3   &1.5983E-04   &5.4 \\
10   &643.4   &1.9575E-08   &7.3   &2.2072E-07   &6.7   &10   &64.3   &8.0117E-07   &4.9   &5.2179E-06   &4.9 \\
			\hline
	\end{tabular}}
	\label{table:QSp9}
\end{table}
%
%
%
%
%
%
%
%
%
%
%
	\section{Proofs of \cref{thm:regular:all:Ckl,thm:regular:interior,thm:regular:Robin:1,thm:corner:1,thm:corner:2,thm:interface,fluxtm2,fluxtm3}}
	\label{sec:proofs}

For the proofs of the theorems in this paper we first need to establish
some auxiliary identities about the solution $u$ of the Helmholtz interface problem in \eqref{Qeques2}.

For $x\in \R$, the floor function $\lfloor x\rfloor$ is defined to be the largest integer less than or equal to $x$.
For an integer $m$, we define
\[
\odd(m):=
\frac{1-(-1)^m}{2}=
\begin{cases}
	0, &\text{if $m$ is even},\\
	1, &\text{if $m$ is odd}.
\end{cases}
\]
Recall that
\[
\begin{split}
&\ind_{M+1}:=\{(m,n)\in \NN^2 \; : \; m+n\le M+1\}, \qquad M+1\in \NN,\\
&\ind_{M+1}^{ 1}:=\{(m,n)\in \ind_{M+1}\; :   m=0,1\},\qquad 
\ind_{M+1}^{ 2}:=\ind_{M+1}\setminus \ind_{M+1}^{ 1},
\end{split}
\]
\[
u^{(m,n)}:=\frac{\partial^{m+n} u}{ \partial^m x \partial^n y}(x_i^*,y_j^*)
\quad\mbox{and}\quad
f^{(m,n)}:=\frac{\partial^{m+n} f}{ \partial^m x \partial^n y}(x_i^*,y_j^*).
\]
Since the function $u$ is a solution to the partial differential equation in \eqref{Qeques2}, all quantities $u^{(m,n)}\in \{ u^{(m,n)}, (m,n)\in \ind_{M+1} \}$ are not independent of each other. The following \eqref{umnV2:relation} and \eqref{umnH2:relation} describe this dependence (\eqref{umnV2:relation} and \eqref{umnH2:relation} can be obtained by the proof of \cite[Lemma~2.1]{FHM21a} and \cite[Lemma~2.1]{FHM21b}),
	\be\label{umnV2:relation}
	u^{(m,n)}=(-1)^{\lfloor \frac{m}{2}\rfloor} \sum_{i=0}^{\lfloor \frac{m}{2}\rfloor} {\lfloor \frac{m}{2}\rfloor \choose i} {\ka}^{2i}
	u^{(\odd(m),2\lfloor \frac{m}{2}\rfloor+n-2i)}+\sum_{i=1}^{\lfloor \frac{m}{2}\rfloor}
	\sum_{j=0}^{i-1}
	(-1)^{i-1} {i-1 \choose j} {\ka}^{2(i-j-1)}f^{(m-2i,n+2j)}
	\ee
	for all $(m,n)\in \ind_{M+1}^{ 2}$,
	\be \label{umnH2:relation}
	u^{(m,n)}=(-1)^{\lfloor \frac{n}{2}\rfloor} \sum_{i=0}^{\lfloor \frac{n}{2}\rfloor} {\lfloor \frac{n}{2}\rfloor \choose i} {\ka}^{2i}
	u^{(2\lfloor \frac{n}{2}\rfloor+m-2i,\odd(n))}
	+\sum_{i=1}^{\lfloor \frac{n}{2}\rfloor}
	\sum_{j=0}^{i-1}
	(-1)^{i-1} {i-1 \choose j} {\ka}^{2(i-j-1)}f^{(m+2j,n-2i)}
	\ee	
	for all $(n,m) \in \ind_{M+1}^{ 2}$, where 	
 $u$ is a smooth function satisfying $\Delta u+{\ka}^2u=f$ in $\Omega\setminus \Gamma$ and the point $(x_i^*,y_j^*)\in \Omega\setminus \Gamma$.
See \cite[Figure 6]{FHM21a} for an illustration of how each $u^{(m,n)}$ with $(m,n) \in \Lambda_{7}$ is categorized based on $\ind_{7}^{j}$ with $j \in \{1,2\}$.

For a smooth function $u$, the values $u(x+x_i^*,y+y_j^*)$ are well approximated by its Taylor polynomial. For $x, y \in (-2h,2h)$,
\be \label{u:approx}
u(x+x_i^*,y+y_j^*)=
\sum_{(m,n)\in \ind_{M+1}} \frac{u^{(m,n)}}{m!n!}x^m y^{n}
+\bo(h^{M+2}), \qquad h \rightarrow 0.
\ee
From \eqref{umnV2:relation}, we have
\be \label{U2:Seperate}
\begin{aligned}
	\sum_{(m,n)\in \ind_{M+1}^{ 2}} \frac{x^m y^{n}}{m!n!} u^{(m,n)}
	&=\overbrace {\sum_{(m,n)\in \ind_{M+1}^{ 2}} \frac{x^m y^{n}}{m!n!} \Big\{ (-1)^{\lfloor \frac{m}{2}\rfloor} \sum_{i=0}^{\lfloor \frac{m}{2}\rfloor} {\lfloor \frac{m}{2}\rfloor \choose i} {\ka}^{2i}
		u^{(\odd(m),2\lfloor \frac{m}{2}\rfloor+n-2i)} \Big\} }^{=:I_1} \\
	& \quad + \underbrace{ \sum_{(m,n)\in \ind_{M+1}^{ 2}} \frac{x^m y^{n}}{m!n!}  \Big\{ \sum_{i=1}^{\lfloor \frac{m}{2}\rfloor}
		\sum_{j=0}^{i-1}
		(-1)^{i-1} {i-1 \choose j} {\ka}^{2i-2j-2}f^{(m-2i,n+2j)}\Big\}}_{=:I_2},
\end{aligned}
\ee
where the first summation $I_1$ above can be expressed as
\begin{align*}
	\nonumber
	I_1&
	=\sum_{(m,n)\in \ind_{M+1}^{ 2}  \atop \ell=\frac{m}{2},\ \text{even } m} \frac{(-1)^{\ell}x^{2\ell} y^{n}}{(2\ell)!n!}  \sum_{i=0}^{\ell} {\ell \choose i} {\ka}^{2i}
	u^{(0,2\ell+n-2i)}
	+ \sum_{(m,n)\in \ind_{M+1}^{ 2}  \atop \ell=\frac{m-1}{2},\ \text{odd } m  } \frac{(-1)^{\ell}x^{2\ell+1} y^{n}}{(2\ell+1)!n!}  \sum_{i=0}^{\ell} {\ell \choose i} {\ka}^{2i}
	u^{(1,2\ell+n-2i)}\\
	\nonumber
	& =\sum_{n=2}^{M+1} \sum_{\ell=1}^{\lfloor \frac{n}{2}\rfloor} \frac{(-1)^{\ell}x^{2\ell} y^{n-2\ell}}{(2\ell)!(n-2\ell)!}  \sum_{i=0}^{\ell} {\ell \choose i} {\ka}^{2i}
	u^{(0,n-2i)}+\sum_{n=2}^{M} \sum_{\ell=1}^{\lfloor \frac{n}{2}\rfloor} \frac{(-1)^{\ell}x^{2\ell+1} y^{n-2\ell}}{(2\ell+1)!(n-2\ell)!}  \sum_{i=0}^{\ell} {\ell \choose i} {\ka}^{2i}
	u^{(1,n-2i)}\\
	& =\sum_{(m,n)\in \ind_{M+1}^{ 1}  \atop n\ge2 } 	\sum_{\ell=1}^{\lfloor \frac{n}{2}\rfloor}
	\frac{(-1)^{\ell}x^{m+2\ell} y^{n-2\ell}}{(m+2\ell)!(n-2\ell)!} \sum_{i=0}^{\ell} {\ell \choose i} {\ka}^{2i}
	u^{(m,n-2i)},
\end{align*}
and the second summation $I_2$ above can be expressed as
\be \label{Qv}
\begin{aligned}
	I_2& =\sum_{(m,n)\in \ind_{M-1} } \sum_{\ell=1}^{1+\lfloor \frac{n}{2}\rfloor}
	\sum_{p=0}^{\ell-1} (-1)^{\ell-1} {\ell-1 \choose p}
	{\ka}^{2(\ell-p-1)} f^{(m,n+2(p+1-\ell))}	 \frac{x^{m+2\ell} y^{n-2\ell+2}}{(m+2\ell)!(n-2\ell+2)!}\\
	& =\sum_{(m,n)\in \ind_{M-1} } \sum_{j \in \{n+2p | p \in \NN, \atop  n+2p\le M-1-m\}} \sum_{\ell=1+\frac{j-n}{2}}^{1+\lfloor\frac{j}{2}\rfloor}  \	(-1)^{\ell-1} {\ell-1 \choose \frac{j-n}{2}}
	{\ka}^{j-n} 	\frac{x^{m+2\ell} y^{j-2\ell+2}}{(m+2\ell)!(j-2\ell+2)!} f^{(m,n)}\\
	& =\sum_{(m,n)\in \ind_{M-1} }H_{M+1,m,n}(x,y) f^{(m,n)}.
\end{aligned}
\ee
where $H_{M+1,m,n}(x,y)$ in \eqref{Hel:HM1mnxy}.

Hence, using the right-hand side of \eqref{u:approx} and the definitions of $\ind_{M+1}^{ 1},\ind_{M+1}^{ 2}$ in \eqref{indV12}, we have
\be \label{Gv}
\begin{aligned}
	I_1 + & \sum_{(m,n)\in \ind_{M+1}^{ 1}} \frac{x^m y^{n}}{m!n!} u^{(m,n)}
	=\sum_{(m,n)\in \ind_{M+1}^{ 1}} \sum_{i=0}^{\lfloor \frac{n}{2}\rfloor}	 \sum_{\ell=i}^{\lfloor \frac{n}{2}\rfloor}
	\frac{(-1)^{\ell}x^{m+2\ell} y^{n-2\ell}}{(m+2\ell)!(n-2\ell)!}  {\ell \choose i} {\ka}^{2i}
	u^{(m,n-2i)}\\
	&  =\sum_{(m,n)\in \ind_{M+1}^{ 1}}
	\sum_{i \in \{n+2p | p \in \NN, \atop n+2p\le M+1-m\}}	 \sum_{\ell=\frac{i-n}{2}}^{\lfloor \frac{i}{2}\rfloor}
	\frac{(-1)^{\ell}x^{m+2\ell} y^{i-2\ell}}{(m+2\ell)!(i-2\ell)!}  {\ell \choose \frac{i-n}{2}} {\ka}^{i-n} u^{(m,n)}\\
	&  =\sum_{(m,n)\in \ind_{M+1}^{ 1}} G_{M+1,m,n}(x,y) u^{(m,n)},
\end{aligned}
\ee
where $G_{M+1,m,n}(x,y)$ in \eqref{Hel:GM1mnxy}.

Suppose $x,y \in (-2h,2h)$. The lowest degree of $h$ for each polynomial $G_{M+1,m,n}(x,y)$ with $(m,n)\in \ind_{M+1}^{ 1}$ in \eqref{Hel:GM1mnxy} is $m+n$. The lowest degree of $h$ for each polynomial $H_{M+1,m,n}(x,y)$ with $(m,n)\in \ind_{M-1}$ in \eqref{Hel:HM1mnxy} is $m+n+2$. Therefore, by \eqref{U2:Seperate}-\eqref{Gv}, the approximation of $u(x+x_i^*,y+y_j^*)$ with $(x,y) \in (-2h,2h)$ in \eqref{u:approx} can be written as 
\be
\label{GQ:V}
u(x+x_i^*,y+y_j^*)  = \sum_{(m,n)\in \ind_{M+1}^{ 1}} u^{(m,n)} G_{M+1,m,n}(x,y) + \sum_{(m,n)\in \ind_{	\tilde{M}-1}} f^{(m,n)} H_{	\tilde{M}+1,m,n}(x,y) + \mathcal{O}(h^{M+2}),\; h\rightarrow 0,
\ee
where $M,	\tilde{M} \in \NN$	and $	 \tilde{M}\ge M$. By a similar calculation, for  $(x,y) \in (-2h,2h)$, we have
\be
\label{GQ:H}
u(x+x_i^*,y+y_j^*)  = \sum_{(n,m)\in \ind_{M+1}^{ 1}} u^{(m,n)} G_{M+1,n,m}(y,x) + \sum_{(m,n)\in \ind_{	\tilde{M}-1}} f^{(m,n)} H_{	\tilde{M}+1,n,m}(y,x) + \mathcal{O}(h^{M+2}), \; h \rightarrow 0.
\ee
Identities \eqref{GQ:V}-\eqref{GQ:H} are critical in finding compact stencils achieving a desired consistency order.

In the rest of this section, we prove the main results stated in
\cref{sec:sixord}. The idea of proofs is to
first construct all possible compact stencils with the maximum consistency order (in the context of methods relying on Taylor expansions and our sort of techniques) and then to minimize the average truncation error of plane waves over the free parameters of stencils to reduce pollution effect.

\begin{proof}[Proof of \cref{thm:regular:all:Ckl}]
	
	Let us consider the following discretization operator at a regular point $(x_i,y_j)$:
	\[
h^{-2}	\mathcal{L}_h u
	:=h^{-2}\sum_{k=-1}^1 \sum_{\ell=-1}^1
	C_{k,\ell} u(x_i+kh,y_j+\ell h)
	\quad \mbox{with}\quad
	C_{k,\ell}=\sum_{p=0}^{M+1} c_{k,\ell,p}(\ka h)^p,
	\]
	where  $c_{k,\ell,p} \in \R$ for all $k,\ell\in\{-1,0,1\}$. Furthermore, we let $C_{-1,-1}=C_{-1,1}=C_{1,-1}=C_{1,1}$ and $C_{-1,0}=C_{0,-1}=C_{0,1}=C_{1,0}$ for symmetry.	
	Approximating $u(x_i+kh,y_j+\ell h)$ as in \eqref{GQ:V} with $x_i^*=x_i$ and $y_j^*=y_j$, we have
	\[
	\begin{split}
	h^{-2}	\mathcal{L}_h u=h^{-2}\sum_{(m,n)\in \ind_{M+1}^{1}}
		u^{(m,n)} I_{m,n}+
		\sum_{(m,n)\in \ind_{	\tilde{M}-1}} f^{(m,n)}
		J_{m,n}=\mathcal{O}(h^{M}), \quad h \rightarrow 0,
	\end{split}
	\]
	where
	\be\label{Helm:regular:JVmn:add}
	I_{m,n}:=\sum_{k=-1}^1 \sum_{\ell=-1}^1 C_{k,\ell} G_{M+1,m,n}(kh, \ell h), \quad \text{and} \quad
	J_{m,n}:=h^{-2} \sum_{k=-1}^1 \sum_{\ell=-1}^1 C_{k,\ell} H_{	 \tilde{M}+1,m,n}(kh, \ell h).
	\ee
	Let
	\be\label{Helm:Lh:uh:regular:HM}
	\begin{split}
	h^{-2}	\mathcal{L}_h u_h:=h^{-2}\sum_{k=-1}^1 \sum_{\ell=-1}^1
		C_{k,\ell} (u_h)_{i+k,j+\ell}=
		\sum_{(m,n)\in \ind_{	\tilde{M}-1}} f^{(m,n)}J_{m,n}.
	\end{split}
	\ee
	Then
	\be\label{proof:orderM}
h^{-2}	\mathcal{L}_h (u-u_h)=
h^{-2}	\mathcal{L}_h u-\sum_{(m,n)\in \ind_{	 \tilde{M}-1}} f^{(m,n)}J_{m,n}
	=\bo(h^{M}),\qquad h\to 0,
	\ee
	if $I_{m,n}$ in \eqref{Helm:regular:JVmn:add} satisfies
	\begin{align}
		&I_{m,n}=\bo(h^{M+2}), \qquad h \to 0, \; \mbox{ for all }\; (m,n)\in \ind_{M+1}^{1}. \label{Helm:stencil:regular:u:V:eq}
	\end{align}
	By calculation, we find that $M=6$ is the maximum positive integer such that the linear system $\eqref{Helm:stencil:regular:u:V:eq}$ has a non-trivial solution.  All such non-trivial solutions for $M=6$ can be uniquely written (up to a constant multiple) as \eqref{CVmn:All}.
So \eqref{CVmn:All}, \eqref {Helm:Lh:uh:regular:HM}, \eqref{proof:orderM} and \eqref{Helm:stencil:regular:u:V:eq}  with $M=6$ and $\tilde{M}\ge 6$	complete the proof of	 \cref{thm:regular:all:Ckl}.
	\end{proof}
		
			\begin{proof}[Proof of \cref{thm:regular:interior}]
			
		Consider a general compact stencil $\{C^{\wa}_{k,\ell}\}_{k,\ell\in\{-1,0,1\}}$ parameterized by $C^{\wa}_{1,1}$, $C^{\wa}_{1,0} \in \R$ satisfying
		\[
		 C^{\wa}_{-1,-1}=C^{\wa}_{-1,1}=C^{\wa}_{1,-1}=C^{\wa}_{1,1}, \quad
		 C^{\wa}_{-1,0}=C^{\wa}_{0,-1}=C^{\wa}_{0,1}=C^{\wa}_{1,0}, \quad
		\text{and}
		\quad
		C^{\wa}_{0,0}=-20,
		\]
where we normalized the stencil by $C^{\wa}_{0,0}=-20$. Take a plane wave solution $u(x,y,\theta):=\exp({\ia}{\ka}(\cos(\theta)x+\sin(\theta)y))$ for any $\theta \in [0,2\pi)$. Clearly, we have $\Delta u+\ka^2u=0$. Hence, the truncation error, multiplied by $h^2$, associated with the general compact stencil coefficients $\{C^{\wa}_{k,\ell}\}_{k,\ell\in\{-1,0,1\}}$ at the grid point $(x_i,y_j)\notin \partial \Omega$ is
		\[
		(T(\theta|\ka h))_{x_i,y_j} := \sum_{k=-1}^1 \sum_{\ell=-1}^1 C^{\wa}_{k,\ell} \exp({\ia}{\ka}(\cos(\theta)(x_i+kh)+\sin(\theta)(y_j+\ell h))).
		\]
		Recall that $\frac{2\pi}{\ka h}$ is the number of points per wavelength. Hence, it is reasonable to choose $\ka h\in[1/4,1]$.
		Without loss of generality, we let $(x_i,y_j)=(0,0)$. Define $S:= \{\tfrac{1}{4}+\tfrac{3s}{4000}:s=0,\dots,1000\}$ and let
		 \be\label{Min:Int:Truncation:Special}
		(\tilde{C}^{\wa}_{1,1}(\ka h),\tilde{C}^{\wa}_{1,0}(\ka h)):=\argmin_{C^{\wa}_{1,1},C^{\wa}_{1,0} \in \R} \int_{0}^{2\pi} |(T(\theta|\ka h))_{0,0}|^2 d \theta, \quad \ka h \in S.
		\ee
		We use the Simpson's $3/8$ rule with 900 uniform sampling points to calculate $\int_{0}^{2\pi} |(T(\theta|\ka h))_{0,0}|^2 d \theta $.
		%
		Now, we link $C_{0,0},C_{1,0},C_{1,1}$ in \eqref{CVmn:All} with $C^{\wa}_{0,0},\tilde{C}^{\wa}_{1,0}(\ka h),\tilde{C}^{\wa}_{1,1}(\ka h)$ in \eqref{Min:Int:Truncation:Special} for $\ka h \in S$. To further simplify the presentation of our stencil coefficients, we set $c_{9}=c_{10}=c_{11}=0$ in \eqref{CVmn:All} so that the coefficients of the polynomials in \eqref{CVmn:All} for degree $7$ are zero.
		Because $C^{\wa}_{0,0}=-20$ is our normalization, we determine the free parameters $c_{i}$ for $i=1,\dots,8$ in \eqref{CVmn:All} by considering the following least-square problem:
			\[
			 (\tilde{c}_1,\tilde{c}_2,\dots,\tilde{c}_8):=	 
			\argmin_{{c_1,c_2,\dots,c_8 \in \R}} \sum_{\ka h \in S} |C_{1,1}(\ka h) + \tfrac{1}{20}\tilde{C}^{\wa}_{1,1}(\ka h)C_{0,0}(\ka h)|^2 +|C_{1,0}(\ka h)+\tfrac{1}{20}\tilde{C}^{\wa}_{1,0}(\ka h)C_{0,0}(\ka h)|^2.
			\]
		For simplicity of presentation, we replace each above calculated coefficient $\tilde{c}_{i}$ with its approximated fractional form $[2^{20}\tilde{c}_{i}]/2^{20}$, where $[\cdot]$ is a rounding operation to the nearest integer. Then we obtain \eqref{Value:C1:C11}. Plugging \eqref{Value:C1:C11} in \eqref{CVmn:All}, we obtain \eqref{stencil:Cv}.  Choosing $	 \tilde{M}=7$ in \eqref{stencil:all:free:Ckl} yields the right-hand side of \eqref{stencil:regular:interior:V}.
	\end{proof}

	\begin{proof}[Proof of \cref{thm:regular:Robin:1}]
	
		We only prove item (1). The proof of item (2) is very similar.
	Since $-u_x-\ia \ka u=g_1$ on $\BB_1:=\{l_{1}\} \times (l_3,l_4)$, we have $u^{(1,n)} = - \ia \ka  u^{(0,n)} - g_1^{(n)}$ for all $n = 0,\dots, 	 \tilde{M}-1$. By \eqref{GQ:V} with $M,	 \tilde{M}$ being replaced by $M-1,	 \tilde{M}-1$ and choosing $	\tilde{M}\ge M$, we have for $x,y \in (-2h,2h)$
	{\small{
			\begin{align*}
				 &u(x+x_i^*,y+y_j^*)\\
				&=	
				\sum_{(m,n)\in \ind_{M}^{1}}
				u^{(m,n)} G_{M,m,n}(x,y) +\sum_{(m,n)\in \ind_{	\tilde{M}-2}}
				f^{(m,n)} H_{	 \tilde{M},m,n}(x,y)+\bo(h^{M+1})\\
				&=
				\sum_{n=0}^{M}
				u^{(0,n)} G_{M,0,n}(x,y)+\sum_{n=0}^{M-1}
				u^{(1,n)} G_{M,1,n}(x,y) +\sum_{(m,n)\in \ind_{	\tilde{M}-2}}
				f^{(m,n)} H_{	 \tilde{M},m,n}(x,y) +\bo(h^{M+1}) \\
				& =
				\sum_{n=0}^{M}
				u^{(0,n)} G_{M,0,n}(x,y)+\sum_{n=0}^{	\tilde{M}-1}
				u^{(1,n)} G_{	 \tilde{M},1,n}(x,y) +\sum_{(m,n)\in \ind_{	 \tilde{M}-2}}
				f^{(m,n)} H_{	 \tilde{M},m,n}(x,y) +\bo(h^{M+1}) \\
				& = \sum_{n=0}^{M}
				u^{(0,n)} G_{M,0,n}(x,y)-\sum_{n=0}^{	\tilde{M}-1}
				\big(  \textsf{i} {\ka}  u^{(0,n)}+g_1^{(n)}  \big) G_{	 \tilde{M},1,n}(x,y) +\sum_{(m,n)\in \ind_{	 \tilde{M}-2}}
				f^{(m,n)} H_{	 \tilde{M},m,n}(x,y) +\bo(h^{M+1}) \\
				&  =u^{(0,M)} G_{M,0,M}(x,y)+\sum_{n=0}^{M-1}
				u^{(0,n)} \Big( G_{M,0,n}(x,y) -\textsf{i} {\ka}  G_{M,1,n}(x,y) \Big) -\sum_{n=0}^{	 \tilde{M}-1}
				g_{1}^{(n)}  G_{	 \tilde{M},1,n}(x,y)\\
				& \quad +\sum_{(m,n)\in \ind_{	\tilde{M}-2}}
				f^{(m,n)} H_{	 \tilde{M},m,n}(x,y)+\bo(h^{M+1}), \quad h \rightarrow 0.
			\end{align*}}}
	We set $C_{k,\ell}:=\sum_{p=0}^{M} (a_{k,\ell,p} + \ia b_{k,\ell,p})(\ka h)^p$, where  $a_{k,\ell,p},b_{k,\ell,p} \in \R$ for all $k\in\{0,1\}$ and $\ell \in \{-1,0,1\}$. Furthermore, we let $C_{0,-1}=C_{0,1}$ and $C_{1,-1}=C_{1,1}$ for symmetry.
	Letting $x_i^*=x_i$ and $y_j^*=y_j$ yields
	\[
	\begin{split}
	h^{-1}	\mathcal{L}_h u  :&= h^{-1}	 \sum_{k=0}^1 \sum_{\ell=-1}^1
		C_{k,\ell} u(x_i+kh,y_j+\ell h)\\
		&=h^{-1}\sum_{n=0}^{M}
		u^{(0,n)}  I_{n}+
		\sum_{(m,n)\in \ind_{	\tilde{M}-2}} f^{(m,n)}
		J_{m,n} +\sum_{n=0}^{	\tilde{M}-1} g_1^{(n)} J_{g_1,n}
		=\bo(h^{M}), \quad h\rightarrow 0,
	\end{split}
	\]
	where
	\be \label{Helm:IB1n:Proof}
	\begin{split}
		&I_{n}:=\sum_{k=0}^1 \sum_{\ell=-1}^1 C_{k,\ell} \left( G_{M,0,n}(kh, \ell h) -\textsf{i} {\ka}  G_{M,1,n}(kh, \ell h) (1-\delta_{n,M}) \right),\\
		& J_{m,n}:=h^{-1} \sum_{k=0}^1 \sum_{\ell=-1}^1 C_{k,\ell} H_{	 \tilde{M},m,n}(kh, \ell h), \quad
		J_{g_1,n}:=-h^{-1} \sum_{k=0}^1 \sum_{\ell=-1}^1 C_{k,\ell} G_{	 \tilde{M},1,n}(kh, \ell h),
	\end{split}
	\ee
	$\delta_{a,a}=1$, and $\delta_{a,b}=0$ for $a \neq b$. Let
	\be\label{Helm:IB1huh}
	\begin{split}
	h^{-1}	\mathcal{L}_h u_h := h^{-1}	 \sum_{k=0}^1 \sum_{\ell=-1}^1
		C_{k,\ell} (u_h)_{i+k,j+\ell}=
		\sum_{(m,n)\in \ind_{	\tilde{M}-2}} f^{(m,n)}
		J_{m,n} +\sum_{n=0}^{	\tilde{M}-1} g_1^{(n)}  J_{g_1,n}.
	\end{split}
	\ee
	We have
	$
	h^{-1}	\mathcal{L}_h (u-u_h)=\bo(h^{M})$, $h \rightarrow 0,$
	if $I_{n}$ for $n=0,\dots,M$ in \eqref{Helm:IB1n:Proof} satisfies
	\be\label{Helm:IB1hM+1}
	I_{n}=\bo(h^{M+1}), \quad h\rightarrow 0.
	\ee
	By calculation, we find that $M=6$ is the maximum positive integer such that the linear system of \eqref{Helm:IB1hM+1} has a non-trivial solution. To further simplify such a solution, we set coefficients associated with $\ka h$ of degrees higher than $4$ to zero; i.e., we now have polynomials of $\ka h$, whose highest degree is 4. All such non-trivial solutions for $M=6$ can be uniquely written (up to a constant multiple) as
	\[
			\begin{aligned}
				&C_{1,1} =1-60(c_{1}{\ia}+2c_{3}{\ia}+c_{4}{\ia}/2+c_{6}{\ia}-4{\ia}/225-c_{8}/2+c_{2}-c_{5}-2c_{7}){\ka h}+12(c_{8}{\ia}-7c_{2}{\ia}/3+7c_{5}{\ia}/3
				\\
				& \qquad
				 +13c_{7}{\ia}/3+7c_{1}/3+13c_{3}/3+c_{4} +7c_{6}/3-4/135)({\ka h})^2 +(c_{2}+c_{6}{\ia})({\ka h})^3 +(c_{3}+c_{7}{\ia})({\ka h})^4,
				\\
				& C_{0,1} = 2 -120(c_{1}{\ia}+2c_{3}{\ia}+c_{4}{\ia}/2+c_{6}{\ia}-29{\ia}/1800-c_{8}/2+c_{2}-c_{5}-2c_{7}){\ka h}+18(c_{8}{\ia}-22c_{2}{\ia}/9
				\\
				& \qquad
				 +22c_{5}{\ia}/9+40c_{7}{\ia}/9+22c_{1}/9+40c_{3}/9+c_{4}+22c_{6}/9
				-11/324)({\ka h})^2 +13(c_{1}{\ia}+20c_{3}{\ia}/13
				\\
				& \qquad
				 +7c_{4}{\ia}/26+12c_{6}{\ia}/13-17{\ia}/1170-7c_{8}/26+12c_{2}/13-c_{5}
				-20c_{7}/13)({\ka h})^3 +(c_{1}+c_{5}{\ia})({\ka h})^4,
				\\
				& C_{1,0} =4
				 -240(c_{1}{\ia}+2c_{3}{\ia}+c_{4}{\ia}/2+c_{6}{\ia}-29{\ia}/1800
				 -c_{8}/2+c_{2}-c_{5}-2c_{7}){\ka h}  	 +36(c_{8}{\ia}-22c_{2}{\ia}/9
				\\
				& \qquad
				 +22c_{5}{\ia}/9+40c_{7}{\ia}/9+22c_{1}/9+40c_{3}/9+c_{4}+22c_{6}/9
				-49/1620)({\ka h})^2  +18(c_{1}{\ia}+4c_{3}{\ia}/3
				\\
				& \qquad
				+c_{4}{\ia}/6 +8c_{6}{\ia}/9-{\ia}/90-c_{8}/6+8c_{2}/9-c_{5}
				-4c_{7}/3)({\ka h})^3  +(c_{4}+c_{8}{\ia})({\ka h})^4,
				\\
				& C_{0,0} =-10 +600(c_{1}{\ia}+2c_{3}{\ia}
                +c_{4}{\ia}/2+c_{6}{\ia}-29{\ia}/4500-c_{8}/2+c_{2}-c_{5}-2c_{7}){\ka h}
                +84(c_{8}{\ia}  -32c_{2}{\ia}/21
                \\
                & \qquad
                +32c_{5}{\ia}/21+74c_{7}{\ia}/21+32c_{1}/21+74c_{3}/21+c_{4}+32c_{6}/21+1/3780)({\ka h})^2
                -80(c_{1}{\ia}+2c_{3}{\ia}
                \\
                & \qquad +c_{4}{\ia}/2+39c_{6}{\ia}/40-7{\ia}/720-c_{8}/2+39c_{2}/40-c_{5}-2c_{7})({\ka h})^3
                -4(c_{8}{\ia}-3c_{2}{\ia}/2+2c_{5}{\ia}+7c_{7}{\ia}/2
                \\
                & \qquad +2c_{1}+7c_{3}/2+c_{4}+3c_{6}/2-1/80)({\ka h})^4,
			\end{aligned}
	\]
	where each $c_i \in \R$ for $i=1,\dots,8$ is a free parameter.  Choosing $	 \tilde{M}=8$ in \eqref{Helm:IB1n:Proof} and \eqref{Helm:IB1huh} yields the right-hand side of \eqref{stencil:regular:interior:V:Robin:1}.

		Next, consider a 6-point stencil $\{C^{\wa}_{k,\ell}\}_{k\in\{0,1\},\ell\in\{-1,0,1\}}$ parameterized by $C^{\wa}_{1,1}, C^{\wa}_{0,1}, C^{\wa}_{1,0} \in \C$ with
		\[
		C^{\wa}_{1,-1}=C^{\wa}_{1,1}, \quad C^{\wa}_{0,-1}=C^{\wa}_{0,1}, \quad \text{and} \quad C^{\wa}_{0,0}=-10,
		\]
where we normalized the general stencil by $C^{\wa}_{0,0}=-10$.
		Take a plane wave solution $u(x,y,\theta):=\exp({\ia}{\ka}(\cos(\theta)x+\sin(\theta)y))$ for any $\theta \in [0,2\pi)$. Clearly, we have $\Delta u+\ka^2u=0$ and $- u_x -\ia \ka u=g_1 \neq 0$ on $\BB_1:=\{l_{1}\} \times (l_3,l_4)$, where $g_1$ and its derivatives are explicitly known by plugging the plane wave solution $u(x,y,\theta)$ into the boundary condition. Hence, the truncation error, multiplied by $h$, associated with the compact general stencil coefficients $\{C^{\wa}_{k,\ell}\}_{k\in\{0,1\},\ell\in\{-1,0,1\}}$ at the grid point $(x_0,y_j) \in \BB_1:=\{l_{1}\} \times (l_3,l_4)$ is
		\[
		\begin{split}
		(T(\theta|\ka h))_{x_0,y_j} :=& \sum_{k=0}^1 \sum_{\ell=-1}^1 C^{\wa}_{k,\ell} \exp({\ia}{\ka}(\cos(\theta)(x_0 + kh)+\sin(\theta)(y_j + \ell h)))\\
		&+\sum_{n=0}^{7}g_{1}^{(n)}  \sum_{k=0}^1 \sum_{\ell=-1}^1 C^{\wa}_{k,\ell} G_{8,1,n}(kh, \ell h).
		\end{split}
		\]
		Without loss of generality, we let $(x_0,y_j)=(0,0)$. Afterwards, we follow a similar minimization procedure as in the proof of \cref{thm:regular:interior} to obtain the concrete stencils in \cref{thm:regular:Robin:1}.
	\end{proof}
	
	\begin{proof}[Proof of \cref{thm:corner:1}]

	By $\B_1u:=\frac{\partial u}{\partial \nv}- \ia \ka u=g_1$ and $\B_3u:=\frac{\partial u}{\partial \nv}=g_3$, we have
		\be\label{g1g3:derivs}
		u^{(1,n)} = - \ia \ka u^{(0,n)} - g_1^{(n)} \quad \text{and} \quad
		u^{(m,1)} = - g_3^{(m)}, \quad \text{for all } m,n \in \NN.
		\ee
		Let $C_{k,\ell}:=c_{k,\ell} + \tilde{c}_{k,\ell}$ for $k,\ell \in \{0,1\}$, where $c_{k,\ell}$ and $ \tilde{c}_{k,\ell}$ are to be determined polynomials of $\ka h$. Note that $x_i^*=x_0$ and $y_j^*=y_0$.
		Approximating $u(x_0 + kh, y_0 + \ell h)$ by \eqref{GQ:V}, \eqref{GQ:H} with $M,	 \tilde{M}$ being replaced by $M-1,	 \tilde{M}-1$,  and using \eqref{g1g3:derivs}, we have
		\begin{align}
		\label{stencil:corner:1:IJK}
		& h^{-1}	\mathcal{L}_h u  := h^{-1} \sum_{k=0}^{1} \sum_{\ell=0}^{1} ( c_{k,\ell} +  \tilde{c}_{k,\ell}) u(x_0 + kh, y_0 + \ell h) = 	h^{-1} \sum_{n=0}^{M} u^{(0,n)} I_{n}\\
		\nonumber
		& \quad + h^{-1} \sum_{m=0}^{M} u^{(m,0)} \tilde{I}_{m}  + \sum_{(m,n) \in \ind_{	 \tilde{M}-2}} f^{(m,n)} S_{m,n} + \sum_{n=0}^{	 \tilde{M}-1} g_1^{(n)} K_{n}  + \sum_{m=0}^{	\tilde{M}-1} g_3^{(m)} \tilde{K}_{m} +\bo(h^{M}), \quad h\rightarrow 0,
		\end{align}
		where
				\begin{align*}
					& \tilde{I}_{m} := \sum_{k=0}^{1} \sum_{\ell=0}^{1} \tilde{c}_{k,\ell} G_{M,0,m}( \ell h, kh),\\
					& I_{n} := \sum_{k=0}^{1} \sum_{\ell=0}^{1} c_{k,\ell}\left( G_{M,0,n}(kh, \ell h) -\textsf{i} {\ka} G_{M,1,n}(kh, \ell h) (1-\delta_{n,M}) \right),\\
					&S_{m,n} :=	h^{-1} \sum_{k=0}^{1} \sum_{\ell=0}^{1} (c_{k,\ell} H_{	 \tilde{M},m,n}(kh,\ell h) + \tilde{c}_{k,\ell} H_{	\tilde{M},n,m}(\ell h, kh)),\\
					& K_{n} := -	h^{-1} \sum_{k=0}^{1} \sum_{\ell=0}^{1} c_{k,\ell} G_{	 \tilde{M},1,n}(kh,\ell h), \quad \text{and} \quad \tilde{K}_{m} :=  -	 h^{-1} \sum_{k=0}^{1} \sum_{\ell=0}^{1} \tilde{c}_{k,\ell} G_{	\tilde{M},1,m}(\ell h, kh).
				\end{align*}
			Let
		\be \label{stencil:corner:1:proof}
\begin{aligned}
	h^{-1}	\mathcal{L}_h u_h  :=
	\begin{aligned}	
		h^{-1}
		\sum_{k=0}^1 \sum_{\ell=0}^1 C_{k,\ell}(u_{h})_{k,\ell}
	\end{aligned}
	= \sum_{(m,n)\in \ind_{\tilde{M}-2}} f^{(m,n)}J_{m,n} + \sum_{n=0}^{\tilde{M}-1}g_{1}^{(n)}J_{g_{1},n} + \sum_{n=0}^{\tilde{M}-1}g_{3}^{(n)}J_{g_{3},n},
\end{aligned}	
\ee
where
\be \label{corner:1:righthand}
\begin{aligned}
	& J_{g_1,2\ell} := K_{2\ell} + \sum_{p=\max\{\ell,1\}}^{\left\lfloor \frac{	 \tilde{M}-1}{2} \right\rfloor} (-1)^{p+1} \binom{p}{\ell} \ka^{2(p-\ell)} \frac{\tilde{I}_{2p+1}}{h} - \frac{\tilde{I}_{1} \delta_{\ell,0}}{h}, \quad \ell=0, \dots, \left\lfloor \tfrac{	\tilde{M}-1}{2} \right\rfloor,\\
	& J_{g_1,2\ell+1}:=  K_{2\ell+1}, \quad \ell=0, \dots, \left\lfloor \tfrac{	 \tilde{M}-2}{2} \right\rfloor,\quad  \text{and} \quad J_{g_3,\ell} := \tilde{K}_{\ell}, \quad \ell=0, \dots, 	 \tilde{M}-1,\\
	& J_{\ell,2j+1} := S_{\ell,2j+1}, \quad \ell=0, \dots, 	\tilde{M}-2j-3, j=0,\dots,\left\lfloor \tfrac{	 \tilde{M}-1}{2} -1\right\rfloor, \\
	& J_{2\ell + \gamma ,2j} := \sum_{p=\max\{j+\ell+1,1\}}^{\left\lfloor \frac{	 \tilde{M}-\gamma}{2} \right\rfloor} (-1)^{p-\ell-1} \binom{p-\ell-1}{j} \ka^{2(p-\ell-j-1)} \frac{\tilde{I}_{2p+\gamma}}{h} + S_{2\ell + \gamma,2j},
\end{aligned}
\ee
$\gamma \in \{0,1\}$, $j=0,\dots,\left\lfloor \tfrac{	 \tilde{M}-\gamma}{2} \right\rfloor-\ell-1$, and $\ell=0, \dots, \left\lfloor \tfrac{\tilde{M}-\gamma}{2} \right\rfloor-1$.
		By replacing $u^{(m,0)}$ for $m=2,\dots,M$ with \eqref{umnV2:relation}, using \eqref{g1g3:derivs}, and rearranging some terms, \eqref{stencil:corner:1:IJK} and \eqref{stencil:corner:1:proof} imply
		{\small
			\begin{align}
				\nonumber
				& \frac{\mathcal{L}_h (u-u_h)}{h} = \frac{u^{(0,0)}}{h} \left(I_{0} + \tilde{I}_{0} - \ia \ka \tilde{I}_{1} + \sum_{p=1}^{\lfloor{\frac{M}{2}}\rfloor}(-1)^p \ka^{2p} \tilde{I}_{2p} + \ia \sum_{p=1}^{\lfloor{\frac{M-1}{2}}\rfloor}(-1)^{p+1} \ka^{2p+1}  \tilde{I}_{2p+1}\right) + \sum_{\ell=0}^{\lfloor{\frac{M-1}{2}}\rfloor} u^{(0,2\ell+1)} \frac{I_{2\ell+1}}{h} \\
				\nonumber
				& \quad + \sum_{\ell=1}^{\lfloor{\frac{M-1}{2}}\rfloor} \frac{u^{(0,2\ell)}}{h} \left( \sum_{p=\max\{\ell,1\}}^{\lfloor{\frac{M}{2}}\rfloor} (-1)^p \binom{p}{\ell} \ka^{2(p-\ell)} \tilde{I}_{2p} + \ia \sum_{p=\max\{\ell,1\}}^{\lfloor{\frac{M-1}{2}}\rfloor} (-1)^{p+1} \binom{p}{\ell} \ka^{2(p-\ell)+1} \tilde{I}_{2p+1} + I_{2\ell}\right)\\
				\nonumber
				& \quad + \frac{u^{(0,2\lfloor \frac{M}{2} \rfloor)}}{h} \left( (-1)^{\lfloor\frac{M}{2}\rfloor} \tilde{I}_{2\lfloor \frac{M}{2} \rfloor} + I_{2\lfloor \frac{M}{2} \rfloor}\right) \left(1-\delta_{\lfloor\frac{M}{2}\rfloor,\lfloor\frac{M-1}{2}\rfloor}\right) + \sum_{\ell=0}^{\lfloor\frac{	 \tilde{M}-2}{2}\rfloor} g_1^{(2\ell+1)} (K_{2\ell+1}- J_{g_1,2\ell+1})\\
				\nonumber
				& \quad + \sum_{\ell=0}^{\lfloor\frac{	 \tilde{M}-1}{2}\rfloor}g_1^{(2\ell)} \left(K_{2\ell} + \sum_{p=\max\{\ell,1\}}^{\lfloor\frac{	 \tilde{M}-1}{2}\rfloor} (-1)^{p+1} \binom{p}{\ell} \ka^{2(p-\ell)} \frac{\tilde{I}_{2p+1}}{h} - \frac{\tilde{I}_{1} \delta_{\ell,0}}{h} - J_{g_1,2\ell}\right) + \sum_{\ell=0}^{	\tilde{M}-1} g_3^{(\ell)} (\tilde{K}_{\ell}-J_{g_3,\ell})\\
				\nonumber
				& \quad  + \sum_{j=0}^{\lfloor\frac{	 \tilde{M}-1}{2}-1\rfloor}\sum_{\ell=0}^{	 \tilde{M}-2j-3} f^{(\ell,2j+1)}(S_{\ell,2j+1} - J_{\ell,2j+1}) + \sum_{\gamma \in \{0,1\}} \sum_{\ell=0}^{\lfloor\frac{	\tilde{M}- \gamma }{2}\rfloor-1} \sum_{j=0}^{\lfloor\frac{	 \tilde{M}-\gamma}{2}\rfloor-\ell-1} f^{(2\ell+\gamma,2j)}\left(\sum_{p=\max\{j+\ell+1,1\}}^{\lfloor\frac{	 \tilde{M}-\gamma}{2}\rfloor} \right.\\
				\nonumber
				& \quad \left. (-1)^{p-\ell-1} \binom{p-\ell-1}{j} \ka^{2(p-\ell-j-1)} \frac{\tilde{I}_{2p+\gamma}}{h} + S_{2\ell+\gamma,2j} - J_{2\ell+\gamma,2j}\right) = \bo (h^{M}), \quad h\rightarrow 0.
			\end{align}
		}
		We set $c_{k,\ell}=\sum_{j=0}^{M} (a_{k,\ell,j} + \ia b_{k,\ell,j})(\ka h)^j$ and $\tilde{c}_{k,\ell}=\sum_{j=0}^{M} (\tilde{a}_{k,\ell,j} + \ia \tilde{b}_{k,\ell,j})(\ka h)^j$, where  $a_{k,\ell,j}$, $b_{k,\ell,j}$, $\tilde{a}_{k,\ell,j}$, $\tilde{b}_{k,\ell,j} \in \R$ for all $k,\ell \in\{0,1\}$. By calculation, $M=6$ is the maximum positive integer such that the linear system, obtained by setting each coefficient of $u^{(0,n)}$ for $n=0,\dots,6$ to be $\bo(h^{7})$ as $h \rightarrow 0$, has a non-trivial solution. Afterwards, to further simplify such a solution, we can set remaining coefficients associated with $(\ka h)^5$ or $(\ka h)^6$ to zero.
		
		By using the minimization procedure described in the proofs of \cref{thm:regular:interior,thm:regular:Robin:1}, we can verify that $c_{0,1}=c_{1,1}=\tilde{c}_{0,0}=\tilde{c}_{1,0}=0$, $c_{0,0}=C_{0,0}$,  $c_{1,0}=C_{1,0}$,  $\tilde{c}_{0,1}=C_{0,1}$, and $\tilde{c}_{1,1}=C_{1,1}$,
		where $\{C_{k,\ell}\}_{k,\ell \in \{0,1\}}$ are defined in \eqref{CR1:Corner:1}. Given these $\{c_{k,\ell}\}_{k,\ell \in \{0,1\}}$ and $\{\tilde{c}_{k,\ell}\}_{k,\ell \in \{0,1\}}$, we set $	\tilde{M}=8$ and plug them into the relations in \eqref{corner:1:righthand}. This completes the proof of \cref{thm:corner:1}.
	\end{proof}
	
	\begin{proof}[Proof of \cref{thm:corner:2}]
		The proof is almost identical to the proof of \cref{thm:corner:1}. Note that we need to replace $u^{(m,1)}=-g_3^{(m)}$ with $u^{(m,1)} = \ia \ka u^{(m,0)} + g_4^{(m)}$ for all $m\in \NN$ in \eqref{g1g3:derivs}.
	\end{proof}
	
	For the following theorems, we note that an identity similar to \eqref{GQ:V} still holds: for $x,y\in (-2h,2h)$,
	\begin{align}
		\label{u:approx:ir:key:2}
		u_\pm (x+x_i^*,y+y_j^*)
		& =\sum_{(m,n)\in \ind_{M}^{1}}
		u_\pm^{(m,n)} G^{\pm}_{M,m,n}(x,y) +\sum_{(m,n)\in \ind_{	\tilde{M}-2}}
		f_\pm ^{(m,n)} H^{\pm}_{	 \tilde{M},m,n}(x,y)+\bo(h^{M+1}),
	\end{align}
	as $h \rightarrow 0$, where $\tilde{M}\ge M$, $\ind_{M}^{1}$ is defined in \eqref{indV12}, $\ind_{	\tilde{M}-2}$ is defined in \eqref{Sk}, $G^{\pm}_{M,m,n}(x,y)$ is obtained by replacing $\ka$ by $\ka_{\pm}$ and $M+1$ by $M$ in \eqref{Hel:GM1mnxy}, $H^{\pm}_{	 \tilde{M},m,n}(x,y)$ is obtained by replacing $\ka$ by $\ka_{\pm}$ and $M+1$ by $\tilde{M}$ in \eqref{Hel:HM1mnxy}.
	
	\begin{proof}[Proof of \cref{thm:interface}]
	The proof closely follows from the proof of \cite[Theorem 2.3]{FHM21a}.
\end{proof}

\begin{proof}[Proof of \cref{fluxtm2}] 
For an irregular point $(x_i,y_j)$, we define
\[
h^{-1}\mathcal{L}_h u:=h^{-1} \sum_{k=-1}^1 \sum_{\ell=-1}^1
C_{k,\ell}  u(x_i+kh,y_j+\ell h).
\]
By \eqref{base:pt:gamma} and \eqref{u:approx:ir:key:2},
we have
		\be\label{ELLIP:u:sum:s1s2}
		\begin{split}
			h^{-1}\mathcal{L}_h u
			&=h^{-1}\sum_{(k,\ell)\in d_{i,j}^+}
			C_{k,\ell}    u(x_i^*+(v_0+k)h,y_j^*+(w_0+\ell) h) \\
			&\qquad +h^{-1} \sum_{(k,\ell)\in d_{i,j}^-}
			C_{k,\ell}  u(x_i^*+(v_0+k)h,y_j^*+(w_0+\ell) h)\\
			&=h^{-1} \sum_{(m,n)\in \ind_{M}^1}
			u_+^{(m,n)}  I^{+}_{m,n}+ \sum_{(m,n)\in \ind_{\tilde{M}-2}}
			f_+ ^{(m,n)} J^{+,0}_{m,n}+h^{-1}\sum_{(m,n)\in \ind_{M}^1}
			u_-^{(m,n)}  I^{-}_{m,n}  \\
			&\qquad + \sum_{(m,n)\in \ind_{\tilde{M}-2}}
			f_- ^{(m,n)} J^{-,0}_{m,n}+\bo(h^{M}),\quad h \rightarrow 0,\\
		\end{split}
		\ee
where
\begin{align} 
	\nonumber
	&C_{k,\ell}:=\sum_{p=0}^{M} c_{k,\ell,p}(\max(\ka_+,\ka_-) h)^p,\; c_{k,\ell,p} \in \R, \qquad I^\pm_{m,n}:=\sum_{(k,\ell)\in d_{i,j}^\pm}
	C_{k,\ell} G^{\pm}_{M,m,n}((v_0+k)h,(w_0+\ell) h),\\
	\label{ELLIP:IJpm}
	&J^{\pm,0}_{m,n}:=h^{-1}\sum_{(k,\ell)\in d_{i,j}^\pm} C_{k,\ell}  H^{\pm}_{\tilde{M},m,n}((v_0+k)h, (w_0+\ell) h).
\end{align}
Using \eqref{transmission:Hel}, we obtain
\begin{align*}
	h^{-1} \sum_{(m',n')\in \ind_{M}^1}
	u_-^{(m',n')}   I^-_{m',n'}
	=&h^{-1} \sum_{(m,n)\in \ind_{M}^1} u_+^{(m,n)}
	J^{u_+,T}_{m,n}+
	 \sum_{(m,n)\in \ind_{M-2}}
	f_+^{(m,n)} J^{+,T}_{m,n}\\
	&+\sum_{(m,n)\in \ind_{M-2}} f_-^{(m,n)} J^{-,T}_{m,n}+\sum_{p=0}^{M} g^{(p)}   J^{g}_p+\sum_{p=0}^{M-1} g^{(p)}_{\Gamma}   J^{g_\Gamma}_p,
\end{align*}
where
\be
\label{ELLIP:JTno}
\begin{split}
	&J^{u_+,T}_{m,n}:=
	\sum_{(m',n')\in \ind_{M}^1} I^-_{m',n'} T^{u_+}_{m',n',m,n}, \quad J^{\pm,T}_{m,n}:=h^{-1}
	\sum_{(m',n')\in \ind_{M}^1} I^-_{m',n'} T^\pm_{m',n',m,n},\\
	&J^{g}_{p}:=h^{-1}
	\sum_{(m',n')\in \ind_{M}^1} I^-_{m',n'} T^{g}_{m',n',p},\quad
	J^{g_\Gamma}_{p}:=h^{-1}
	\sum_{(m',n')\in \ind_{M}^1} I^-_{m',n'} T^{g_\Gamma}_{m',n',p},
\end{split}
\ee
$T^{u_+}_{m',n',m,n}$, $T^{\pm}_{m',n',m,n}$, $T^{g}_{m',n',p}$, $T^{g_{\Gamma}}_{m',n',p}$ are transmission coefficients in \eqref{transmission:Hel}.
Let
	\be \label{ELLIP:stencil:irregular}
	\begin{split}
	h^{-1}	\mathcal{L}_h u_h & :=h^{-1}\sum_{k=-1}^1 \sum_{\ell=-1}^1
		C_{k,\ell} (u_h)_{i+k,j+\ell}\\
		&=\sum_{(m,n)\in \ind_{\tilde{M}-2}} \left(f_-^{(m,n)}J^-_{m,n}
		+f_+^{(m,n)}J^+_{m,n}\right)+ \sum_{p=0}^{\tilde{M}} g^{(p)} J^{g}_p+\sum_{p=0}^{\tilde{M}-1} g^{(p)}_{\Gamma} J^{g_\Gamma}_p,
	\end{split}
	\ee
where  $\tilde{M}\ge M$ and $J^{\pm}_{m,n}:= J_{m,n}^{\pm,0}+J^{\pm,T}_{m,n}$.
%
%
Hence, $h^{-1} \mathcal{L}_h (u-u_h)=\bo(h^M)$, $h\to 0$, if the following holds
%
%
\be \label{ELLIP:stencil:regular:u:ir}
I_{m,n}:=I^+_{m,n}+J^{u_+,T}_{m,n}=\bo(h^{M+1}), \quad h\to 0, \; \mbox{ for all }\; (m,n)\in \ind_{M}^1.
\ee
Since $\ka_{+}=\ka_{-}$, by \eqref{transmission:Hel:equal}, \eqref{ELLIP:JTno}, and the definition of $I_{m,n}$,  we have $J^{u_+,T}_{m,n}= I^-_{m,n}$, and so, $I_{m,n} =I^+_{m,n}+I^-_{m,n}$.
%
%
%
Similar to the existence of the nontrivial
solution $\{C_{k,\ell}\}_{k,\ell=-1,0,1}$ for \eqref{Helm:stencil:regular:u:V:eq} and the proof of \cref{thm:regular:interior},
we can say
the largest $M$ such that the nontrivial
solution $\{C_{k,\ell}\}_{k,\ell=-1,0,1}$ exists for \eqref{ELLIP:stencil:regular:u:ir} with  $\ka_{+}=\ka_{-}$ is $M=7$. The coefficients $\{C_{k,\ell}\}_{k,\ell=-1,0,1}$ in \eqref{stencil:Cv} yields the left-hand side of \eqref{irregular:same:k}. Letting $M=\tilde{M}=7$ in \eqref{ELLIP:IJpm}--\eqref{ELLIP:stencil:irregular} 
yields the right-hand side of \eqref{irregular:same:k} and \eqref{irregular:same:k:right}.
\end{proof}	
\begin{proof}[Proof of \cref{fluxtm3}]
	Similar to the proof of \cref{fluxtm2}, $\{C_{k,\ell}\}_{k,\ell=-1,0,1}$ in \eqref{irregular:diff:k} are obtained by solving
\be \label{System:Ckl}
\begin{split}
&\sum_{(k,\ell)\in d_{i,j}^+}
 C_{k,\ell} G^{+}_{M,m,n}((v_0+k)h,(w_0+\ell) h)\\
&+	\sum_{(m',n')\in \ind_{M}^1} \sum_{(k,\ell)\in d_{i,j}^-}
C_{k,\ell} G^{-}_{M,m',n'}((v_0+k)h,(w_0+\ell) h) T^{u_+}_{m',n',m,n} =\bo(h^{M+1}), \quad h \rightarrow 0,
\end{split}
\ee
 for all $(m,n)\in \ind_{M}^1$, where $ T^{u_+}_{m',n',m,n}$ are the transmission coefficients in \eqref{transmission:Hel} and\\ $C_{k,\ell}:=\sum_{p=0}^{M} c_{k,\ell,p}(\max(\ka_+,\ka_-) h)^p$, $c_{k,\ell,p} \in \R$.
 By calculation, $M=5$ is the largest positive integer such that the linear system of \eqref{System:Ckl} has a non-trivial solution.
Letting $M=\tilde{M}=5$ in \eqref{ELLIP:IJpm}--\eqref{ELLIP:stencil:irregular} yields the right-hand side of \eqref{irregular:diff:k} and \eqref{irregular:diff:k:right}.
\end{proof}

\end{document}